\documentclass[final,3p,times,fleqn]{elsarticle}
\usepackage{mathtools,booktabs}
\usepackage{bm}
\usepackage{bbm}
\usepackage{makecell}
\usepackage{multirow}
\usepackage{amsmath}
\usepackage{amsthm}
\usepackage{gensymb} 
\usepackage{footnote} 
\allowdisplaybreaks[4] 
\usepackage{epstopdf}
\usepackage[T1]{fontenc}
\usepackage{color} 
\usepackage{graphicx}
\usepackage{epstopdf}
\usepackage{subfig} 
\usepackage{subfloat}
\usepackage{hyperref}
\usepackage{stmaryrd}
\usepackage{rotating}

\usepackage{threeparttable}

\biboptions{sort&compress}
 
\begin{document}
	\begin{frontmatter}
		\title{A unified fourth-order Bhatnagar-Gross-Krook lattice Boltzmann model for high-dimensional linear hyperbolic equations}
\author[a]{Ying Chen}  
\author[a,b,c]{Zhenhua Chai \corref{cor1}}
\ead{hustczh@hust.edu.cn}
\author[a,b,c]{Baochang Shi} 
\address[a]{School of Mathematics and Statistics, Huazhong University of Science and Technology, Wuhan, 430074, China}
\address[b]{Institute of Interdisciplinary Research for Mathematics and Applied Science, Huazhong University of Science and Technology, Wuhan 430074, China}
\address[c]{Hubei Key Laboratory of Engineering Modeling and
	Scientific Computing, Huazhong University of Science and Technology,
	Wuhan 430074, China}
\cortext[cor1]{Corresponding author.}
\begin{abstract}
	In this work, we first develop a unified Bhatnagar-Gross-Krook lattice Boltzmann (BGK-LB) model for the $d$($d\geq 1$)-dimensional linear hyperbolic equation (L-HE), where the natural moments and the D$d$Q$(2d^2+1)$ [($2d^2+1$) discrete velocities in $d$-dimensional space] lattice structure are considered. Subsequently, at the acoustic scaling, we conduct an accuracy analysis on the developed BGK-LB model by the direct Taylor expansion (DTE) method, and present the second- and third-order moments of the equilibrium distribution functions (EDFs) to ensure that the BGK-LB model can be fourth-order consistent with the L-HE. And on this basis, when considering the  Dirichlet boundary condition, the fourth-order full-way and half-way boundary schemes are proposed to approximate the unknown distribution functions to ensure that the BGK-LB model can be overall fourth-order accurate. Thereafter, based on the kinetic entropy theory, we derive the conditions that the fourth-order moments of the EDFs should satisfy to ensure the microscopic entropy stability of the BGK-LB model. In addition, with the aid of the von Neumann stability analysis, we also discuss the $L^2$ stability of the BGK-LB model and numerically plot the  stability regions. In particular, from a numerical perspective, we find that the region of microscopic entropy stability is identical to that of $L^2$ stability. Finally, we carry out some numerical experiments to test the accuracy and stability of the BGK-LB model, and the numerical results are in agreement with our theoretical analysis. In addition, we compare the developed full-way and half-way boundary schemes for the Dirichlet boundary condition, which shows that the full-way boundary scheme is more stable.	
\end{abstract}
\end{frontmatter}
\section{Introduction}
Over the past three decades, the mesoscopic lattice Boltzmann (LB) method--a highly efficient second-order kinetic theory-based approach for fluid flow problems governed by the Navier-Stokes equations \cite{Guo2013,Kruger2017,Wang2019}--has been successfully extended to the study of solving macroscopic partial differential equations (PDEs), including diffusion equations \cite{HUBER20107956,ANCONA1994107,Suga,Lin2022,SILVA2023105735,Chen2023}, convection-diffusion equations \cite{VANDERSMAN2000766,GINZBURG20051171,RASIN2005453,PhysRevE.79.016701,Chai2013,Chai2016,Chen2023-CDE}, and Burgers equations \cite{Bruce,Ma2005,Elton1996,Velivelli2006,Duan2006,Zhang2008,Li2018,Chen2024-Burger}. In particular, for the collision terms in the LB method, the LB models can be classified into three basic types: the Bhatnagar-Gross-Krook (BGK) or single-relaxation-time (SRT) model \cite{YHQian1992}, the two-relaxation-time (TRT) model \cite{GINZBURG20051171}, and the multiple-relaxation-time (MRT) model \cite{Chai2020}. However, although the TRT and MRT models have the potential for greater accuracy and/or stability by using additional adjustable parameters, the optimal selection approach for these adjustable parameters remains unclear. Thus, compared to the TRT-LB and MRT-LB models, the SRT-LB model is the most efficient and universal. For this reason, we consider the BGK-LB model in this work.

It is widely known to us that through some asymptotic analysis approaches, including the Chapman-Enskog analysis \cite{chapman1990}, Maxwell iteration \cite{ikenberry1956,Yong2016}, direct Taylor expansion (DTE) \cite{HOLDYCH2004,Wagner2006}, recurrence equations \cite{d2009viscosity,Ginzburg2012} and equivalent equations \cite{DUBOIS2008,Dubois2009,Dubois2022}, it is straightforward to validate that the LB method has a second-order accuracy.   However, there is relatively little research to study the higher-order consistency of the LB model with the PDE to be solved.  In fact, in order to effectively control the overall error in numerical simulations and thus to more accurately describe the physical information in PDEs, developing the higher-order LB model  is necessary.  In recent years, several researchers have made significant contributions to the development of higher-order LB models for some specific PDEs. Chen et. al. \cite{Chen2023-CDE} developed a fourth-order MRT-LB model with the D1Q3 lattice structure for the one-dimensional convection-diffusion equation. Suga \cite{Suga} focused on the one-dimensional diffusion equation and proposed a fourth-order BGK-LB model with the D1Q3 lattice structure. Subsequently, Lin et. al. \cite{Lin2022} further developed a sixth-order MRT-LB model. Silva \cite{SILVA2023105735} considered the TRT collision model and obtained a fourth-order LB model for the one-dimensional equation with a linear source term. For the two-dimensional diffusion equation, Chen et. al.  \cite{Chen2023} proposed a fourth-order MRT-LB model with the D2Q5 lattice structure. It is noteworthy that the aforementioned high-order LB models are restricted to the one- and two-dimensional cases. In a recent study, for the $d$($d\geq 1$)-dimensional coupled Burgers’ equations  based on the Cole-Hopf transformation \cite{cole1951}, Chen et. al. \cite{Chen2024-Burger} proposed a unified fourth-order MRT-LB model with the natural moments and the D$d$Q($2d^2+1$) lattice structure  employed. 

From the works mentioned above, one can find that the PDEs that can be solved with high-order accuracy by the LB model are some special cases of the general convection-diffusion equations. In particular, when we consider that the diffusion coefficient is equal to zero, the general convection-diffusion equation  reduces to the hyperbolic equation. However, the high-order LB model for the hyperbolic equation has not been extensively investigated. This is largely because it is enough to study the LB model for the general convection-diffusion equation with a non-zero diffusion term at the diffusion scaling, while the LB model for the hyperbolic equation must be analyzed at the acoustic scaling. Under this premise, the second-order accuracy of the LB model for the hyperbolic equation can be achieved only when the relaxation parameter related to the diffusion coefficient is equal to 2 \cite{GRAILLE2014,ANANDAN2024,guillon2024}. This challenges our conventional understanding of the LB method, as we always consider that the instability will occur when the relaxation parameter approaches  2 \cite{YOSHIDA2010,Chai2014}. Nonetheless, recently, Bellotti et. al. \cite{bellotti2023-weakly} considered the one-dimensional scalar linear hyperbolic equation (L-HE) and proposed a fourth-order BGK-LB model with the relaxation parameter equal to 2, followed by an analysis of the $L^2$ stability. What is more, they also developed a fourth-order entropy-stable LB model for the multidimensional  nonlinear hyperbolic equation \cite{bellotti2024-entropy}. In this work, a time-splitting method is employed in the LB model, which means that the developed LB model in Ref. \cite{bellotti2024-entropy} requires some additional time steps compared to the standard LB model. In fact, developing a high-order LB model without any additional treatments to the standard LB model is extremely difficult, especially for the nonlinear problems. 

In this paper, we will only consider the L-HE and further focus on higher-dimensional cases. Similar to the previous work in Ref. \cite{Chen2024-Burger}, we will adopt the natural moments and the D$d$Q($2d^2 + 1$) lattice structure to develop a unified standard BGK-LB model for the high-dimensional L-HE. On this basis, when developing a higher-order LB model, we would like to point out that the fundamental idea of constructing higher-order LB models for the given PDEs in the above works \cite{Chen2023-CDE,Suga,Lin2022,SILVA2023105735,Chen2023,Chen2024-Burger} is as follows. 

 -- Give the general form of the EDFs (e.g., the commonly used linear or quadratic form \cite{Chai2016}) with the weight coefficients being free.

-- Adopt the weight coefficients and the free relaxation parameters to eliminate the higher-order truncation errors of the macroscopic equation recovered from the LB model when it is compared with the PDE to be solved. 

Unlike the works mentioned in Refs.  \cite{Chen2023-CDE,Suga,Lin2022,SILVA2023105735,Chen2023,Chen2024-Burger}, in this paper, we will consider the moments of the EDFs as free parameters to develop a unified higher-order standard BGK-LB model for the L-HE instead of defining the form of the EDFs in advance. In particular, the form of the EDFs determined from the moments in our developed high-order LB model will be entirely different from the commonly used EDFs.

The paper is organized as follows. In Sec. \ref{Sec-unified-BGK-LB}, we introduce the L-HE addressed in this paper and propose the corresponding  unified BGK-LB model.  In Sec. \ref{Sec-fourth-unified-BGK-LB}, we perform an accuracy analysis on the BGK-LB model and determine the high-order moments of the  EDFs to ensure that the BGK-LB model can be fourth-order consistent with the L-HE. In Sec. \ref{Sec-stability-BGK-LB}, the microscopic entropy stability  and  $L^2$  stability of the BGK-LB model are analyzed. In Sec. \ref{Numer}, some numerical experiments are performed to confirm the theoretical analysis of the developed fourth-order BGK-LB model. Finally, some conclusions are presented in Sec. \ref{Conclusion}.
\section{The unified Bhatnagar-Gross-Krook lattice Boltzmann model for $d$-dimensional linear  hyperbolic equations}\label{Sec-unified-BGK-LB}
In this section, we will present a unified BGK-LB model for the $d$-dimensional L-HE with the natural moments and the D$d$Q($2d+1$) lattice structure considered.  
\subsection{The $d$-dimensional linear  hyperbolic equation}
We aim at approximating the solution of the following L-HE:
\begin{align}\label{HE}
\begin{cases}
	\partial_t\phi(\mathbf{x},t)+\mathbf{u}\cdot\nabla\phi(\mathbf{x},t)=R,&\:\:\:\:\:\:\:\:\:\mathbf{x}\in \mathbb{R}^d,\:\: t\in(0,T],\\
	\phi(\mathbf{x},0)=\phi^{0}(\mathbf{x}),&\:\:\:\:\:\:\:\:\:\:\:\:\:\:\:\:\:\:\:\:\mathbf{x}\in \Omega \subset \mathbb{R}^d,
\end{cases}
\end{align}
where the scalar variable $\phi$ and the initial datum $\phi^0$ are assumed to be given smooth functions in this work. $\mathbf{u}\in \mathbb{R}^d$ is the transport velocity, and $R$ is the constant source term. $\Omega$ is the computational domain with the straight boundaries considered in this work. In particular, when the source term $R\equiv0$, Eq. (\ref{HE}) follows the hyperbolic conservation law, and the quadratic entropy and corresponding entropy flux can be given by $S(\phi)=\phi^2/2$ and $G(\phi)=\phi^2\mathbf{u}/2$, respectively.
\subsection{Spatial and temporal discretization}
To simplify the subsequent analysis, we utilize a uniform Cartesian lattice $\mathcal{L}_{\mathbf{x}}=\Delta x\mathbb{Z}^d$ with a lattice spacing of $\Delta x > 0$ to discretize the $d$-dimensional space in the LB method. Here, for simplicity, the more general rectangular lattice structure \cite{Chai2023-changfangxing} is not taken into consideration. The time is uniformly discretized as $t_n=n\Delta t\in\mathcal{T}$ (where $\mathcal{T}=\Delta t\mathbb{N}$ is the temporal lattice and $n\in\mathbb{N}$), and $\Delta t$ is the time step. Subsequently, the so-called lattice velocity is defined as $c=\Delta x/\Delta t$. In particular, for the L-HE (\ref{HE}), we consider the acoustic scaling, that is, $\Delta x$ is proportional to $\Delta t$ ($\Delta x\propto\Delta t$).
\subsection{The Bhatnagar-Gross-Krook lattice Boltzmann model}
Here, we would like to point out that although the more general MRT-LB model would have good accuracy and/or stability in the study of complex problems \cite{Chai2013,Chai2016,chai2016comparative}, so far only the concept of microscopic entropy stability of the BGK-LB model has been proposed \cite{GRAILLE2014,ANANDAN2024,guillon2024,bellotti2024-entropy}. Therefore, in order to develop an entropy-stable LB model for the L-HE (\ref{HE}), we consider the more efficient BGK-LB model. The corresponding evolution equation is written as
\begin{align}\label{BGK-LB-Evolution}
f_k(\mathbf{x}+\mathbf{c}_k\Delta t,t+\Delta t)=f_k(\mathbf{x},t)-s\big[f_k(\mathbf{x},t)-f_k^{eq}(\mathbf{x},t)\big]+\Delta tR_k,\:\: k\in \llbracket1,q\rrbracket \footnotemark[1],
\end{align}  where $f_k(\mathbf{x},t)$, $f_k^{eq}(\mathbf{x},t)$, and $R_k$  represent the distribution functions, EDFs, and discrete source terms at position $\mathbf{x}$ and time $t$, respectively. $q$ indicates the number of discrete velocities in the D$d$Q$q$ lattice structure. In the following, we will provide a unified lattice structure with $q=2d^2+1$ for the $d$-dimensional L-HE (\ref{HE}) [see Eqs. (\ref{d1q3}), (\ref{d2q9}), (\ref{d3q19}), and (\ref{ddqq}) below]. $s$ is the relaxation parameter, due to the fact that there is no  diffusion term in the L-HE (\ref{HE}), we take $s=2$  in Eq. (\ref{BGK-LB-Evolution}). Then, Eq. (\ref{BGK-LB-Evolution}) becomes
\footnotetext[1]{We shall consistently use the notation $\llbracket{1,q}\rrbracket:=\{1,2,\ldots,q\}$.}
\begin{align}\label{BGK-LB-Evolution-seq2}
f_k(\mathbf{x}+\mathbf{c}_k\Delta t,t+\Delta t)=-f_k(\mathbf{x},t)+2f_k^{eq}(\mathbf{x},t)+\Delta tR_k,\:\: k\in \llbracket{1,2d^2+1}\rrbracket.
\end{align}
In the above Eq. (\ref{BGK-LB-Evolution-seq2}), the D$d$Q($2d^2+1$) lattice structure with the discrete velocity $\mathbf{c}=\big(\mathbf{c}_k\big)_{k=1}^{2d^2+1}$ and the transform matrix $\mathbf{M}$ based on the natural moments are given by 
\begin{subequations}\label{DdQq-M}
	\begin{align}
	&	d=1: \begin{cases}\label{d1q3}
			\mathbf{c}^{x_1}=\big(0,1,-1\big)^Tc,&\\
			\mathbf{c}=\mathbf{c}^{x_1} ,&\\ 
			\mathbf{M} =\Big(I_{3},\mathbf{c}^{x_1},(\mathbf{c}^{x_1})^{.2})\Big)^T, 
		\end{cases} \\
	 		&	d=2:	\begin{cases}\label{d2q9}
				\mathbf{c}^{x_1} =\big(0,1,0,-1,0,1,-1,-1,1\big)^Tc,&\\
				\mathbf{c}^{x_2} =\big(0,0,1,0,-1,1,1,-1,-1\big)^Tc,& \\
				\mathbf{c}=\big(\mathbf{c}^{x_1},\mathbf{c}^{x_2}\big)^T,&\\
				\mathbf{M}  =\Big(I_{ 9},\mathbf{c}^{x_1},\mathbf{c}^{x_2},(\mathbf{c}^{x_1})^{.2},(\mathbf{c}^{x_2})^{.2},\mathbf{c}^{x_1}.\mathbf{c}^{x_2},(\mathbf{c}^{x_1})^{.2}\mathbf{c}^{x_2},\mathbf{c}^{x_1}.(\mathbf{c}^{x_2})^{.2},(\mathbf{c}^{x_1})^{.2}.(\mathbf{c}^{x_2})^{.2}\Big)^T, 
			\end{cases}  \\
			&d=3:	\begin{cases}\label{d3q19}
				\mathbf{c}^{x_1} =\big(0,1,0,0,-1,0,0,1,-1,-1,1,1,-1,-1,1,0,0,0,0\big)^Tc,&\\
				\mathbf{c}^{x_2} =\big(0,0,1,0,0,-1,0,1,1,-1,-1,0,0,0,0,1,-1,-1,1\big)^Tc,&\\
				\mathbf{c}^{x_3} =\big(0,0,0,1,0,0,-1,0,0,0,0,1,1,-1,-1,1,1,-1,-1\big)c,&\\
				\mathbf{c}=\big(\mathbf{c}^{x_1},\mathbf{c}^{x_2},\mathbf{c}^{x_3}\big)^T,&\\
				\mathbf{M} =\Big(I_{19},\mathbf{c}^{x_1},\mathbf{c}^{x_2},\mathbf{c}^{x_3},&\\
				\qquad\quad (\mathbf{c}^{x_1})^{.2},(\mathbf{c}^{x_2})^{.2},(\mathbf{c}^{x_3})^{.2},\mathbf{c}^{x_1}.\mathbf{c}^{x_2},\mathbf{c}^{x_1}.\mathbf{c}^{x_3},\mathbf{c}^{x_2}.\mathbf{c}^{x_3},\\
				\qquad\quad(\mathbf{c}^{x_1})^{.2}.\mathbf{c}^{x_2},(\mathbf{c}^{x_1})^{.2}.\mathbf{c}^{x_3},
				\mathbf{c}^{x_1}.(\mathbf{c}^{x_2})^{.2},(\mathbf{c}^{x_2})^{.2}.\mathbf{c}^{x_3},
				\mathbf{c}^{x_1}.(\mathbf{c}^{x_3})^{.2},	\mathbf{c}^{x_2}.(\mathbf{c}^{x_3})^{.2},&\\
				\qquad\quad (\mathbf{c}^{x_1})^{.2}.(\mathbf{c}^{x_2})^{.2},
				(\mathbf{c}^{x_1})^{.2}.(\mathbf{c}^{x_3})^{.2},
				(\mathbf{c}^{x_2})^{.2}.(\mathbf{c}^{x_3})^{.2}\Big)^T,
			\end{cases} \\
			&d>3:	\begin{cases} \label{ddqq}
				\mathbf{c}^{x_1} =\big(Q_1,\underbrace{J_1^{d-1}}_{4(d-1)},O_{\eta_1}\big)^Tc,  \footnotemark[2]&\\
				\mathbf{c}^{x_2} =\big(Q_2,\underbrace{J_2,O_{4(d-2)}}_{4(d-1)},J_1^{d-2},O_{\eta_2}\big)^Tc,\\
				\qquad\qquad\qquad\qquad\qquad\vdots&\\  
				\mathbf{c}^{x_k} =\big(Q_k,\underbrace{O_{4(k-2)},J_2,O_{4(d-k)}}_{4(d-1)},\underbrace{O_{4(k-3)},J_2,O_{4(d-k)}}_{4(d-2)},\ldots,\underbrace{O_{4(k-l-1)},J_2,O_{4(d-k)}}_{4(d-l)},\ldots,J_1^{d-k},O_{\eta_k}\bigg)^Tc,&\\
				\qquad\qquad\qquad\qquad\qquad\vdots&\\  
				\mathbf{c}^{x_d} =\bigg(Q_{d},\underbrace{O_{4(d-2)},J_2}_{4(d-1)},\underbrace{O_{4(d-3)},J_2}_{4(d-2)},\ldots,\underbrace{O_{4(d-j-1)]},J_2}_{4(d-j)},\ldots,\underbrace{O_{4},J_2}_{4\times 2},J_2\bigg)^Tc,&\\
				\mathbf{c}=\big(\mathbf{c}^{x_1},\mathbf{c}^{x_2},\ldots,\mathbf{c}^{x_d}\big)^T,&\\
				\mathbf{M} =\Big(I_{1+2d^2},\mathbf{c}_{x_1}^T,\mathbf{c}_{x_2}^T,\ldots,\mathbf{c}_{x_d}^T,(\mathbf{c}_{x_1}^{.2})^T,(\mathbf{c}_{x_2}^{.2})^T,\ldots,(\mathbf{c}_{x_d}^{.2})^T,&\\
				\qquad(\mathbf{c}_{x_1}.\mathbf{c}_{x_2})^T,(\mathbf{c}_{x_1}.\mathbf{c}_{x_3})^T,\ldots,(\mathbf{c}_{x_1}.\mathbf{c}_{x_d})^T,\ldots,(\mathbf{c}_{x_{i-1}}.\mathbf{c}_{x_d})^T,(\mathbf{c}_{x_i}.\mathbf{c}_{x_{i+1}})^T,\ldots,(\mathbf{c}_{x_i}.\mathbf{c}_{x_d})^T,\ldots,(\mathbf{c}_{x_{d-1}}.\mathbf{c}_{x_d})^T,&\\
				\quad\quad(\mathbf{c}_{x_1}^{.2}.\mathbf{c}_{x_2})^T,(\mathbf{c}_{x_1}^{.2}.\mathbf{c}_{x_3})^T,\ldots,(\mathbf{c}_{x_1}^{.2}.\mathbf{c}_{x_d})^T,\ldots,&\\
				\quad\quad(\mathbf{c}_{x_i}^{.2}.\mathbf{c}_{x_1})^T,(\mathbf{c}_{x_i}^{.2}.\mathbf{c}_{x_2})^T,\ldots,(\mathbf{c}_{x_i}^{.2}.\mathbf{c}_{x_{i-1}})^T,(\mathbf{c}_{x_i}^{.2}.\mathbf{c}_{x_{i+1}})^T,\ldots,(\mathbf{c}_{x_i}^{.2}.\mathbf{c}_{x_d})^T,\ldots,&\\
				\quad\quad(\mathbf{c}_{x_d}^{.2}.\mathbf{c}_{x_1})^T,(\mathbf{c}_{x_d}^{.2}.\mathbf{c}_{x_2})^T,\ldots,(\mathbf{c}_{x_d}^{.2}.\mathbf{c}_{x_{d-1}})^T,&\\
				\quad\quad(\mathbf{c}_{x_1}^{.2}.\mathbf{c}_{x_2}^{.2})^T,(\mathbf{c}_{x_1}^{.2}.\mathbf{c}_{x_3}^{.2})^T,\ldots,(\mathbf{c}_{x_1}^{.2}.\mathbf{c}_{x_d}^{.2})^T,\ldots,(\mathbf{c}_{x_{i-1}}^{.2}.\mathbf{c}_{x_d}^{.2})^T,(\mathbf{c}_{x_i}^{.2}.\mathbf{c}_{x_{i+1}}^{.2})^T,\ldots,(\mathbf{c}_{x_i}^{.2}.\mathbf{c}_{x_d}^{.2})^T,\ldots,(\mathbf{c}_{x_{d-1}}^{.2}.\mathbf{c}_{x_d}^{.2})^T\Big), 
			\end{cases} 
	\end{align}
	\end{subequations}
where $ k\in\llbracket 1,d\rrbracket$, $l\in\llbracket1,k-1\rrbracket$, $j\in\llbracket1,d-1\rrbracket$ and
\footnotetext[2]{The notation $\underbrace{J}_{p}$ indicates that the number of elements in the vector $J$ is equal to $p$.}
\begin{subequations}
\begin{align}
	&I_{p}=\big(1,1,\cdots,1)^T\in \mathbb{R}^{p\times 1}, \\
	&(\mathbf{c}^{x_1})^{.p_1}.(\mathbf{c}^{x_2})^{.p_1}.\cdots.(\mathbf{c}^{x_d})^{.p_d}=\bigg(\prod_{i=1}^d (c^{x_i}_1)^{p_i},\prod_{i=1}^d(c^{x_i}_2)^{p_i},\ldots,\prod_{i=1}^d (c^{x_i}_{2d^2+1})^{p_i}\bigg),\:\: p_i\in \llbracket{0,2}\rrbracket,\:\:i\in \llbracket{1,d}\rrbracket,\\
	&Q_p=(0,0,\ldots,\underbrace{1}_{{(p+1)}_{th}},0,\ldots,0,\ldots,\underbrace{-1}_{(d+p+1)_{th}},0,\ldots,0)\in \mathbb{R}^{(2d+1)\times 1},\\
	&J_1^p=(\underbrace{J_1,J_1,\ldots,J_1}_{4p})\in \mathbb{R}^{(4p)\times1},\:\: p\in \llbracket{1,d-1}\rrbracket,\\ &J_1=(1,-1,-1,1),\:\:J_2=(1,1,-1,-1),\\
	&O_p=\big(0,0,\ldots,0)^T\in \mathbb{R}^{p\times 1},\\ 
	&\eta_i=2d(d-1)-2i(2d-1-i),		
\end{align}
\end{subequations} 
here $c^{x_k}_i$ represents the $i_{th}$ element of the column vector $\mathbf{c}^{x_k}$. In order to correctly recover the L-HE (\ref{HE}) from the BGK-LB model (\ref{BGK-LB-Evolution-seq2}), the zeroth- and first-order moments of the EDFs $f_k^{eq}$, and the zeroth-order moment of the discrete source terms $R_k$, should satisfy the following conditions:   
\begin{subequations}\label{0-1-meq}
\begin{align}
	&m_1^{eq}=\sum_{k=1}^{2d^2+1}f_k^{eq}=\phi=\sum_{k=1}^{2d^2+1}f_i,\\
	&m_l^{eq}=\sum_{k=1}^{2d^2+1}c_{kl}f_k^{eq}=\phi u_{k-1},\:\: l\in \llbracket{2,d+1}\rrbracket,\\
	&m^R_1=R,
\end{align}
\end{subequations} 
where moments $m_k^{eq}$ and $m^R_k$ ($ k\in \llbracket{1,2d^2+1}\rrbracket$) are the $k_{th}$ elements of $\mathbf{m}^{eq}=\mathbf{Mf}^{eq}$ and $\mathbf{m}^R=\mathbf{M}\mathbf{R}$ [$\mathbf{f}^{eq}:=(f_k^{eq})_{k=1}^{2d^2+1}$ and $\mathbf{R}:=(R_k)_{k=1}^{2d^2+1}$], respectively. In this work, we assume that $\mathbf{m}^{eq}=\bm{\varepsilon}\phi$ and $\mathbf{m}^R=\bm{\varepsilon}R$ with   
\begin{align}\label{epsilon}
	\bm{\varepsilon}=\bigg[1,\mathbf{u},\mathbf{\dot{m}}^{eq|2}_{\mathbf{x}^{.2}},\mathbf{\dot{m}}^{eq|2}_{(\mathbf{xx})_{\alpha<\beta}},\mathbf{\dot{m}}^{eq|3}_{\mathbf{x}^{.2}\mathbf{\overline{x}}},\mathbf{\dot{m}}^{eq|4}_{(\mathbf{x}^{.2}\mathbf{x}^{.2})_{\alpha<\beta}}\bigg]^T,
\end{align}
here 
\begin{subequations}\label{epsilon-elements}
	\begin{align} 
	&\mathbf{\dot{m}}^{eq|2}_{\mathbf{x}^{.2}}=\Big(\dot{m}^{eq|2}_{x_1^2},\dot{m}^{eq|2}_{x_2^2},\ldots,\dot{m}^{eq|2}_{x_d^2}\Big),\\
	&\mathbf{\dot{m}}^{eq|2}_{(\mathbf{xx})_{\alpha<\beta}}=\Big(\dot{m}^{eq|2}_{x_1x_2},\dot{m}^{eq|2}_{x_1x_3},\ldots,\dot{m}^{eq|2}_{x_1x_d},\dot{m}^{eq|2}_{x_2x_3},\dot{m}^{eq|2}_{x_2x_4},\ldots,\dot{m}^{eq|2}_{x_2x_d},\ldots,\dot{m}^{eq|2}_{x_ix_{i+1}},\ldots,\dot{m}^{eq|2}_{x_ix_d},\ldots,\dot{m}^{eq|2}_{x_{d-1}x_d}\Big),\\
&\mathbf{\dot{m}}^{eq|3}_{\mathbf{x}^{.2}\mathbf{\overline{x}}}=\Big(\dot{m}^{eq|3}_{x_1^2x_2},\dot{m}^{eq|3}_{x_1^2x_3},\ldots,\dot{m}^{eq|3}_{x_1^2x_d},\dot{m}^{eq|3}_{x_2^2x_1},\dot{m}^{eq|3}_{x_2^2x_3},\ldots,\dot{m}^{eq|3}_{x_2^2x_d},\ldots,\dot{m}^{eq|3}_{x_i^2x_1},\ldots,\dot{m}^{eq|3}_{x_i^2x_{i-1}},\dot{m}^{eq|3}_{x_i^2x_{i+1}},\ldots,\dot{m}^{eq|3}_{x_i^2x_d},\ldots,\dot{m}^{eq|3}_{x_d^2x_{d-1}}\Big),\\
&\mathbf{\dot{m}}^{eq|4}_{(\mathbf{x}^{.2}\mathbf{x}^{.2})_{\alpha<\beta}}=\Big(\dot{m}^{eq|4}_{x_1^2x_2^2},\dot{m}^{eq|4}_{x_1^2x_3^2},\ldots,\dot{m}^{eq|4}_{x_1^2x_d^2},\dot{m}^{eq|4}_{x_2^2x_3^2},\dot{m}^{eq|4}_{x_2^2x_4^2},\ldots,\dot{m}^{eq|4}_{x_2^2x_d^2},\ldots,\dot{m}^{eq|4}_{x_i^2x_{i+1}^2},\ldots,\dot{m}^{eq|4}_{x_i^2x_d^2},\ldots,\dot{m}^{eq|4}_{x_{d-1}^2x_d^2}\Big),
	\end{align}
\end{subequations}
where the element  $\dot{m}^{eq|r}_{x_i^{r_i}x_j^{r_j}}$ for each $r_{i,j}\in\mathbb{N}$ with $r_i+r_j=r\in\llbracket2,4\rrbracket$, which are related to the transport velocity  $\mathbf{u}$ and lattice velocity $c$, corresponds to one of the $r_{th}$-order moments of the EDFs $f_k^{eq}$ and is given by. 
  \begin{align}
 	\dot{m}^{eq|r}_{x_i^{r_i}x_j^{r_j}}=\Big((\mathbf{c}^{x_i})^{.r_i}.(\mathbf{c}^{x_j})^{.r_j}\Big)^T\bm{\omega},
 \end{align}
here the column vector $\bm{\omega}=\mathbf{M}^{-1}\bm{\varepsilon}$  satisfies $\mathbf{f}^{eq}=\bm{\omega}\phi$.

According to Eq. (\ref{0-1-meq}), it is easy to validate that the BGK-LB model (\ref{BGK-LB-Evolution-seq2}) can recover the L-HE (\ref{HE})  with a second-order accuracy through some frequently-used asymptotic methods \cite{chapman1990,ikenberry1956,Yong2016,HOLDYCH2004,Wagner2006,d2009viscosity,Ginzburg2012,DUBOIS2008,Dubois2009,Dubois2022}. Meanwhile, the second- to fourth-order moments of the EDFs $f_k^{eq}$ remain unconstrained and can be adjusted to change certain features of the BGK-LB model (\ref{BGK-LB-Evolution-seq2}), particularly in terms of accuracy and stability. (See Secs. \ref{Sec-fourth-unified-BGK-LB} and \ref{Sec-stability-BGK-LB} for details.)
\section{The accuracy analysis of the unified Bhatnagar-Gross-Krook lattice Boltzmann model}\label{Sec-fourth-unified-BGK-LB}
 In this section, a detailed accuracy analysis will be carried out on the BGK-LB model (\ref{BGK-LB-Evolution-seq2}) using the DTE method. Subsequently, the conditions of the second- and third-orders moments of the EDFs $f_k^{eq}$ will be presented, which ensure that at the acoustic scaling, the BGK-LB model (\ref{BGK-LB-Evolution-seq2}) can be consistent with the L-HE (\ref{HE}) up to the order of $O(\Delta t^4)$. Then, to guarantee that the overall accuracy of the developed BGK-LB model (\ref{BGK-LB-Evolution-seq2}) can reach fourth-order, we will discuss how to approximate the unknown distribution functions $f_k$ for the initialization and Dirichlet boundary conditions of the L-HE (\ref{HE}).

\subsection{The direct Taylor expansion} \label{dte-analysis}
In this part, we will utilize the DTE method \cite{HOLDYCH2004,Wagner2006,Chai2020} to derive Eq. (\ref{HE}) from the above BGK-LB model (\ref{BGK-LB-Evolution-seq2}) with a fourth-order truncation error. To this end, we first apply the Taylor expansion to Eq. (\ref{BGK-LB-Evolution-seq2}) up to the order of $O(\Delta t^5)$, where the acoustic scaling is considered. To simplify the following derivation, we denote $\psi=\psi(\mathbf{x},t)$, $(\mathbf{x},t)\in(\mathcal{L}_{\mathbf{x}},\mathcal{T})$ ($\psi = f_k$, $f^{eq}_k$, and $R_k$). And one can get 
\begin{align}\label{taylor-expansion-eq}
	&f_k+\Delta t\mathbf{c}_k\cdot \nabla f_k+\frac{\Delta t^2}{2}\mathbf{c}_k^{.2}\overset{2}{\cdot} \nabla^2 f_k+\frac{\Delta t^3}{6}\mathbf{c}_k^{.3}\overset{3}{\cdot} \nabla^3 f_k
	+\frac{\Delta t^4}{24}\mathbf{c}_k^{.4}\overset{4}{\cdot}\nabla^4 f_k\notag\\
	&\quad=-f_k
	+\Delta t\partial_tf_k
	-\frac{\Delta t^2}{2}\partial_t^2f_k
	+\frac{\Delta t^3}{6}\partial_t^3f_k-\frac{\Delta t^4}{24}\partial_t^4f_k\notag\\
	&\qquad+2f_k^{eq}-2\Delta t\partial_tf_k^{eq}+\Delta t^2\partial_t^2f_k^{eq}-\frac{\Delta t^3}{3}\partial_t^3f_k^{eq}
	+\frac{\Delta t^4}{12}\partial_t^4f_k^{eq}+ \Delta t  R_k+O(\Delta t^5),
\end{align} 
where   $\mathbf{c}_k\overset{i}{\cdot} \nabla :=\sum_{\alpha_i=1}^{d}\cdots\sum_{\alpha_1=1}^{d}\big(c_{k\alpha_1}\cdots c_{k\alpha_i}\nabla^i_{\alpha_1\ldots\alpha_i}\big)$ with $i\in \llbracket2,4\rrbracket$. Based on the relation of $f_k=f^{eq}_k+f_k^{ne}$, from Eq. (\ref{taylor-expansion-eq}), one can obtain the following first- to fourth-order truncation equations of the distribution functions $f_k$:
 \begin{subequations}\label{trun-or-1-4}
	\begin{align}
		&f_k^{ne}=O(\Delta t),\label{ME-1} \\
		&f_k=f_k^{eq}-\frac{\Delta t}{2}\mathbf{c}_k\cdot\nabla f_k^{eq}-\frac{\Delta t}{2}\partial_tf_k^{eq}+\frac{\Delta t}{2}R_k+O(\Delta t^2),\label{ME-2} \\
		&f_k=f_k^{eq}-\Delta t\partial_tf_k^{eq}+\frac{\Delta t^2}{4}\partial_t^2f_k^{eq}-\frac{\Delta t}{2}\mathbf{c}_k\cdot\nabla f_k-\frac{\Delta t^2}{4}\mathbf{c}_k^{.2}\overset{2}{\cdot}\nabla^2f_k^{eq}+\frac{\Delta t}{2}\partial_tf_k +\frac{\Delta t}{2}R_k+O(\Delta t^3), \label{ME-3}\\
		&f_k=f_k^{eq}-\Delta t\partial_tf_k^{eq}+\frac{\Delta t^2}{2}\partial_t^2f_k^{eq}-\frac{\Delta t^2}{4}\partial_t^2f_k-\frac{\Delta t}{2}\mathbf{c}_k\cdot\nabla f_k-\frac{\Delta t^2}{4}\mathbf{c}_k^{.2}\overset{2}{\cdot}\nabla^2f_k^{eq}+\frac{\Delta t}{2}\partial_tf_k\notag\\
		&\qquad-\frac{\Delta t^3}{12}\partial_t^3f_k^{eq}-\frac{\Delta t^3}{12}\mathbf{c}_k^{.3}\overset{3}{\cdot} \nabla^3 f_k+\frac{\Delta t}{2}R_k +O(\Delta t^4)\label{ME-4}.
	\end{align}
\end{subequations}  
Substituting Eq. (\ref{ME-2}) into Eq. (\ref{taylor-expansion-eq}) and  taking the zeroth-order moment, one can obtain the following modified equation which is second-order consistent with the L-HE (\ref{HE}): 
\begin{align}\label{eq-ode2-simp}
	\partial_t\phi+\mathbf{u}\cdot\nabla\phi=R+O(\Delta t^2).
\end{align}	
To analyze the higher-order truncation errors of the BGK-LB model (\ref{BGK-LB-Evolution-seq2}) in comparison with the L-HE (\ref{HE}), we further consider the high-order expansions of the distribution functions $f_k$ as given in Eqs.  (\ref{ME-3}) and (\ref{ME-4}). In particular, according to Eq. (\ref{ME-2}), Eq. (\ref{ME-3}) can be rewritten as
\begin{align}\label{eq-od3}
	f_k=&f_k^{eq}-\Delta t\partial_tf_k^{eq}+\frac{\Delta t^2}{4}\partial_t^2f_k^{eq}-\frac{\Delta t}{2}\mathbf{c}_k\cdot\nabla \Big(f_k^{eq}-\frac{\Delta t}{2}\mathbf{c}_k\cdot\nabla f_k^{eq}-\frac{\Delta t}{2}\partial_tf_k^{eq}\Big)\notag\\
	&-\frac{\Delta t^2}{4}\mathbf{c}_k^{.2}\overset{2}{\cdot}\nabla^2f_k+\frac{\Delta t}{2}\Big(\partial_tf_k^{eq}-\frac{\Delta t}{2}\mathbf{c}_k\cdot\nabla \partial_tf_k^{eq}-\frac{\Delta t}{2}\partial_t^2f_k^{eq}\Big)+O(\Delta t^3)\notag\\
	&=f_k^{eq}-\frac{\Delta t}{2}\mathbf{c}_k\cdot\nabla f_k^{eq}-\frac{\Delta t}{2}\partial_tf_k^{eq}+\frac{\Delta t}{2}R_k+O(\Delta t^3),
\end{align}	
by comparing Eq. (\ref{eq-od3}) with Eq. (\ref{ME-2}), it can be observed that the second- and third-order expansions of the distribution functions $f_k$ are identical when the second- and third-order terms, i.e., $O(\Delta t^2)$ in Eq. (\ref{ME-2}) and $O(\Delta t^3)$ in Eq. (\ref{eq-od3}), are neglected.  Substituting the above Eq. (\ref{eq-od3}) into Eq. (\ref{taylor-expansion-eq}) yields
\begin{align}\label{eq-ode3-simp}
	&2f_k+\Delta t\mathbf{c}_k\cdot\nabla \Big(f_k^{eq}-\frac{\Delta t}{2}\mathbf{c}_k\cdot\nabla f_k^{eq}-\frac{\Delta t}{2}\partial_tf_k^{eq}+\frac{\Delta t}{2}R_k\Big)+\frac{\Delta t^2}{2}\mathbf{c}_k^{.2}\overset{2}{\cdot}\nabla^2 \Big(f_k^{eq}-\frac{\Delta t}{2}\mathbf{c}_k\cdot\nabla f_k^{eq}-\frac{\Delta t}{2}\partial_tf_k^{eq}+\frac{\Delta t}{2}R_k\Big)\notag\\
	&\quad+\frac{\Delta t^3}{6}\mathbf{c}_k^{.3}\overset{3}{\cdot} \nabla^3 f_k^{eq}+O(\Delta t^4)\notag\\
	&=\Delta t\partial_tf_k
	-\frac{\Delta t^2}{2}\partial_t^2f_k
	+\frac{\Delta t^3}{6}\partial_t^3f_k +2f_k^{eq}-2\Delta t\partial_tf_k^{eq}+\Delta t^2\partial_t^2f_k^{eq}-\frac{\Delta t^3}{3}\partial_t^3f_k^{eq}+\Delta tR_k+O(\Delta t^4),
\end{align}
then taking the zeroth-order moment of Eq. (\ref{eq-ode3-simp}), one can obtain the following third-order modified equation compared to the L-HE (\ref{HE}): 
 	 \begin{align} 
 	 	&\nabla\cdot\mathbf{u}\phi-\frac{\Delta t}{2}\mathbf{u}\cdot\nabla\phi_t
 		-\frac{\Delta t^2}{4}\mathbf{q}^2\overset{2}{\cdot}\nabla^2\phi_t-\frac{\Delta t^2}{12}\mathbf{q}^3\overset{3}{\cdot}\nabla^3\phi  =-\partial_t\phi
 		+\frac{\Delta t}{2}\partial_t^2\phi
 		-\frac{\Delta t^2}{6}\partial_t^3\phi+R+O(\Delta t^3), 
 	\end{align} 
where the symmetric tensors $\mathbf{q}^r$ ($r=2,3$) are given by
 \begin{subequations}
\begin{align}
	&d=1:\:\mathbf{q}^2=\dot{m}^{eq}_{x_1^2},\:\:	\nabla^2=\partial_{x_1}^2,\:\:
	\mathbf{q}^3=\dot{m}^{eq}_{x_1^3},\:\: \nabla^3=\partial_{x_1}^3,\\
	&d=2:\:\mathbf{q}^2=\left(
	\begin{matrix}
		\dot{m}^{eq}_{x_1^2}&\dot{m}^{eq}_{x_1x_2}\\
		\dot{m}^{eq}_{x_1x_2}&\dot{m}^{eq}_{x_2^2}
	\end{matrix}
	\right),\:\:
	\nabla^2=
	\left(
	\begin{matrix}
		\partial_{x_1}^2&\partial_{x_1x_2}^2\\
		\partial_{x_1x_2}^2&\partial_{x_2}^2
	\end{matrix}
	\right),\:\:
	\mathbf{q}^3=\left(
	\begin{matrix}
		\dot{m}^{eq}_{x_1^3}&3\dot{m}^{eq}_{x_1^2x_2}\\
		3\dot{m}^{eq}_{x_1x_2^2}&\dot{m}^{eq}_{x_2^3}
	\end{matrix}
	\right),\:\:	\nabla^3=
	\left(
	\begin{matrix}
		\partial_{x_1}^3&\partial_{x_1^2x_2}^3\\
		\partial_{x_1x_2^2}^3&\partial_{x_2}^3
	\end{matrix}
	\right),\\
	&d\geq 3:\:\big(\mathbf{q}^r\big)_{i_1i_2\ldots i_r}=\dot{m}^{eq|r}_{x_{1}^{n_1}x_{2}^{n_2}\ldots x_{d}^{n_d}},\:\:	\big(\nabla^r\big)_{i_1i_2\ldots i_r}=\partial^r_{x_{1}^{n_1}x_{2}^{n_2}\ldots x_{d}^{n_d}},\:\: \sum_{l=1}^dn_l=r,
\end{align}
 \end{subequations}
here, for  each $l\in\llbracket1,d\rrbracket$, $n_l$ represents the number  of the elements  in the vector $(i_1,i_2,\ldots,i_r)$ that are equal to $l$.  

At the acoustic scaling, with the aid of the free parameter, i.e., the second- to fourth-order moments of the EDFs $f_k^{eq}$, one can obtain a unified BGK-LB model (\ref{BGK-LB-Evolution-seq2}) for the L-HE (\ref{HE}) with third-order accuracy  once  the following condition is satisfied: 
\begin{align} \label{third-condition}
	&
	-3\mathbf{q}^2\overset{2}{\cdot}\nabla^2\big(\mathbf{u}\cdot\nabla\phi\big)+\mathbf{q}^3\overset{3}{\cdot}\nabla^3\phi +2\mathbf{u}^3\overset{3}{\cdot}\phi=0, 
\end{align} 
where Eq. (\ref{eq-ode2-simp}) has been used. After some algebraic manipulations, the tensors $\mathbf{q^2}$ and $\mathbf{q}^3$ that are related to the second- and third-order moments of the EDFs $f_k^{eq}$, respectively, i.e., the solutions of the above condition (\ref{third-condition}), can be determined as
\begin{subequations}\label{third-condition-2}
\begin{align} 
	&
	\dot{m}^{eq|2}_{x_i^2}=\frac{c^2+2u_i^2}{3},\:\: i\in \llbracket{1,d}\rrbracket,\\
	&\dot{m}^{eq|2}_{x_ix_j}=\frac{2u_iu_j}{3},\:\: i,j\in \llbracket{1,d}\rrbracket,\: i\neq j,\\
	&\dot{m}^{eq|3}_{x_i^3}=cu_i,\:\: i\in \llbracket{1,d}\rrbracket,\\
	&\dot{m}^{eq|3}_{x_i^2x_j}=\frac{u_j}{3},\:\: j\in \llbracket{1,d}\rrbracket,\: i\neq j, \\
	&\dot{m}^{eq|3}_{x_ix_jx_l}=0,\:\: i,j,l\in \llbracket{1,d}\rrbracket,\:i\neq j,\: i\neq l,\: j\neq l.
\end{align} 
\end{subequations} 
Thus, as long as the conditions in Eq. (\ref{third-condition-2}) hold, one can conclude  
\begin{align}\label{third-macro}
	\partial_t\phi+\mathbf{u}\cdot\nabla\phi=R+O(\Delta t^3).
\end{align}
Here, we would like to point out that $\dot{m}^{eq|4}_{x_i^2x_j^2}$ ($i,j\in\llbracket1,d\rrbracket$, $i\neq j$), which correspond to the fourth-order moments of the EDFs $f_k^{eq}$ in the BGK-LB model (\ref{BGK-LB-Evolution-seq2}), still remain free.

Similar to the derivation of Eq. (\ref{eq-od3}), according to Eq. (\ref{eq-od3}), Eq. (\ref{ME-4}) can be rewritten as
\begin{align}\label{eq-od4}
	f_k=&f_k^{eq}-\Delta t\partial_tf_k^{eq}+\frac{\Delta t^2}{2}\partial_t^2f_k^{eq}-\frac{\Delta t^2}{4}\partial_t^2\Big(f_k^{eq}-\frac{\Delta t}{2}\mathbf{c}_k\cdot\nabla f_k^{eq}-\frac{\Delta t}{2}\partial_tf_k^{eq}+\frac{\Delta t}{2}R_k\Big)\notag\\
	&-\frac{\Delta t}{2}\mathbf{c}_k\cdot\nabla \Big(f_k^{eq}-\frac{\Delta t}{2}\mathbf{c}_k\cdot\nabla f_k^{eq}-\frac{\Delta t}{2}\partial_tf_k^{eq}+\frac{\Delta t}{2}R_k\Big)-\frac{\Delta t^2}{4}\mathbf{c}_k^{.2}\overset{2}{\cdot}\nabla^2\Big(f_k^{eq}-\frac{\Delta t}{2}\mathbf{c}_k\cdot\nabla f_k^{eq}-\frac{\Delta t}{2}\partial_tf_k^{eq}+\frac{\Delta t}{2}R_k\Big)\notag\\
	&+\frac{\Delta t}{2}\partial_t\Big(f_k^{eq}-\frac{\Delta t}{2}\mathbf{c}_k\cdot\nabla f_k^{eq}-\frac{\Delta t}{2}\partial_tf_k^{eq}+\frac{\Delta t}{2}R_k\Big) -\frac{\Delta t^3}{12}\partial_t^3f_k^{eq}-\frac{\Delta t^3}{12}\mathbf{c}_k^{.3}\overset{3}{\cdot} \nabla^3 f_k^{eq} +\frac{\Delta t}{2}R_k+O(\Delta t^4)\notag\\
	=&f_k^{eq}-\frac{\Delta t}{2}\mathbf{c}_k\cdot\nabla f_k^{eq}-\frac{\Delta t}{2}\partial_tf_k^{eq}+\frac{\Delta t^3}{8}\mathbf{c}_k\cdot\nabla \partial_t^2f_k^{eq}\notag\\
	&+\frac{\Delta t^3}{8}\mathbf{c}_k^{.2}\overset{2}{\cdot}\nabla^2\partial_tf_k^{eq}+\frac{\Delta t^3}{24}\partial_t^3f_k^{eq}+\frac{\Delta t^3}{24}\mathbf{c}_k^{.3}\overset{3}{\cdot} \nabla^3 f_k^{eq}+\frac{\Delta t}{2}R_k +O(\Delta t^4),
\end{align}	
 substituting Eq. (\ref{eq-od4}) into Eq. (\ref{taylor-expansion-eq}) yields
\begin{align}\label{eq-ode4-simp}
	&2f_k+\Delta t\mathbf{c}_k\cdot \nabla \Big(f_k^{eq} -\frac{\Delta t}{2}\partial_tf_k^{eq}+\frac{\Delta t^3}{8}\mathbf{c}_k\cdot\nabla \partial_t^2f_k^{eq}+\frac{\Delta t^3}{24}\mathbf{c}_k^{.2}\overset{2}{\cdot}\nabla^2\partial_tf_k^{eq}+\frac{\Delta t^3}{24}\partial_t^3f_k^{eq}\Big)\notag\\
	&\quad-\frac{\Delta t^3}{4}\mathbf{c}_k^{.2}\overset{2}{\cdot} \nabla^2\partial_tf_k^{eq}-\frac{\Delta t^3}{12}\mathbf{c}_k^{.3}\overset{3}{\cdot}\nabla^3  f_k^{eq} +O(\Delta t^5)\notag\\
	&=\Delta t\partial_tf_k
	-\frac{\Delta t^2}{2}\partial_t^2f_k
	+\frac{\Delta t^3}{6}\partial_t^3f_k
	-\frac{\Delta t^4}{24}\partial_t^4f_k\notag\\
	&\quad+2f_k^{eq}-2\Delta t\partial_tf_k^{eq}+\Delta t^2\partial_t^2f_k^{eq}-\frac{\Delta t^3}{3}\partial_t^3f_k^{eq}
	+\frac{\Delta t^4}{12}\partial_t^4f_k^{eq}+\Delta tR_k+O(\Delta t^5),
\end{align} 
then taking the zeroth-order moment of Eq. (\ref{eq-ode4-simp}) gives the following fourth-order modified equation compared  to the L-HE (\ref{HE}): 
\begin{align}\label{fourth-me}
	&\nabla\cdot\mathbf{u}\phi-\frac{1}{2}\mathbf{u}\cdot\nabla\phi_t+\frac{\Delta t^3}{24}\partial_t\Big(3\mathbf{q}^2\nabla^2 \partial_t+\mathbf{q}^3\nabla^3+ \mathbf{u}\cdot\nabla\partial_t^2\Big)\phi	-\frac{\Delta t^2}{12}\Big(3\mathbf{q}^{2}\overset{2}{\cdot} \nabla^2\partial_t\phi+\mathbf{q}^3\overset{3}{\cdot} \nabla^3\Big)\phi  \notag\\
	&=-\partial_t\phi
	+\frac{\Delta t}{2}\partial_t^2\phi
	-\frac{\Delta t^2}{6}\partial_t^3\phi+\frac{\Delta t^2}{24}\partial_t^4\phi+O(\Delta t^4),
\end{align} 
with the help of Eqs. (\ref{third-condition-2}) and (\ref{third-macro}), Eq. (\ref{fourth-me}) can be simplified as
  \begin{align}\label{fourth-me-2}
  	\partial_t\phi+\mathbf{u}\cdot\nabla\phi=R +O(\Delta t^4).
  \end{align} 
 This implies that when the tensors $\mathbf{q^2}$ and $\mathbf{q}^3$, which are related to the second- and third-order moments of the EDFs $f_k^{eq}$, satisfy the conditions in Eq. (\ref{third-condition-2}), a fourth-order unified BGK-LB model  with the fourth-order moments of the EDFs $f_k^{eq}$ remaining free can be derived. It should be noted that the choices of the fourth-order moments of the EDFs $f_k^{eq}$ are not arbitrary but are closely related to the stability of the BGK-LB model (see Sec. \ref{Sec-stability-BGK-LB} for details).
\subsection{The fourth-order initialization and boundary schemes}
The previous Part \ref{dte-analysis} only shows that the modified equation of the BGK-LB model (\ref{BGK-LB-Evolution-seq2}) can be consistent with the L-HE (\ref{HE}) up to the order of $O(\Delta t^4)$ when Eq. (\ref{third-condition-2}) holds. To ensure that the BGK-LB model (\ref{BGK-LB-Evolution-seq2}) can be completely fourth-order accurate in numerical simulations, the initialization and boundary schemes of the distribution functions $f_k$ in the BGK-LB model (\ref{BGK-LB-Evolution-seq2}) must be properly determined. In this work, we assume that the initial and boundary data are sufficiently smooth and only consider the Dirichlet [see Eq. (\ref{Diri}) below] and periodic boundary conditions. With the aid of the above accuracy analysis on the BGK-LB model (\ref{BGK-LB-Evolution-seq2}), for the given initial condition in the L-HE (\ref{HE}), we adopt the following scheme to approximate the unknown distribution functions $f_k(\mathbf{x},t)$:
\begin{equation}\label{ini}
	\begin{aligned}
		\:\: f_k(\mathbf{x},0)= f_k^{eq}(\mathbf{x},0)-\frac{\Delta t}{2}\mathbf{c}_k\cdot\nabla f_k^{eq}(\mathbf{x},0)-\frac{\Delta t}{2}\partial_tf_k^{eq}(\mathbf{x},0)+\frac{\Delta t}{2}R_k,\:\: \mathbf{x}\in\Omega,
		\end{aligned}
\end{equation}
here 
\begin{align}
	\partial_tf_k^{eq}(\mathbf{x},0)=R_k-\mathbf{u}\cdot\nabla f_k^{eq}(\mathbf{x},0),
\end{align} 
where Eq. (\ref{fourth-me-2}) has been used. Then,  Eq. (\ref{ini}) can be rewritten as
	\begin{align}\label{ini-simp}
		f_k(\mathbf{x},0)= f_k^{eq}(\mathbf{x},0)+\frac{\Delta t}{2}(\mathbf{u}-\mathbf{c}_k)\cdot\nabla f_k^{eq}(\mathbf{x},0),\:\: \mathbf{x}\in\Omega,
	\end{align}
from Eq. (\ref{eq-od3}), it can be seen that the initialization scheme in Eq.  (\ref{ini-simp}) is a third-order approximation of the distribution functions $f_k(\mathbf{x},0)$. This is sufficient and does not reduce the overall accuracy of the fourth-order BGK-LB model (\ref{BGK-LB-Evolution-seq2}) \cite{bellotti2023-weakly,BOGHOSIAN2024,strikwerda2004finite}. We would also like to point out that the commonly used initialization scheme $f_k = f_k^{eq}$ is only first-order accurate at the acoustic scaling. 

For the Dirichlet boundary condition of the L-HE (\ref{HE}), in this work, we consider 
\begin{align}\label{Diri}
	\phi(\mathbf{x}_b,t)=h(\mathbf{x}_b,t),\;\:\mathbf{x}_b\in\partial \Omega,
\end{align}
where $h(\mathbf{x}_b,t):$ $(\mathbb{R}^d,\mathbb{R})\rightarrow \mathbb{R}$ is a given continuous and smooth function at the boundary point $\mathbf{x}_b$ and time $t$, and $\partial \Omega $ is the boundary. Specifically,  $\mathbf{x}_b=\big(x_{1},\ldots,x_{j-1},x_{j},x_{j+1},\ldots,x_{d})$ with the $j_{th}$ $(j\in\llbracket1,d\rrbracket)$ element fixed while the others varying continuously in the computational domain $\Omega$. Now we develop the fourth-order full-way and half-way boundary schemes, which are respectively used in the cases where the boundary is located at the lattice node and in the middle of two lattice nodes (see Fig. \ref{node} below), to approximate the unknown distribution functions $f_k$. 
		
\begin{figure}[htbp]    
		\centering           
		{
			\includegraphics[width=0.7\textwidth]{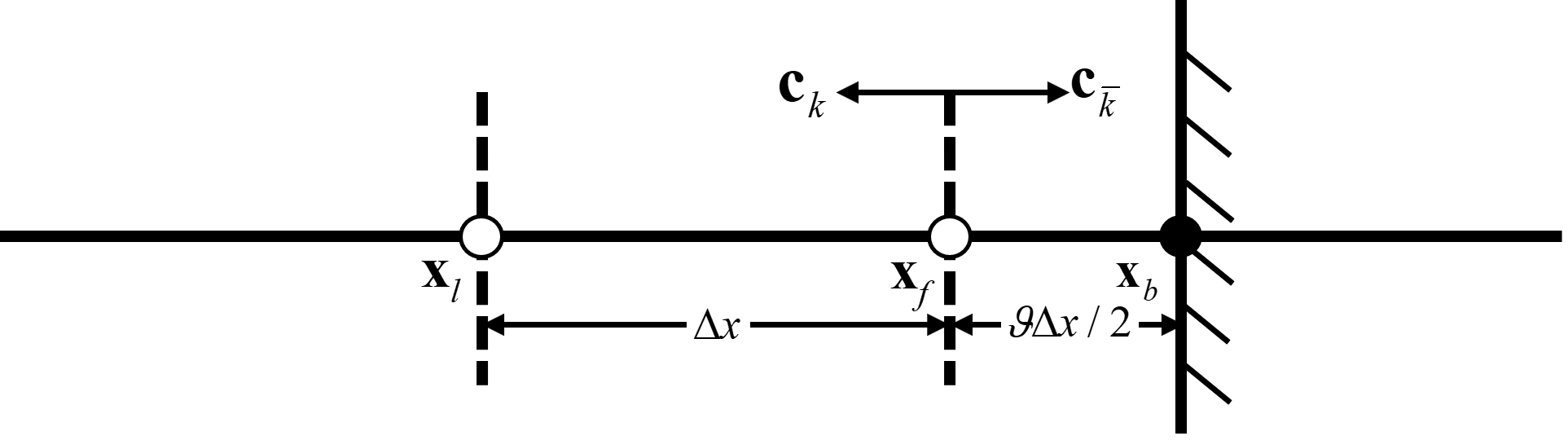} 
		}  
		\caption{The dotted straight line is the lattice line and the solid straight line is the boundary. White circles (◦) are the lattice nodes in the computational domain, the black circle (•) is the intersection of the boundary and the lattice line.  $\vartheta=0$ and $\vartheta=1$ represent the full-way and half-way cases, respectively.}       \label{node}   
	\end{figure} 
\begin{itemize}
	\item Full-way 
	
	In the full-way case, the internal lattice node $\mathbf{x}_f$ coincides with the boundary point $\mathbf{x}_b$. The corresponding  boundary scheme at  different accuracy orders can be designed as 
	\begin{subequations}\label{bun1}
		\begin{align} 
			&O(\Delta t):\: f_k(\mathbf{x}_b,t)=f_k^{eq}(\mathbf{x}_b,t),\\
			&O(\Delta t^3): \: f_k(\mathbf{x}_b,t)=f_k^{eq}(\mathbf{x}_b,t)-\frac{\Delta t}{2}\mathbf{c}_k\cdot\nabla f_k^{eq}(\mathbf{x}_b,t)-\frac{\Delta t}{2}\partial_tf_k^{eq}(\mathbf{x}_b,t)+\frac{\Delta t}{2}R_k,\\
			&O(\Delta t^4):\: f_k(\mathbf{x}_b,t)=   f_k^{eq}(\mathbf{x}_b,t)-\frac{\Delta t}{2}\mathbf{c}_k\cdot\nabla f_k^{eq}(\mathbf{x}_b,t)-\frac{\Delta t}{2}\partial_tf_k^{eq}(\mathbf{x}_b,t)+\frac{\Delta t^3}{8}\mathbf{c}_k\cdot\nabla \partial_t^2f_k^{eq}(\mathbf{x}_b,t)\notag\\
			&\qquad\qquad\qquad\qquad+\frac{\Delta t^3}{8}\mathbf{c}_k^{.2}\overset{2}{\cdot}\nabla^2\partial_tf_k^{eq}(\mathbf{x}_b,t)+\frac{\Delta t^3}{24}\partial_t^3f_k^{eq}(\mathbf{x}_b,t)+\frac{\Delta t^3}{24}\mathbf{c}_k^{.3}\overset{3}{\cdot} \nabla^3 f_k^{eq}(\mathbf{x}_b,t)+\frac{\Delta t}{2}R_k,\label{bun-1}
		\end{align} 
	\end{subequations} 
	where Eqs.  (\ref{ME-1}),  (\ref{eq-od3}), and (\ref{eq-od4}) have been used.

	\item Half-way  
	
	In the half-way case, $\mathbf{x}_f$ is the internal lattice node near the boundary point $\mathbf{x}_b$ with a distance of $\Delta x/2$ (see Fig. \ref{node} above). The corresponding first- to third-order half-way boundary schemes  can be designed as (the detailed construction can be found in \ref{Appendix-half-way})
	\begin{subequations}\label{bun2}
		\begin{align} 
			&O(\Delta t):\: f_k(\mathbf{x}_f,t)=  \big(f_k^{eq}+f_{\overline{k}}^{eq}\big) (\mathbf{x}_b,t),\\
			&O(\Delta t^2):\: f_k(\mathbf{x}_f,t)= \big(f_k^{eq}+f_{\overline{k}}^{eq}\big) (\mathbf{x}_b,t)+\frac{\Delta t}{2}\partial_t\big(f_k^{eq}+f_{\overline{k}}^{eq}\big)(\mathbf{x}_b,t)
			+\frac{\Delta t}{2}\big(R_k+ R_{\overline{k}}\big),\\
			&O(\Delta t^3):\:
			f_k(\mathbf{x}_f,t)=   \big(f_k^{eq}+f_{\overline{k}}^{eq}\big) (\mathbf{x}_b,t)+\frac{\Delta t}{2}\partial_t\big(f_k^{eq}+f_{\overline{k}}^{eq}\big)(\mathbf{x}_b,t)
			-\frac{\Delta t^2}{8}\mathbf{c}_k^{.2}\overset{2}{\cdot}\nabla^2\big(f_k^{eq}+f_{\overline{k}}^{eq}\big)(\mathbf{x}_b,t)
			\notag\\
			&\qquad\qquad\qquad\qquad-\frac{\Delta t^2}{4}\mathbf{c}_k\cdot\nabla\partial_t\big(f_k^{eq}-f_{\overline{k}}^{eq}\big)(\mathbf{x}_b,t) +\frac{\Delta t}{2}\big(R_k+ R_{\overline{k}}\big),\\
			&O(\Delta t^4):\: 	f_k(\mathbf{x}_f,t)=   \big(f_k^{eq}+f_{\overline{k}}^{eq}\big) (\mathbf{x}_b,t)+\frac{\Delta t}{2}\partial_t\big(f_k^{eq}+f_{\overline{k}}^{eq}\big)(\mathbf{x}_b,t)
			-\frac{\Delta t^2}{8}\mathbf{c}_k^{.2}\overset{2}{\cdot}\nabla^2\big(f_k^{eq}+f_{\overline{k}}^{eq}\big)(\mathbf{x}_b,t)
			\notag\\
			&\qquad\qquad\qquad\qquad-\frac{\Delta t^2}{4}\mathbf{c}_k\cdot\nabla\partial_t\big(f_k^{eq}-f_{\overline{k}}^{eq}\big)(\mathbf{x}_b,t)-\frac{\Delta t^3}{8}\mathbf{c}_k\cdot\nabla \partial_t^2\big(f_k^{eq}-f_{\overline{k}}^{eq}\big)(\mathbf{x}_b,t)
			\notag\\
			&\qquad\qquad\qquad\qquad
			-\frac{\Delta t^3}{16}\mathbf{c}_k^{.2}\overset{2}{\cdot}\nabla^2\partial_t\big(f_k^{eq}+f_{\overline{k}}^{eq}\big)(\mathbf{x}_b,t)-\frac{\Delta t^3}{24}\partial_t^3\big(f_k^{eq}+f_{\overline{k}}^{eq}\big)(\mathbf{x}_b,t)
			+\frac{\Delta t}{2}\big(R_k+ R_{\overline{k}}\big),\label{bun-2}
		\end{align}
	\end{subequations}  
	where $\overline{k}$ is defined such that $\mathbf{c}_{\overline{k}}=-\mathbf{c}_k$, $f_{\overline{k}}^{\star}$ represents the post-collision distribution function, i.e.,
	\begin{align}\label{collision}
		f_k^{\star}(\mathbf{x},t)=-f_k(\mathbf{x},t)+2f_k^{eq}(\mathbf{x},t)+\Delta tR_k,\:\: k\in \llbracket{1,2d^2+1}\rrbracket.
	\end{align}	 
	\end{itemize}

For the developed full-way and half-way boundary schemes for the  Dirichlet boundary condition (\ref{Diri}) of the L-HE (\ref{HE}), we give two remark as follows.\\
\noindent\textbf{Remark 1.}  From both the full-way and half-way  boundary schemes, i.e., Eqs. (\ref{bun1}) and (\ref{bun2}), it can be observed that there exists one spatial direction at the boundary $\partial\Omega$ for which the terms $\nabla f_k^{eq}$, $\nabla\partial_{t^2}^2f_k^{eq}$, $\nabla^2f_k^{eq}$, $\nabla^2\partial_tf_k^{eq}$, and $\nabla^3f_k^{eq}$  are unknown. Without loss of generality, we consider the $j_{th}$ spatial direction here. And according to Eq. (\ref{fourth-me-2}), neglecting the fourth-order truncation error, these terms in the $j_{th}$  spatial direction can be replaced by the following spatial and temporal derivatives of the EDFs $f_k^{eq}$:
\begin{subequations}\label{bun-1-spatial-deriva}
	\begin{align}
		&\partial_{x_j}f_k^{eq}(\mathbf{x}_b,t)=\Big(R_k-\partial_t-\sum_{l=1,l\neq j}^{d}u_{l}\partial_{x_l}\Big)f_k^{eq}(\mathbf{x}_b,t) ,\\
		&\partial^3_{t^2x_j}f_k^{eq}(\mathbf{x}_bt)=-\Big(\partial^3_{t^3}+\sum_{l=1,l\neq j}^{d}u_{l}\partial^3_{t^2x_l} \Big)f_k^{eq}(\mathbf{x}_b,t),\\ 	&\partial^2_{x_j}f_k^{eq}(\mathbf{x}_b,t)=\Bigg[\partial^2_{t}+2\sum_{l=1,l\neq j}^du_l\partial^2_{tx_l}+\sum_{l=1,l\neq j}^du_l\partial_{x_l}\bigg(\sum_{p=1,p\neq j}^du_p\partial_{x_p}\bigg)\Bigg]f_k^{eq}(\mathbf{x}_b,t),\\
		&\partial^3_{tx_j^2}f_k^{eq}(\mathbf{x}_b,t)=\Bigg[\partial^3_{t}+2\sum_{l=1,l\neq j}^du_l\partial^3_{t^2x_l}+\sum_{l=1,l\neq j}^du_l\partial^2_{tx_l}\bigg(\sum_{p=1,p\neq j}^du_p\partial_{x_p}\bigg)\Bigg]f_k^{eq}(\mathbf{x}_b,t),\\
		&\partial^3_{x_j}f_k^{eq}(\mathbf{x}_b,t)=-\Bigg[\partial^3_{t}+3\sum_{l=1,l\neq j}^{d}u_{l}\partial^3_{t^2x_l}
		+3\sum_{l=1,l\neq j}^{d}u_{l}\partial^2_{tx_l}\bigg(\sum_{p=1,p\neq j}^{d}u_{p}\partial_{x_p}\bigg)\notag\\
		&\qquad\qquad\qquad\quad+\sum_{l=1,l\neq j}^du_l\partial_{x_l}\bigg(\sum_{p=1,p\neq j}^du_p\partial_{x_p}\bigg)\bigg(\sum_{r=1,r\neq j}^{d}u_{r}\partial_{x_r}\bigg)\Bigg]f_k^{eq}(\mathbf{x}_b,t),
	\end{align} 
\end{subequations} 
where for the EDFs $f_k^{eq}$, the temporal derivative  and the spatial derivatives in the $l_{th}$ ($l\in \llbracket{1,d}\rrbracket$, $l\neq j$) directions   can be determined from the boundary condition in Eq. (\ref{Diri}). \\
\noindent\textbf{Remark 2.} One can find that the full-way boundary schemes in Eq. (\ref{bun1}) at different accuracy orders are actually the first- to fourth-order expansions of the distribution functions $f_k$ derived from the BGK-LB model (\ref{BGK-LB-Evolution-seq2}), i.e., Eqs. (\ref{ME-1}), (\ref{eq-od3}), and (\ref{eq-od4}). This implies that the certain properties, such as the stability, of the full-way boundary schemes in Eq. (\ref{bun1}) would be consistent with the bulk BGK-LB model (\ref{BGK-LB-Evolution-seq2}).  However, in the half-way case, since the boundary is located in the middle of two lattice nodes, the idea of the half-way anti-bounce-back scheme \cite{Zhao2019-2,ZHANG2019} is introduced in the construction of the half-way boundary scheme (see  \ref{Appendix-half-way} for details), and it can be observed that the half-way boundary schemes in Eq. (\ref{bun2}) are no longer consistent with Eqs. (\ref{ME-1}), (\ref{eq-od3}), and (\ref{eq-od4}). Thus, when the bulk BGK-LB model (\ref{BGK-LB-Evolution-seq2}) can be stabilized, compared to the half-way boundary schemes in Eq. (\ref{bun2}), the full-way boundary schemes in Eq. (\ref{bun1}), which are consistent with the bulk BGK-LB model (\ref{BGK-LB-Evolution-seq2}), would be more stable, as can be found  in Sec. \ref{Numer}.
\section{The stability analysis of the unified Bhatnagar-Gross-Krook lattice Boltzmann model}\label{Sec-stability-BGK-LB}
In this section, we will first conduct a theoretical analysis on the the microscopic entropy stability of the BGK-LB model  (\ref{BGK-LB-Evolution-seq2}). Then, we will further discuss the $L^2$ stability using the von Neumann analysis and numerically plot the stability regions of the BGK-LB model  (\ref{BGK-LB-Evolution-seq2}). 
\subsection{The entropy stability analysis}
In this part, we focus on the conservative L-HE (\ref{HE}), i.e., the source term $R\equiv0$, and the corresponding macroscopic entropy is chosen as the commonly used $S(\phi)=\phi^2/2$. For the developed BGK-LB model (\ref{BGK-LB-Evolution-seq2}), based on the entropy theory of the kinetic representation \cite{Bouchut1999,Bouchut2003}, we introduce the microscopic entropy given by the sum of the kinetic entropies as \cite{Dubois2013}
\begin{align}\label{entropy}
\varpi(\mathbf{f}):=\sum\Big(f_1,\ldots,f_{2d^2+1}\Big)=\sum_{k=1}^{2d^2+1}s_k(f_k),
\end{align}
where the microscopic entropies $s_k$ ($k\in\llbracket{1,2d^2+1}\rrbracket$) are the strictly convex functions of the distribution functions $f_k$ on $\mathbb{R}$. In particular, we have
\begin{align}
	S(\phi)=\min_{\phi=\sum_{k=1}^{2d^1+1}f_k}\sum\Big(f_1,\ldots,f_{2d^2+1}\Big)=\sum\Big(f_1^{eq},\ldots,f_{2d^2+1}^{eq}\Big),
\end{align} 
which indicates that the macroscopic entropy $S(\phi)$ is a constrained optimization of the microscopic entropy $\varpi(\mathbf{f})$, and in which the minimum is reached at the EDFs $f_k^{eq}$. Furthermore, the EDFs $f_k^{eq}$ and the kinetic entropies $s_k$ have the following relation:
\begin{align}\label{feq-sk}
	f_k^{eq}\big(\phi(\phi^*)\big)=\nabla_{\phi^{*}}s_k^{*}(\phi^*),
\end{align}
where $\phi^{*}$ is called the dual  or entropy variable. $s_k^{*}$ represent the conjugate or dual kinetic entropies.  In the regular case, when the entropy $s(\upsilon)$ is smooth and
strictly convex on $\mathbb{R}$, we have
\begin{subequations}\label{phi-s*}
	\begin{align}
		&\upsilon^*=\nabla s\big(\upsilon(\upsilon^*)\big),\:\:\upsilon=\nabla_{\upsilon^*} s^*\big(\upsilon^*(\upsilon)\big),\:\:\upsilon^{**}=\upsilon,\\
		&s^*(\upsilon^*(\upsilon))=\upsilon^*\cdot\big[\upsilon(\upsilon^*)\big]-s(\upsilon\big(\upsilon^*)\big),\:\:s(\upsilon(\upsilon^*))=\upsilon\cdot\big[\upsilon^*(\upsilon)\big]-s^*(\upsilon^*\big(\upsilon)\big),\\ 
		&s(\upsilon(\upsilon^*))=s^{**}\big(\upsilon^{**}(\upsilon^{*})\big)=\upsilon^{**}\cdot\big[\upsilon^*(\upsilon^{**})\big]-s^*(\upsilon^*\big(\upsilon^{**})\big).
	\end{align}
\end{subequations}
According to Eqs. (\ref{feq-sk}) and (\ref{phi-s*}), we can obtain 
\begin{align}\label{sk*}
\nabla_{(f_k^{eq})^*}s_k^{*}\big((f_k^{eq})^*\big)=f_k^{eq}=	\nabla_{\phi^*}s_k^{*}(\phi^*),\:\: k\in\llbracket1,2d^2+1\rrbracket,
\end{align}
which means that $(f_k^{eq})^*=\phi$ and $s_k^{*}(\phi)=\omega_k\phi^2/2+C$ with $C$ is a constant.
Then, we have
\begin{align}\label{sk}
	&s_k(f_k^{eq})=f_k^{eq}\phi-s_k^{*}=\frac{(f_k^{eq})^2}{2\omega_k}-C,\:\: k\in\llbracket1,2d^2+1\rrbracket, 
\end{align}
which satisfies $\sum_{k=1}^{2d^2+1}s_k(f_k^{eq})=S(\phi)$ provided that  $C=0$. Thus, the kinetic entropies  $s_k(f_k)$  are given by
\begin{align}\label{skfk}
	s_k(f_k)=\frac{f_k^2}{2\omega_k},\:\: k\in\llbracket1,2d^2+1\rrbracket.
\end{align} 
It should be noted that the kinetic entropies $s_k(f_k)$ must be strictly convex, which is also known as the sub-characteristic conditions. In this work, we refer to them as the microscopic entropy stability conditions. And from Eq. (\ref{skfk}), it is obvious that the  microscopic entropy stability conditions are $\omega_k>0$ ($k\in\llbracket{1,2d^2+1}\rrbracket$).

In  Part \ref{dte-analysis}, we know that $\dot{m}^{eq|4}_{x_i^2x_j^2}$ ($i,j\in \llbracket{1,d}\rrbracket,\: i\neq j$)  in the BGK-LB model (\ref{BGK-LB-Evolution-seq2})  still remain free. In the following, we will present the relation between $\dot{m}^{eq|4}_{x_i^2x_j^2}$ and the microscopic entropy stability of the  BGK-LB model (\ref{BGK-LB-Evolution-seq2}).
 
 \noindent
 \textbf{Proposition 1.} The BGK-LB model (\ref{BGK-LB-Evolution-seq2}) for the L-HE (\ref{HE}) is microscopic entropy-stable once the following constraints of $\dot{m}^{eq|4}_{x_i^2x_j^2}$ ($i,j\in \llbracket{1,d}\rrbracket,\: i\neq j$), which correspond to the fourth-order moment   of the EDFs $f_k^{eq}$, hold:
 \begin{subequations}\label{entropy-stable-condition}
 	\begin{align}
 		&\sum_{i,j=1,j>i}^d \dot{m}^{eq|4}_{x_i^2x_j^2}>\frac{(d-3)+2\sum_{i=1}^d \sigma_i^2}{3}c^4,\:\:i\in \llbracket{1,d}\rrbracket,\\
 		&\sum_{j=1,j\neq i}^d\dot{m}^{eq|4}_{x_i^2x_j^2}<\frac{-|(4-d)\sigma_i|+1+2\sigma_i^2}{3} c^4,\:\: i\in \llbracket{1,d}\rrbracket,\\ 
 		&\dot{m}^{eq|4}_{x_i^2x_j^2}>\max\Bigg\{\frac{-2\sigma_i\sigma_j+|\sigma_i+\sigma_{j}|}{3},\frac{2\sigma_i\sigma_j+|\sigma_i-\sigma_{j}|}{3}\Bigg\}c^4 ,\:\: i,j\in\llbracket{1,d}\rrbracket,\: j>i.
 	\end{align} 
 \end{subequations}
where $\sigma_{i}=u_i/c$ ($i\in \llbracket{1,d}\rrbracket$) with $d<6$.\\
\noindent
\textbf{Proof.} Firstly, we denote the inverse of the transform matrix $\mathbf{M}$ in Eq. (\ref{DdQq-M}) as $\mathbf{\tilde{M}}$, and further denote $\mathbf{\tilde{M}}_l^{\dot{m}^{eq|\tilde{n}_p}_{\tilde{p}}}$ with $\tilde{p}=\prod_{p=1}^dx_p^{n_p}$ and $\tilde{n}_p=\sum_{p=1}^d n_p$ ($n_p\in \llbracket{0,2}\rrbracket$) as the $l_{th}$  column of the matrix $\mathbf{\tilde{M}}$. Specifically,   $\mathbf{\tilde{M}}_l^{\dot{m}^{eq}_{\tilde{p}}}$ ($l\in \llbracket{1,2d^2+1}\rrbracket$) are expressed as (see  \ref{inverse-M} for the detailed derivation)
\begin{subequations}\label{m-inverse}
\begin{align}
	&\mathbf{\tilde{M}}_1^{\dot{m}^{eq|0}_{\prod_{p=1}^dx_p^0}}=\mathbf{e}_1^T,\:\: \mathbf{e}_1=(1,\mathbf{0})\in \mathbb{R}^{(2d^2+1)\times 1},\\
	&\mathbf{\tilde{M}}_{k+1}^{\dot{m}^{eq|1}_{x_k}}=\Big(\mathbf{0},\underbrace{1/2}_{(i+1)_{th}},\mathbf{0},\underbrace{-1/2}_{(i+d+1)_{th}},\mathbf{0}\Big)^T/c,\:\: k\in\llbracket{1,d}\rrbracket,\footnotemark[3]\\
	&\mathbf{\tilde{M}}_{k+d+1}^{m^{eq|2}_{x_k^2}}=\Big(-1,\mathbf{0},\underbrace{1/2}_{(i+1)_{th}},\mathbf{0},\underbrace{1/2}_{(i+d+1)_{th}},\mathbf{0}\Big)^T/c^2,\:\: k\in \llbracket{1,d}\rrbracket,\\
	&\mathbf{\tilde{M}}_{k}^{\dot{m}^{eq|2}_{x_ix_j}}=\Big(\mathbf{0},\underbrace{1/4}_{l_{th}},
	\underbrace{-1/4}_{(l+1)_{th}},\underbrace{1/4}_{(l+2)_{th}},\underbrace{-1/4}_{(l+3)_{th}},\mathbf{0}\Big)^T/c^2,\:\: k=1+d+\tilde{n}+r,\:\: l=\tilde{d}+4(\tilde{n}+r), \: j>i,\\
	&\mathbf{\tilde{M}}_{k}^{\dot{m}^{eq|3}_{x_i^2x_j}}=\Big(\mathbf{0},\underbrace{-1/2}_{(j+1)_{th}},
	\mathbf{0},\underbrace{1/2}_{(d+1+j)_{th}},\mathbf{0},\underbrace{{\makebox{sgn}}\big(\mathbf{c}^{x_j}_{l},\mathbf{c}^{x_j}_{l+1},\mathbf{c}^{x_j}_{l+2},\mathbf{c}^{x_j}_{l+3}\big)/4}_{l_{th}\sim(l+3)_{th}},\mathbf{0}\Big)^T/c^3,\:\: k=\overline{d}-d(d-1)+\tilde{n}+|r|,\:\: l=\tilde{d}+4(\tilde{n}+|r|),\\
	&\mathbf{\tilde{M}}_{k}^{\dot{m}^{eq|4}_{x_i^2x_j^2}}=\Big(1,\mathbf{0},\underbrace{-1/2}_{(i+1)_{th}},
	\mathbf{0},\underbrace{-1/2}_{(j+1)_{th}},\mathbf{0},\underbrace{-1/2}_{(d+1+i)_{th}},\mathbf{0},\underbrace{-1/2}_{(d+j+1)_{th}},\mathbf{0},\underbrace{\mathbf{I}_{1\times 4}/4}_{k_{th}\sim (k+3)_{th}},\mathbf{0}\Big)^T/c^3,\:\: k=\overline{d}+\tilde{n}+r,\:\: l=\tilde{d}+4(\tilde{n}+r), \: j>i,
\end{align} 
\end{subequations}
where $\tilde{n}=\sum_{p=0}^{\min(i,j)-1}(d-p)$, $r=j-i$ ($i,j\in \llbracket{1,d}\rrbracket$), $\tilde{d}=-2(1+d)$, and $\overline{d}=(3d^2-d+2)/2$.  
\footnotetext[3]{The notation $J=\big(\ldots,\underbrace{a}_{p_{th}}\ldots\big)$ indicates that $a$ is the $p_{th}$ element of the vector $J$.}

Then, based on the relation of  $\bm{\omega}=\mathbf{\tilde{M}}\bm{\varepsilon}$, according to Eq. (\ref{epsilon}) and the unified expression of the matrix $\mathbf{\tilde{M}}$ in Eq.  (\ref{m-inverse}), we have
\begin{subequations}\label{omega_k-1}
	\begin{align}
		&\omega_1=1-\frac{\sum_{i=1}^d\dot{m}^{eq|2}_{x_i^2}}{c^2}+\frac{\sum_{i,j=1,i< j}^{d}\dot{m}^{eq|4}_{x_i^2x_j^2}}{c^2},\\
		&\omega_{k}=\frac{{\makebox{sgn}}\big(\mathbf{c}_{k}^{x_i}\big)u_i}{2c}+\frac{\dot{m}^{eq|2}_{x_i^2}}{2c^2}-\frac{{\makebox{sgn}}\big(\mathbf{c}_{k}^{x_i}\big)\sum_{j=1,j\neq i}^d\dot{m}^{eq|3}_{x_j^2x_i}}{2c}-\frac{\sum_{j=1,j\neq i}^d\dot{m}^{eq|4}_{x_i^2x_j^2}}{2c^4},\:\: k=i+1,\:\: i\in \llbracket{1,d}\rrbracket,\\
	&\omega_{k}=\frac{{\makebox{sgn}}\big(\mathbf{c}_{k}^{x_i}\big)u_i}{2c}+\frac{\dot{m}^{eq|2}_{x_i^2}}{2c^2}-\frac{{\makebox{sgn}}\big(\mathbf{c}_{k}^{x_i}\big)\sum_{j=1,j\neq i}^d\dot{m}^{eq|3}_{x_j^2x_i}}{2c}-\frac{\sum_{j=1,j\neq i}^d\dot{m}^{eq|4}_{x_i^2x_j^2}}{2c^4},\:\: k=i+1+d,\:\: i\in \llbracket{1,d}\rrbracket,\\
	&\omega_{k}=\frac{\dot{m}^{eq|4}_{x_i^2x_j^2}}{4c^4}+\frac{{\makebox{sgn}}\big(\mathbf{c}^{x_i}_k\mathbf{c}^{x_j}_k\big)\dot{m}^{eq|2}_{x_ix_j}}{4c^2}+\frac{{\makebox{sgn}}\big(\mathbf{c}_{k}^{x_i}\big)\dot{m}^{eq|3}_{x_j^2x_i}+
	{\makebox{sgn}}\big(\mathbf{c}_{k}^{x_j}\big)\dot{m}^{eq|3}_{x_i^2x_j}}{4c},\: k=l+1,\: l=\tilde{d}+4(\tilde{n}+|r|),\: i,j\in \llbracket{1,d}\rrbracket,\: j>i,
	\end{align} 
\end{subequations}
in particular, from Eqs. (\ref{0-1-meq}) and (\ref{third-condition-2}), we know that the moment vector $\mathbf{m}^{eq}=\bm{\varepsilon}\phi$ with $\bm{\varepsilon}$  being  expressed  as
 \begin{align}\label{meq}
\bm{\varepsilon}=\bigg[1,\mathbf{u},\frac{c^2+2\mathbf{u}^{.2}}{3},\frac{2(\mathbf{uu})_{\alpha<\beta}}{3},\frac{\mathbf{\overline{u}}_1}{3},\frac{\mathbf{\overline{u}}_2}{3},\ldots,\frac{\mathbf{\overline{u}}_d}{3},\mathbf{\dot{m}}^{eq|4}_{(\mathbf{x}^{.2}\mathbf{x}^{.2})_{\alpha<\beta}}\bigg]^T,
\end{align}  
here  
\begin{subequations}
	\begin{align}
		&(\mathbf{uu})_{\alpha<\beta}=\big(u_1u_2,u_1u_3,\ldots,u_1u_d,u_2u_3,u_2u_4,\ldots,u_2u_d,\ldots,u_{d-1}u_d\big),\\
		&\mathbf{\overline{u}}_i=\big(u_1,\ldots,u_{i-1},u_{i+1},\ldots,u_d\big),\:\: i\in\llbracket 1,d\rrbracket.
	\end{align}
\end{subequations}   After substituting Eq. (\ref{meq}) into the Eq. (\ref{omega_k-1}), we can obtain
\begin{subequations}\label{omega_k-2}
	\begin{align}
		&\omega_1=\frac{3c^4-dc^4-2c^2\sum_{i=1}^d u_i^2+3\sum_{i,j=1,j>i}^d\dot{m}^{eq|4}_{x_i^2x_j^2}}{3c^4},\\
		&\omega_{k}=\frac{c^3(4-d)u_i+c^4+2c^2u_i^2-3\sum_{j=1,j\neq i}^d\dot{m}^{eq|4}_{x_i^2x_j^2}}{6c^4},\:\: k=i+1,\:\: i\in \llbracket{1,d}\rrbracket,\label{omega_k-2-k}\\
		&\omega_{k}=\frac{-c^3(4-d)u_i+c^4+2c^2u_i^2-3\sum_{j=1,j\neq i}^d \dot{m}^{eq|4}_{x_i^2x_j^2}}{6c^4},\:\: k=i+1+d,\:\: i\in\llbracket{1,d}\rrbracket,\\
		&\omega_{k}=\frac{3\dot{m}^{eq|4}_{x_i^2x_j^2}+2c^2{\makebox{sgn}}\big(\mathbf{c}_{k}^{x_i}\mathbf{c}_{k}^{x_j}\big)u_iu_j+c^3\Big[{\makebox{sgn}}\big(\mathbf{c}_{k}^{x_i}\big)u_i+
			{\makebox{sgn}}\big(\mathbf{c}_{k}^{x_j}\big)u_j\Big]}{12c^4} ,\: k=1+\tilde{d}+4(\tilde{n}+|r|),\: i,j\in\llbracket{1,d}\rrbracket,\: j>i.
	\end{align} 
\end{subequations}
Thus, the BGK-LB model (\ref{BGK-LB-Evolution-seq2}) is microscopic entropy-stable as long as the following constraints of $\dot{m}^{eq|4}_{x_i^2x_j^2}$ are satisfied:
\begin{subequations}\label{omega>0}
	\begin{align}
		&\sum_{i,j=1,j>i}^d\dot{m}^{eq|4}_{x_i^2x_j^2}>\frac{(d-3)c^4+2c^2\sum_{i=1}^d u_i^2}{3},\\
		&\sum_{j=1,j\neq i}^d\dot{m}^{eq|4}_{x_i^2x_j^2}<\frac{c^3(4-d)u_i+c^4+2c^2u_i^2}{3c^4},\:\: k=i+1,\:\: i\in \llbracket{1,d}\rrbracket,\\
		&\sum_{j=1,j\neq i}^d\dot{m}^{eq|4}_{x_i^2x_j^2}<\frac{-c^3(4-d)u_i+c^4+2c^2u_i^2}{3c^4},\:\: k=i+d+1,\:\: i\in \llbracket{1,d}\rrbracket,\\
		&\dot{m}^{eq|4}_{x_i^2x_j^2}>\frac{2c^2{\makebox{sgn}}\big(\mathbf{c}_{k}^{x_i}\mathbf{c}_{k}^{x_j}\big)u_iu_j+c^3\Big[{\makebox{sgn}}\big(\mathbf{c}_{k}^{x_i}\big)u_i+
			{\makebox{sgn}}\big(\mathbf{c}_{k}^{x_j}\big)u_j\Big]}{3c^4} ,\: k=1+\tilde{d}+4(n+|r|),\: i,j\in\llbracket{1,d}\rrbracket\: j>i.
	\end{align} 
\end{subequations}
According to the lattice structure in Eq. (\ref{DdQq-M}), the above  Eq. (\ref{omega>0}) can be simplified as  Eq. (\ref{entropy-stable-condition}). Here, we would like to point out that to ensure that   Eq. (\ref{entropy-stable-condition}) makes sense, the constraints in Eq. (\ref{entropy-stable-condition}) must hold at least when $\sigma_i = 0$ ($i\in\llbracket{1,d}\rrbracket$), i.e., 
\begin{subequations}\label{sigma=0} 
	\begin{align}
		&\sum_{i,j=1,j>i}^d\dot{m}^{eq|4}_{x_i^2x_j^2}>\frac{(d-3)}{3}c^4,\:\:i\in \llbracket{1,d}\rrbracket,\\
		&\sum_{j=1,j\neq i}^d\dot{m}^{eq|4}_{x_i^2x_j^2}<\frac{1}{3}c^4,\:\: i\in \llbracket{1,d}\rrbracket,\\ 
		&\dot{m}^{eq|4}_{x_i^2x_j^2}>0 ,\:\: i,j\in\llbracket{1,d}\rrbracket,\: j>i.\label{meq-eq0}
	\end{align} 
\end{subequations}
From the above Eq. (\ref{sigma=0}), one can conclude that the spatial dimension $d$ must be less than 6, i.e., $d<6$, and the proof of this proposition is completed.
	
\begin{figure}     
	\centering     
		\subfloat[$\sigma_3=0$]    
		{
			\includegraphics[width=0.3\textwidth]{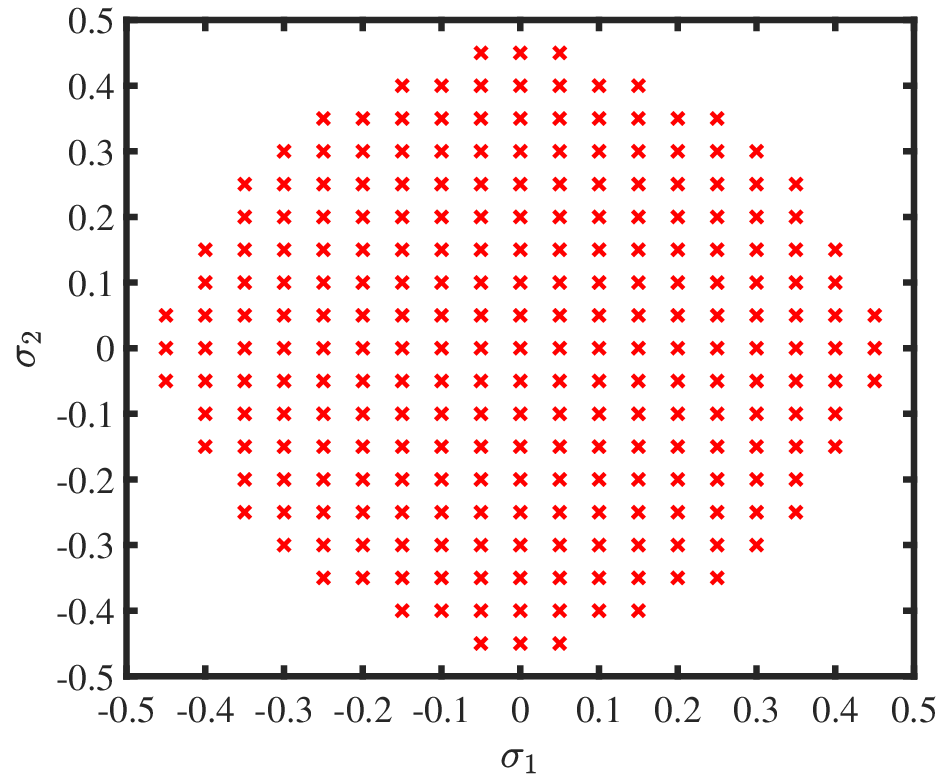} 
		}
		\subfloat[$\sigma_3=0.1$]    
		{
			\includegraphics[width=0.3\textwidth]{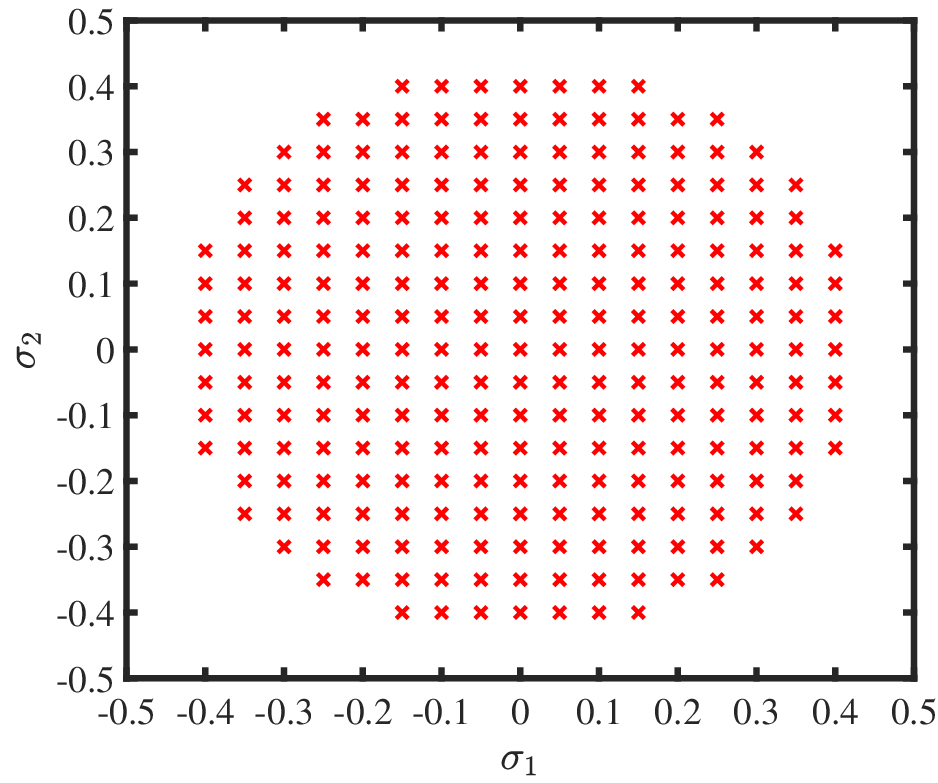} 
		} 
		
		\subfloat[$\sigma_3=0.2$]    
		{
			\includegraphics[width=0.3\textwidth]{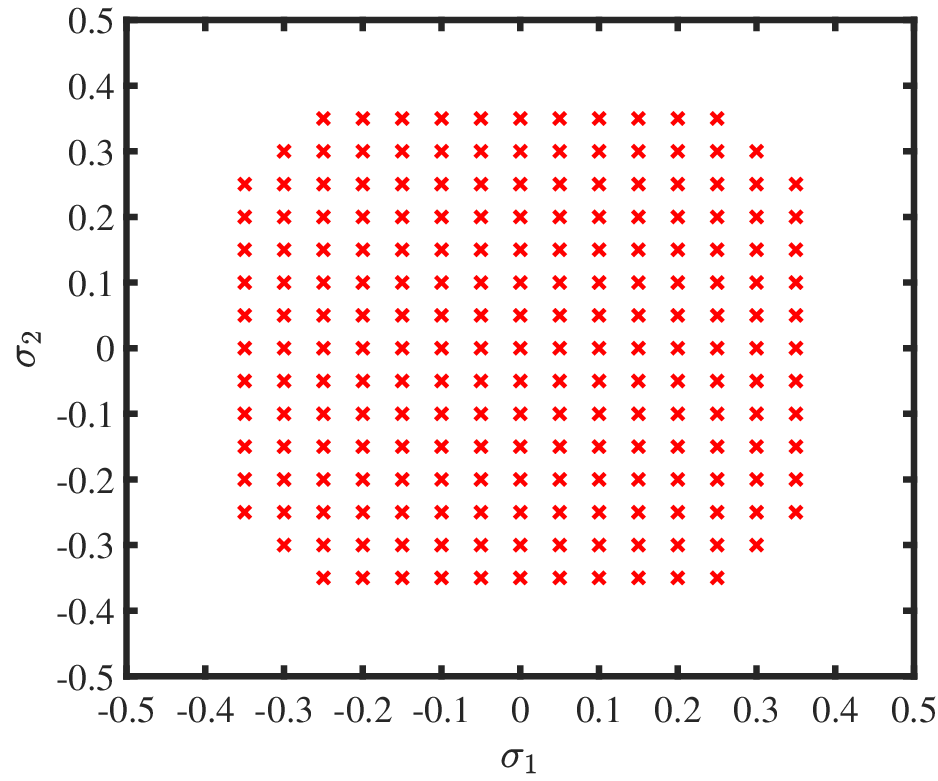} 
		}
		\subfloat[$\sigma_3=0.4$]    
		{
			\includegraphics[width=0.3\textwidth]{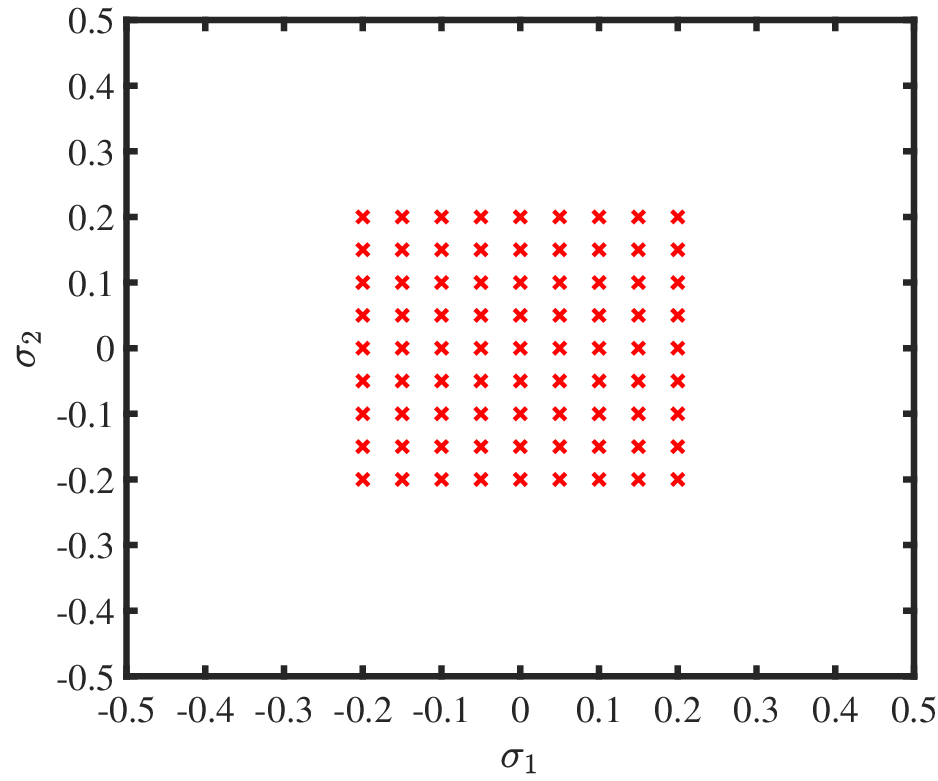} 
		}
		\caption{The microscopic entropy stability region  of the BGK-LB model  for the three-dimensional L-HE. } 
		\label{entropy-stable-fig}  
\end{figure}

In the following, we give some remarks on the above Proposition 1.\\
\noindent\textbf{Remark 3.} For the L-HE (\ref{HE}), which is a special case of the general convection-diffusion equations, the commonly used EDFs in the LB method are designed by the following linear form  
 \begin{align}\label{feq-1rd}
 	f_k^{eq}=\zeta_k\phi\bigg(1+\frac{\mathbf{c}^{x_k}\cdot\mathbf{u}}{c_s^2}\bigg),\:\: k\in\llbracket1,2d^2+1\rrbracket,
 \end{align}
 or the quadratic form \cite{Chai2016}
 \begin{align}\label{feq-2rd}
 	f_k^{eq}=\zeta_k\phi\Bigg(1+\frac{\mathbf{c}^{x_k}\cdot\mathbf{u}}{c_s^2}+\frac{\mathbf{uu}:\Big(\mathbf{c}^{x_k}\mathbf{c}^{x_k}-c_s^2\mathbf{I}_{d}\Big)}{2c_s^4}\Bigg),\:\: k\in\llbracket1,2d^2+1\rrbracket,
 \end{align}
 where $c_s$ is a model parameter related to the lattice velocity $c$, $\zeta_k$ are the weight coefficients, and $\mathbf{I}_{d}\in\mathbb{R}^{d\times d}$ is the identity matrix. From Eq. (\ref{omega_k-2}), it is obvious that the EDFs $f_k^{eq}$ proposed in our work are different from Eqs. (\ref{feq-1rd}) and (\ref{feq-2rd}), which means that the commonly used linear and quadratic EDFs $f_k^{eq}$ in the LB model are no longer applicable in the higher-order LB model for the L-HE (\ref{HE}). What is more, from the above discussion,  it is evident that the form of the EDFs $f_k^{eq}$ in the LB method is not unique, this is attributed to the fact that the number of required moment conditions for the EDFs $f_k^{eq}$ is generally lower than the number of the EDFs $f_k^{eq}$  themselves.  \\ 
\noindent\textbf{Remark 4.} For the one- and two-dimensional cases, we here prove that the BGK-LB model (\ref{BGK-LB-Evolution-seq2}) can be microscopic entropy-stable as long as $|\sigma_i|<1/2$  ($i\in\llbracket{1,d}\rrbracket$).
\begin{itemize}
	\item For the one-dimensional case, $\omega_k$ ($k\in\llbracket 1,3\rrbracket$) are expressed as
	\begin{subequations}
		\begin{align}
			&\omega_1=\frac{2-2\sigma^2_1}{3},\\
			&\omega_2=\frac{\sigma_1}{2}+\frac{2\sigma_1^2+1}{6},\\
			&\omega_3=-\frac{\sigma_1}{2}+\frac{2\sigma_1^2+1}{6},
		\end{align}
	\end{subequations}
	it is obvious that $\bm{\omega}>0$ as long as $|\sigma_1|<1/2$. 
	\item For the two-dimensional case, from Eq. (\ref{entropy-stable-condition}), we have
	\begin{subequations}
		\begin{align} 
			&\dot{m}^{eq|4}_{x_1^2x_2^2}<\min\Bigg\{\frac{-|2\sigma_1|+1+2\sigma_1^2}{3},\frac{-|2\sigma_2|+1+2\sigma_2^2}{3}\Bigg\}c^4,\\ 
			&\dot{m}^{eq|4}_{x_1^2x_2^2}>\max\Bigg\{\frac{-1+2(\sigma_1^2+\sigma_2^2)}{3},\frac{-2\sigma_1\sigma_2+|\sigma_1+\sigma_2|}{3},\frac{2\sigma_1\sigma_2+|\sigma_1-\sigma_2|}{3}\Bigg\}c^4 .
		\end{align} 
	\end{subequations}
	which is equivalent to requiring that $\sigma_{1,2}$ must satisfy
		\begin{subequations}
		\begin{align} 
			1+2\min\Bigg\{-|\sigma_1|+\sigma_1^2,-|\sigma_2|+\sigma_2^2\Bigg\}>\max\Bigg\{2(\sigma_1^2+\sigma_2^2)-1,|\sigma_1+\sigma_2|-2\sigma_1\sigma_2,|\sigma_1-\sigma_2|+2\sigma_1\sigma_2\Bigg\} ,
		\end{align} 
	\end{subequations}
	i.e.,
		\begin{subequations}\label{four-condition}
		\begin{align} 
			2\sigma_1^2-2|\sigma_1|+1-|\sigma_1+\sigma_2|+2\sigma_1\sigma_2>0,\label{eq-1}\\ 
			2\sigma_1^2-2|\sigma_1|+1-|\sigma_1-\sigma_2|-2\sigma_1\sigma_2>0,\label{eq-2}\\ 
			2\sigma_2^2-2|\sigma_2|+1-|\sigma_1+\sigma_2|+2\sigma_1\sigma_2>0,\label{eq-3}\\ 
			2\sigma_2^2-2|\sigma_2|+1-|\sigma_1-\sigma_2|-2\sigma_1\sigma_2>0,\label{eq-4}\\ 
			1-|\sigma_1|-\sigma_2^2>0,\label{eq-5}\\ 
			1-|\sigma_2|-\sigma_1^2>0,\label{eq-6}
		\end{align} 
	\end{subequations}
	and it is obvious that $|\sigma_{1,2}|<1$ are necessary for Eqs. (\ref{eq-5}) and (\ref{eq-6}). 
	
	To derive the regions of $\sigma_{1,2}$ which ensure that Eq. (\ref{four-condition}) holds, the detailed analysis can be divided into the following two cases.

\textbf{	Case 1.} We first consider the case of $|\sigma_{1,2}|\neq 1/2$, from Eqs. (\ref{eq-1}) and (\ref{eq-3}), we have
\begin{align}
	\begin{cases} 
		\sigma_1+\sigma_2\geq0:&2\sigma_1^2-2|\sigma_1|+1-\sigma_1-\sigma_2+2\sigma_1\sigma_2>0,\\ 
		\sigma_1+\sigma_2\leq0:&2\sigma_1^2-2|\sigma_1|+1+\sigma_1+\sigma_2-2\sigma_1\sigma_2>0, \\	\sigma_1+\sigma_2\geq0:&2\sigma_2^2-2|\sigma_2|+1-\sigma_1-\sigma_2+2\sigma_1\sigma_2>0,\\ 
		\sigma_1+\sigma_2\leq0:&2\sigma_2^2-2|\sigma_2|+1+\sigma_1+\sigma_2-2\sigma_1\sigma_2>0,
	\end{cases} 
\end{align}
which is equivalent to  
\begin{align}
	\begin{cases} 
			\sigma_1+\sigma_2\geq0:&\big(A_1\cup A_2\big)\cap\big(A_3\cup A_4\big),\\ 
		\sigma_1+\sigma_2\leq0:&\big(B_1\cup B_2\big)\cap\big(B_3\cup B_4\big),
	\end{cases} 
\end{align}
where
\begin{subequations}\label{set-all}
\begin{align} 
		&A_1=\Bigg\{-\sigma_1\leq \sigma_2\leq \frac{2\sigma_1^2-2|\sigma_1|-\sigma_1+1}{1-2\sigma_1},-1< \sigma_1< 1/2\Bigg\},\\
		&A_2=\Bigg\{1-\sigma_1\leq \sigma_2\leq 1,1/2< \sigma_1<1\Bigg\},\\ 
		&A_3=\Bigg\{\frac{2\sigma_1^2-2|\sigma_1|+\sigma_1+1}{-1-2\sigma_1}\leq \sigma_2\leq -\sigma_1,-1/2< \sigma_1< 1\Bigg\},\\
		&A_4=\Bigg\{-1< \sigma_2\leq -1-\sigma_1,-1< \sigma_1<-1/2\Bigg\}, \\
		&B_1=\Bigg\{-\sigma_2\leq \sigma_1\leq \frac{2\sigma_2^2-2|\sigma_2|-\sigma_2+1}{1-2\sigma_2},-1< \sigma_2< 1/2\Bigg\},\\
		&B_2=\Bigg\{1-\sigma_2\leq \sigma_1\leq 1,1/2< \sigma_2<1\Bigg\},\\ 
		&B_3=\Bigg\{\frac{2\sigma_2^2-2|\sigma_2|+\sigma_2+1}{-1-2\sigma_2}\leq \sigma_1\leq -\sigma_2,-1/2< \sigma_2< 1\Bigg\},\\
		&B_4=\Bigg\{-1< \sigma_1\leq -1-\sigma_2,-1< \sigma_2<-1/2\Bigg\}. 
\end{align}
\end{subequations}
From the above Eq. (\ref{set-all}), we would like to point out that 
\begin{subequations}
	\begin{align} 
		&A_{1\cap4}:=A_1\cap A_4=\Bigg\{1-\sigma_2=\sigma_1<1/2,1/2<\sigma_2<1\Bigg\},\\
		&A_{2\cap3}:=A_2\cap A_3=\Bigg\{1-\sigma_1=\sigma_2<1/2,1/2<\sigma_1<1\Bigg\},\\
		&B_{1\cap4}:B_1\cap B_4=\Bigg\{-1/2<\sigma_1=-\sigma_2-1,-1<\sigma_2<-1/2\Bigg\},\\
		&B_{2\cap3}:B_2\cap B_3=\Bigg\{-1/2<\sigma_2=-\sigma_1-1,-1<\sigma_1<-1/2\Bigg\},\\ 
		&A_{1\cup3}:=A_1\cup A_3=\Bigg\{-\sigma_1\leq \sigma_2<1/2,-\sigma_2\leq \sigma_1<1/2\Bigg\},\\
		&A_{2\cup4}:=A_2\cup A_4=\Bigg\{1/2<\sigma_{1,2}<1\Bigg\},\\
		&B_{1\cup3}:=B_1\cup B_3=\Bigg\{-1/2<\sigma_2\leq -\sigma_1,-1/2<\sigma_1<-\sigma_2\Bigg\},\\
		&B_{2\cup4}:=B_2\cup B_4=\Bigg\{-1<\sigma_{1,2}<-1/2\Bigg\}.
	\end{align}
\end{subequations}
Thus, Eqs. (\ref{eq-1}) and (\ref{eq-3}) hold in the case of  $|\sigma_{1,2}|\neq 1/2$ as long as $|\sigma_{1,2}|\in \overline{AB}$, where the set $\overline{AB}$ is given by
\begin{align} 
\overline{AB}&=\Big[\big(A_1\cup A_2\big)\cap\big(A_3\cup A_4\big)\Big]\cup\Big[\big(B_1\cup B_2\big)\cap\big(B_3\cup B_4\big)\Big]\notag\\
&=\bigg(\big(A_{1\cap 3}\cup A_{1\cap 4}\big)\cup\big(A_{2\cap 3}\cup A_{2\cup 4}\big)\bigg)\cup\bigg(\big(B_{1\cap 3}\cup B_{1\cap 4}\big)\cup\big(B_{2\cap 3}\cup B_{2\cup 4}\big)\bigg)\notag\\
&=C_1\cup C_2\cup C_3\cup C_4,
\end{align}
here
\begin{subequations}
	\begin{align} 
		&C_{1}=\bigg\{-1/2<\sigma_{1,2}<1/2\bigg\},\\
		&C_{2}=\bigg\{1/2<\sigma_{1,2}<1,-1<\sigma_{1,2}<-1/2\bigg\},\\
		&C_{3}=\bigg\{\sigma_1+\sigma_2=1,0\leq\sigma_1<1,\sigma_1\neq 1/2\bigg\},\\
		&C_{4}=\bigg\{\sigma_1+\sigma_2=-1,-1<\sigma_1\leq 0,\sigma_1\neq -1/2\bigg\}. 
	\end{align}
\end{subequations}
It is straightforward to validate that when $\sigma_{1,2}\in C_1$, Eqs. (\ref{eq-2}), (\ref{eq-4}), (\ref{eq-5}), and (\ref{eq-6}) hold. However, when $\sigma_{1,2}\in C_2\cup  C_3\cup C_4$, we find that Eqs. (\ref{eq-2}) and (\ref{eq-4}) cannot hold (see \ref{Appen-entropy-2D} for details). Thus, for the case of $|\sigma_{1,2}|\neq 1/2$, Eq. (\ref{four-condition}) holds as long as $|\sigma_{1,2}|<1/2$.
  
\textbf{Case 2.} For the case of $|\sigma_{1,2}|=1/2$, here we only consider $\sigma_1=1/2$ (the analysis on the other cases of $\sigma_1=-1/2$, $\sigma_2=1/2$, and $\sigma_2=-1/2$ is similar), Eqs. (\ref{eq-1}), (\ref{eq-2}), (\ref{eq-5}), and (\ref{eq-6}) become
  	\begin{subequations}\label{four-condition-u-1/2}
  	\begin{align} 
  		1/2-|1/2+\sigma_2|+\sigma_2>0,\label{eq-21}\\ 
  		1/2-|1/2-\sigma_2|-\sigma_2>0,\label{eq-22}\\  
  		1/2-\sigma_2^2>0,\label{eq-25}\\ 
  		3/4-|\sigma_2|>0,\label{eq-26}
  	\end{align} 
  \end{subequations}
it is obvious that the above Eqs. (\ref{eq-21}), (\ref{eq-22}), (\ref{eq-25}), and (\ref{eq-26})   cannot hold simultaneously. Thus, $|\sigma_{1,2}|\neq 1/2$.
 
In combination with the conclusions of cases 1 and 2, we prove that when the spatial dimension $d=2$, the BGK-LB model (\ref{BGK-LB-Evolution-seq2}) is  microscopic entropy-stable as long as $|\sigma_{1,2}|<1/2$.
\end{itemize}
\noindent\textbf{Remark 5.} It should be noted that, when the spatial dimension $d\geq 3$,  the number of the  fourth-order moments of the EDFs $f_k^{eq}$ is greater than 1, which makes it difficult to provide a unified approach to determine $\dot{m}^{eq|4}_{x_i^2x_j^2}$  $(i,j\in\llbracket{1,d}\rrbracket,j>i)$ from Eq. (\ref{entropy-stable-condition}). In this work, for the sake of  simplicity, let us assume that $\dot{m}^{eq|4}_{x_i^2x_j^2}$ are equal for every $i,j\in\llbracket{1,d}\rrbracket\:(j>i)$. Under this premise, Eq. (\ref{entropy-stable-condition}) can be reduced to
\begin{align}\label{psi-4-equal}
	\max\Bigg\{\frac{|\sigma_{i}\pm\sigma_{j}|\mp2\sigma_{i}\sigma_{j}}{3},
\frac{-6+2d+4\sum_{i=1}^d\sigma_{i}^2}{3d(d-1)}\Bigg\}<\dot{m}^{eq|4}_{x_i^2x_j^2}<\min\Bigg\{\frac{1+2\sigma_{i}^2-|(4-d)\sigma_{i}|}{3(d-1)}\Bigg\}, \:\: i,j\in \llbracket{1,d}\rrbracket,\: d<6,
\end{align} 
and in this case, $\dot{m}^{eq|4}_{x_i^2x_j^2}$  $(i,j\in\llbracket{1,d}\rrbracket,j>i)$ can be given by
\begin{align}\label{meq-4rd}
	\dot{m}^{eq|4}_{x_i^2x_j^2}=\max\Bigg\{\frac{|\sigma_{x_i}\pm\sigma_{x_j}|\mp2\sigma_{x_i}\sigma_{x_j}}{3},
	\frac{-6+2d+4\sum_{i=1}^d\sigma_{x_i}^2}{3d(d-1)}\Bigg\}gc^4+\min\Bigg\{\frac{1+2\sigma_{x_i}^2-|(4-d)\sigma_{x_i}|}{3(d-1)}\Bigg\}(1-g)c^4,
\end{align}
where $g\in(0,1)$. In particular, the BGK-LB model (\ref{BGK-LB-Evolution-seq2}) is  microscopic entropy-stable as long as $\sigma_i$ ($i\in\llbracket1,d\rrbracket$) satisfy
\begin{align}\label{region-d>=3}
	\max\Bigg\{\frac{|\sigma_{i}\pm\sigma_{j}|\mp2\sigma_{i}\sigma_{j}}{3},
	\frac{-6+2d+4\sum_{i=1}^d\sigma_{i}^2}{3d(d-1)}\Bigg\}<\min\Bigg\{\frac{1+2\sigma_{i}^2-|(4-d)\sigma_{i}|}{3(d-1)}\Bigg\},\:\: i,j\in \llbracket{1,d}\rrbracket,\: d<6.
\end{align} 
Due to the nonlinearity and coupling of each $\sigma_i$ ($i\in\llbracket1,d\rrbracket$) in Eq. (\ref{region-d>=3}), here  we take the   three-dimensional case as an example, and numerically plot the microscopic entropy stability region under different values of $\sigma_3$ in Fig. \ref{entropy-stable-fig}. From the figure, it can be found that the microscopic entropy stability region  of the  BGK-LB model (\ref{BGK-LB-Evolution-seq2})  becomes smaller as the value of $|\sigma_3|$ increases.\\ 

\noindent\textbf{Remark 6.} When the microscopic entropy stability  conditions in Eq. (\ref{entropy-stable-condition}) hold,  
we can obtain the following conclusions:
\begin{itemize}
\item The conservation of the total microscopic entropy  $\sum_{\mathbf{x}\in\Omega}\varpi(\mathbf{f})$ 

According to the work in \cite{bellotti2024-entropy}, we can conclude that $\sum_{\mathbf{x}\in\Omega}\varpi(\mathbf{f})$ in the computational domain $\Omega$ of the BGK-LB model (\ref{BGK-LB-Evolution-seq2}) is conserved when considering the periodic boundary condition. It is worth noting that the relaxation parameter $s = 2$ is necessary for the total microscopic entropy conservation. This is because the kinetic entropies $s_k$ in  Eq. (\ref{skfk}) are quadratic functions, and their respective minima under the conservation constraint are given at the EDFs $f_k^{eq}$. Thus, this is the reason why we adopt the BGK-LB model rather than the more general MRT-LB model.
\item  The  BGK-LB model (\ref{BGK-LB-Evolution-seq2}) is weighted $L^2$ stable 

According to the definition of the stability structure proposed in Ref. \cite{stability-structure}, for the developed BGK-LB model (\ref{BGK-LB-Evolution-seq2})  with the  microscopic entropy stability  conditions in Eq. (\ref{entropy-stable-condition}), i.e.,  $\bm{\omega}>\mathbf{0}$, being satisfied, when considering the periodic boundary condition, one can obtain  
\begin{align}\label{stability-structure}
	\Big	\|\sum_{\mathbf{x}\in\Omega}\mathbf{f}(\mathbf{x},t_n)\Big\|_{\mathbf{P}}\leq\Big\|\sum_{\mathbf{x}\in\Omega}\mathbf{f}(\mathbf{x},0)\Big\|_{\mathbf{P}},\:\: \forall n\geq 0,
\end{align}
 where the norm $\|\cdot\|_{\mathbf{P}}$ is defined as
  \begin{align}
	\|\mathbf{P}\mathbf{f}(\mathbf{x},t_n)\|^2=\sum_{k=0}^{2d^2+1}a_i^2f_k^2(\mathbf{x},t_n)
\end{align}where $\mathbf{P}=$\textbf{diag}$\big(a_1,a_2,\ldots,a_{2d^2+1}\big)$ with $a_k\in\mathbb{R}/\{0\}$. Eq. (\ref{stability-structure})  indicates that the norm of $\sum_{\mathbf{x}\in\Omega}\mathbf{f} $ does not increase during collision, i.e., the developed BGK-LB model (\ref{BGK-LB-Evolution-seq2}) for the L-HE (\ref{HE}) is weighted $L^2$ stable. In particular, the matrix $\mathbf{P}$ can be chosen as the identity matrix $\mathbf{I}_{2d^2+1}$, in this case, the  BGK-LB model (\ref{BGK-LB-Evolution-seq2}) is $L^2$ stable. However, we would like to point out that for the BGK-LB model (\ref{BGK-LB-Evolution-seq2}), the weighted $L^2$ stability presented in Eq. (\ref{stability-structure}) is with respect to the distribution functions    $f_k $ rather than the conserved quantity $\phi$ discussed in Part \ref{l2-stability} below.
\end{itemize}

\subsection{The $L^2$ stability analysis}\label{l2-stability}
In this part, we will discuss the $L^2$ stability of the BGK-LB model (\ref{BGK-LB-Evolution-seq2}) with the help of the von Neumann analysis. Without loss of generality, let us neglect the constant source term $R$. To this end, we first introduce the spatial shift operator on the spatial and temporal lattices $\mathcal{L}_{\mathbf{x}}$ and $\mathcal{T}$ as
\begin{align}\label{operator}
	T^{\mathbf{z}}_{\Delta t}h(\mathbf{x},t)=h(\mathbf{x}+\mathbf{z}\Delta t,t),\:\: \mathbf{x}\in\mathcal{L}_{\mathbf{x}},\: t\in\mathcal{T},
\end{align}
where vector $\mathbf{z} \in\mathbb{R}^d$, and $h: (\mathbb{R}^d,\mathbb{R})\rightarrow\mathbb{R}$ is an arbitrary smooth function. Based on this operator (\ref{operator}), one can rewrite the BGK-LB model (\ref{BGK-LB-Evolution-seq2}) in the moment space, i.e.,
 \begin{align}
 	\mathbf{m}(\mathbf{x},t)=\mathbf{MT}\mathbf{\tilde{M}}\Big(2\bm{\varepsilon}\otimes \mathbf{e}_1m_1-\mathbf{I}_{2d^2+1}\mathbf{m}\Big)(\mathbf{x},t),
 \end{align}
 where $\mathbf{m}=\mathbf{Mf}$, $m_1=m_1^{eq}=\phi$ is the conserved quantity, and  the shift matrix $\mathbf{T}$ is given by
 \begin{align}\label{moment-sapce}
 	\mathbf{T}={\makebox{\textbf{diag}}}\Big(T^{-\mathbf{c}_1}_{\Delta t},T^{-\mathbf{c}_2}_{\Delta t},\ldots,T^{-\mathbf{c}_{2d^2+1}}_{\Delta t}\Big).
 \end{align} 
Based on the works in Refs. \cite{Be2023,STERLING1996}, we take the discrete Fourier transform to Eq. (\ref{moment-sapce}) and can obtain the amplification matrix $\mathbf{G}$ of the BGK-LB model (\ref{BGK-LB-Evolution-seq2}) as
\begin{align}\label{G}
	\mathbf{G}=\mathbf{M\hat{T}}\mathbf{\tilde{M}}\Big(2\bm{\varepsilon}\otimes \mathbf{e}_1-\mathbf{I}_{2d^2+1}\Big),
\end{align}
where the shift matrix $\mathbf{T}$ in the Fourier space denoted by $\mathbf{\hat{T}}$ is expressed as
\begin{align}
	\mathbf{\hat{T}}={\makebox{\textbf{diag}}}\Big(e^{-\mathbf{c}_1\cdot\bm{\xi}/c},e^{-\mathbf{c}_2\cdot\bm{\xi}/c},\ldots,e^{-\mathbf{c}_{2d^2+1}\cdot\bm{\xi}/c}\Big),\:\: \bm{\xi}\in\big[0,2\pi\big)^d.
\end{align}
And the BGK-LB model (\ref{BGK-LB-Evolution-seq2}) is $L^2$ stable as long as the minimal polynomial of the amplification matrix $\mathbf{G}$ in Eq. (\ref{G}) is a von Neumann polynomial \cite{bellotti2023-weakly}. One can find that, for the D$d$Q($2d^2+1$) lattice structure, the amplification matrix $\mathbf{G}\in\mathbb{R}^{(2d^2+1)\times (2d^2+1)}$, which means that the explicit $L^2$ stability condition of the BGK-LB model (\ref{BGK-LB-Evolution-seq2})  for the high-dimensional case is difficult to analyze. The only theoretical work to analyze the $L^2$ stability is for the one-dimensional case \cite{bellotti2023-weakly} , where the authors have proved that the BGK-LB model (\ref{BGK-LB-Evolution-seq2}) is $L^2$ stable as long as $|\sigma_1|<1/2$, which is consistent with the microscopic entropy stability condition presented in Remark 4 of our work. It should be noted that for the case of $|\sigma_1|=1/2$, we can prove that the BGK-LB model (\ref{BGK-LB-Evolution-seq2}) is still $L^2$ stable (see \ref{app-sigma=1/2} for details). From the work in Ref. \cite{bellotti2023-weakly}, we can find that the $L^2$ stability analysis of the BGK-LB model (\ref{BGK-LB-Evolution-seq2})  is divided into two cases, i.e., the case of  $\xi=0$ and  the case of  $\xi\neq0$. Similarly, for the high-dimensional case addressed in our work, when  $\bm{\xi}=\mathbf{0}$, the  amplification matrix  $\mathbf{G}$ in Eq. (\ref{G}) reduces
\begin{align}\label{G-3}
	\mathbf{G}=\Big(2\bm{\varepsilon}\otimes \mathbf{e}_1-\mathbf{I}_{2d^2+1}\Big),
\end{align}
and it is easy to validate that $\mathbf{G}^2=\mathbf{I}_{2d^2+1}$. Thus, the polynomial $p(\lambda)=\lambda^2-1$ is the minimal polynomial of the  amplification matrix  $\mathbf{G}$, and the two roots of the  polynomial $p(\lambda)=\lambda^2-1$ lying on the  unit circle are different. Thus, from the $L^2$ stability definition in the framework of the LB method,  we know  that for the case of $\bm{\xi}=\mathbf{0}$,  the BGK-LB model (\ref{BGK-LB-Evolution-seq2}) is stable without the constraint on $\sigma_i$ ($i\in\llbracket 1,d\rrbracket$). 

\begin{figure} 
		\centering   
			\includegraphics[width=0.35\textwidth]{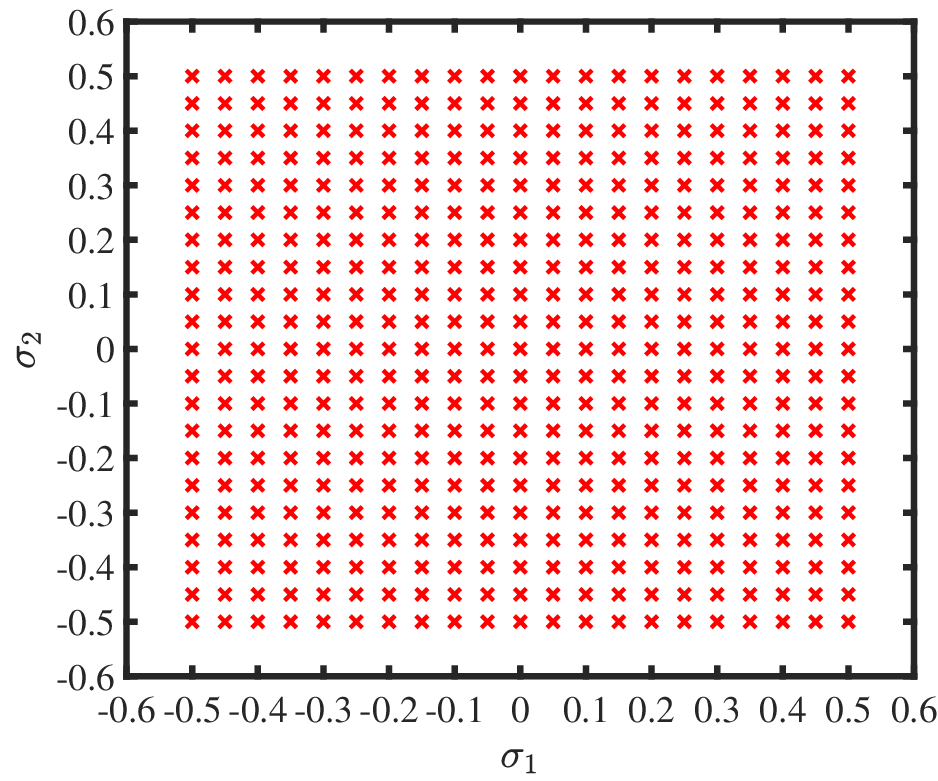} 
		\caption{The $L^2$ stability region  of the BGK-LB model   for the two-dimensional L-HE  ($\bm{\xi}\neq \mathbf{0}$). } 
		\label{l2-stable-2d}  
\end{figure}

\begin{figure} 
	\centering
		\subfloat[$\sigma_3=0$]    
		{
			\includegraphics[width=0.3\textwidth]{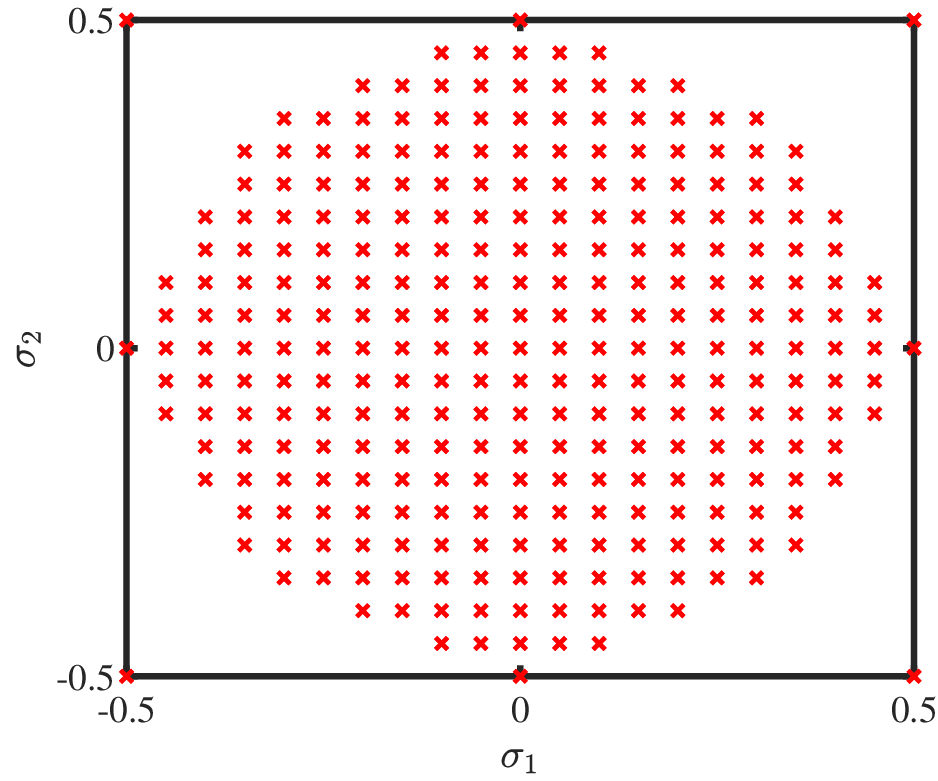} 
		}       
		\subfloat[$\sigma_3=0.1$]    
		{
			\includegraphics[width=0.3\textwidth]{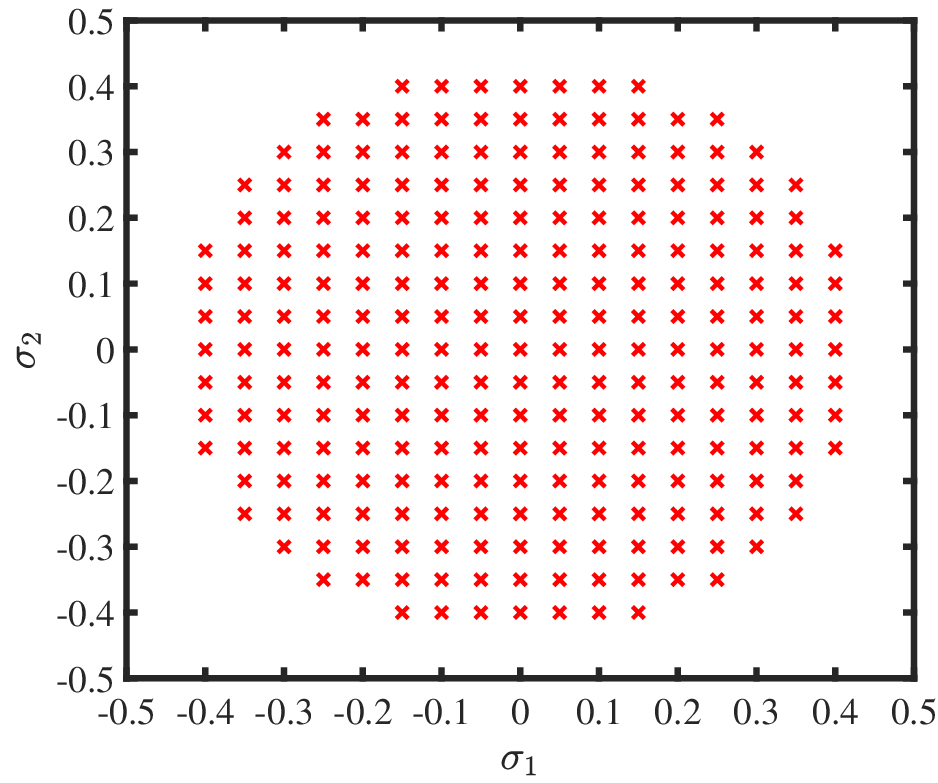} 
		}
	 
		\subfloat[$\sigma_3=0.2$]    
		{
			\includegraphics[width=0.3\textwidth]{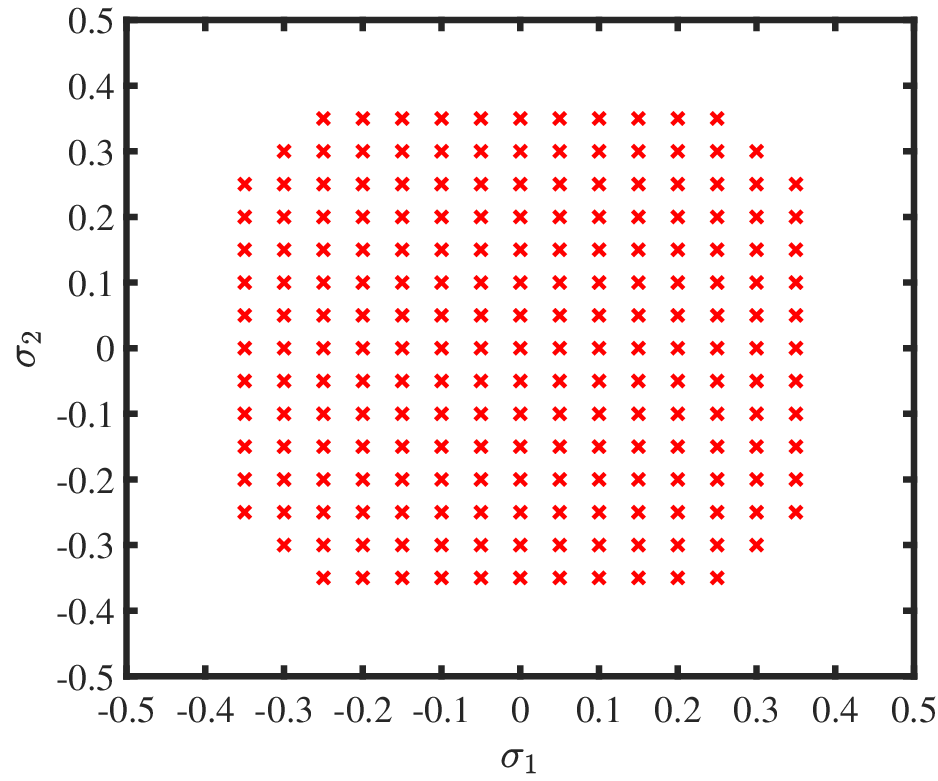} 
		}
		\subfloat[$\sigma_3=0.4$]    
		{
			\includegraphics[width=0.3\textwidth]{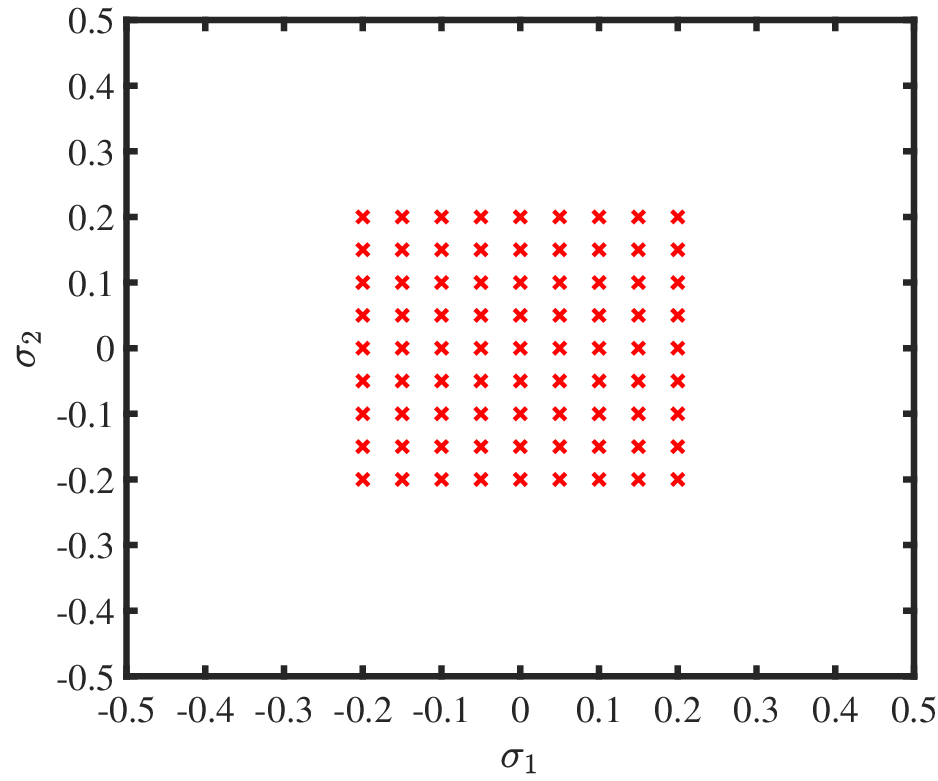} 
		}
		\caption{The $L^2$ stability  region of the BGK-LB model for the three-dimensional L-HE ($\bm{\xi}\neq \mathbf{0}$). } 
		\label{l2-stable-3d}  
\end{figure}

For the case of $\bm{\xi}\neq\mathbf{0}$, the theoretical analysis becomes extremely complicated. Here, we numerically plot the $L^2$ stability regions of the two- and three-dimensional cases in Figs. \ref{l2-stable-2d} and \ref{l2-stable-3d}, respectively, where $g$ in Eq. (\ref{meq-4rd}) is taken as 1/2 to determine $\dot{m}^{eq}_{x_i^2x_j^2}$. From these figures, according to  Remark 4 and the numerical microscopic entropy stability regions shown in Fig. \ref{entropy-stable-fig}, one can observed that the numerical $L^2$ and microscopic entropy stability regions are identical, expect at the points where the kinetic entropies $s_k$ are not defined, i.e., $\omega_k=0$. Thus, from a numerical perspective, it can be concluded that for the the BGK-LB model (\ref{BGK-LB-Evolution-seq2}), under the premise of $\omega_k\neq 0$, the $L^2$ stability condition, which is difficult to analyze theoretically, is same as the microscopic entropy stability condition, which indicates that when preforming numerical simulations,  it is reasonable to determine the parameters in the BGK-LB model  (\ref{BGK-LB-Evolution-seq2})   based on the microscopic entropy stability.

\section{Numerical results and discussion}\label{Numer}
In this section, we will perform some simulations to test the present BGK-LB model (\ref{BGK-LB-Evolution-seq2}) for the L-HE (\ref{HE}). And the lattice velocity is taken as $c=1$, the EDFs $f_k^{eq}$ are determined from Eqs. (\ref{meq}) and (\ref{meq-4rd}), in particular, without loss of generality, we take $g=1/2$ in Eq. (\ref{meq-4rd}). To measure the accuracy of the BGK-LB model  (\ref{BGK-LB-Evolution-seq2}),  the following root-mean-square error (RMSE) is adopted \cite{Kruger2017},
\begin{align}\label{rmse}
	{\makebox{RMSE}}=\sqrt{\frac{\sum_{i_1=1}^{\overline{N}_1}\ldots\sum_{i_d=1}^{\overline{N}_d}\Big(\phi(\mathbf{x}_{i_1,i_2,\ldots,i_d},t_n)-\phi^{\star}(\mathbf{x}_{i_1,i_2,\ldots,i_d},t_n)\Big)^2}{\prod_{l=1}^d\overline{N}_l}},
\end{align}
where $\prod_{l=1}^d\overline{N}_l$ is the number of the lattice points, $\overline{N}_l=L_l/\Delta x+1-\vartheta$ is the  number of the lattice points in the $l_{th}$ spatial direction with $L_l=b_l-a_l$ being the characteristic length, $\mathbf{x}_{i_1,i_2,\ldots,i_d}$ is the lattice point with $x_{i_l}=a_l+[i_l-1/(\vartheta+1)]\Delta x$. $\phi$ and $\phi^{\star}$ are the numerical and analytical solutions respectively. Based on the definition of RMSE (\ref{rmse}), once can estimate the convergence rate (CR) of the BGK-LB model  (\ref{BGK-LB-Evolution-seq2}).\\ 
\noindent\textbf{Example 1.} We first consider the two-dimensional L-HE (\ref{HE}) with the following initial and Dirichlet boundary conditions:
\begin{subequations}
	\begin{align}
		&\phi(x_1,x_2,0)=\exp\big[-(x_1+x_2)\big],\:\: (x_1,x_2)\in[0,1]^2,\\
		&\phi(x_1,0,t)=\exp (-x_1)+Rt,\:\: \phi(x_1,1,t)=\exp\big[-(1+x_1)\big]+Rt,\:\: x_1\in[0,1],\: t>0,\\
		&\phi(0,x_2,t)=\exp (-x_2)+Rt,\:\: \phi(1,x_2,t)=\exp\big[-(1+x_2)\big]+Rt,\:\: x_2\in[0,1],\: t>0, 
			\end{align}
	\end{subequations}
and obtain the analytical solution as
\begin{align}
	\phi^{\star}(x_1,x_2,t)=\exp\big[-(x_1+x_2)\big]+Rt,\:\: x\in[0,1]^2,\: t>0,
\end{align}
where we take $R=5$ in the following simulations.

\begin{figure} 
	\centering
	\subfloat[Or.-Boun.: $O(\Delta t)$]  
	{
		\includegraphics[width=0.308\textwidth]{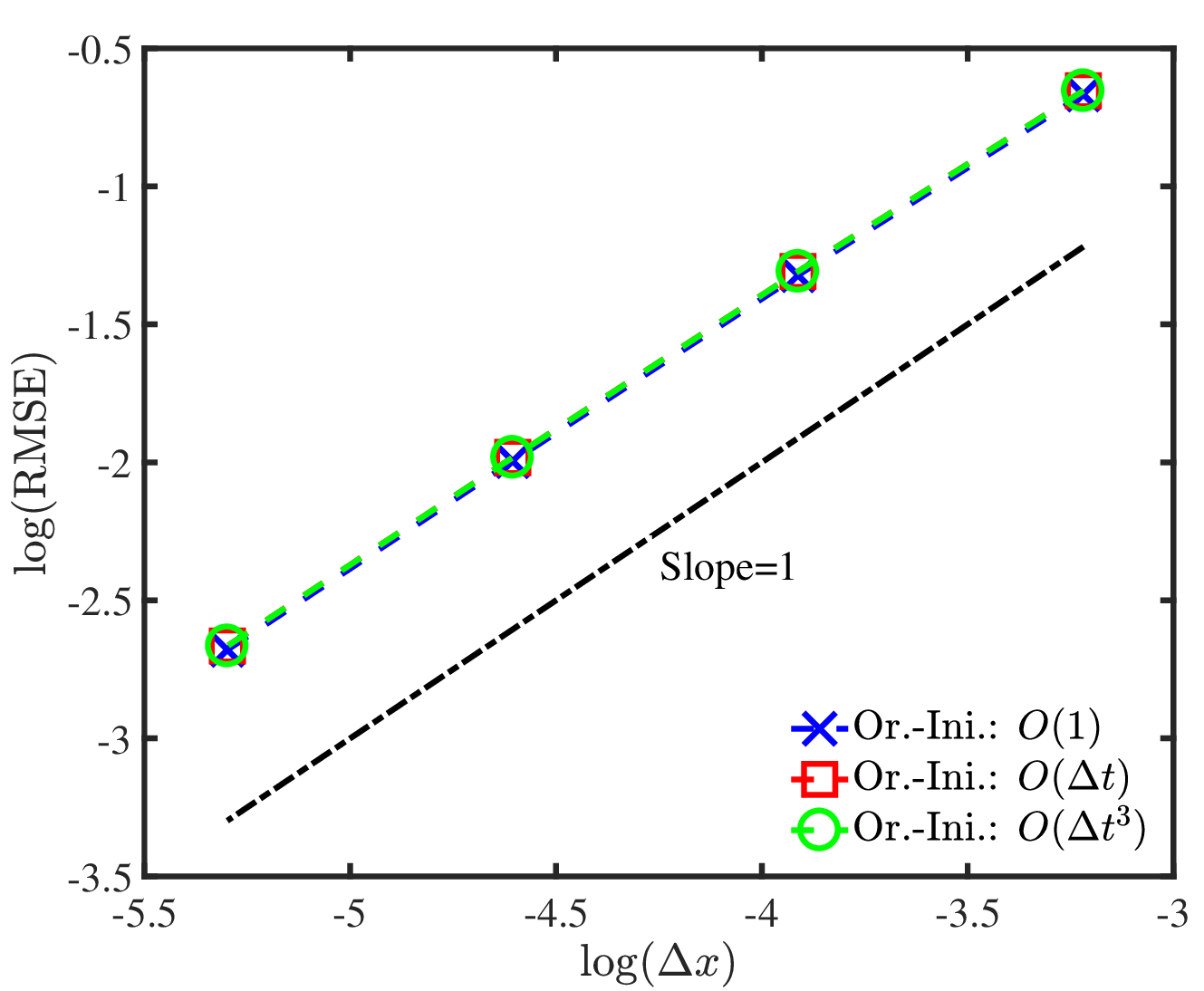} 
	} 
	\subfloat[Or.-Boun.: $O(\Delta t^3)$]    
	{
		\includegraphics[width=0.3\textwidth]{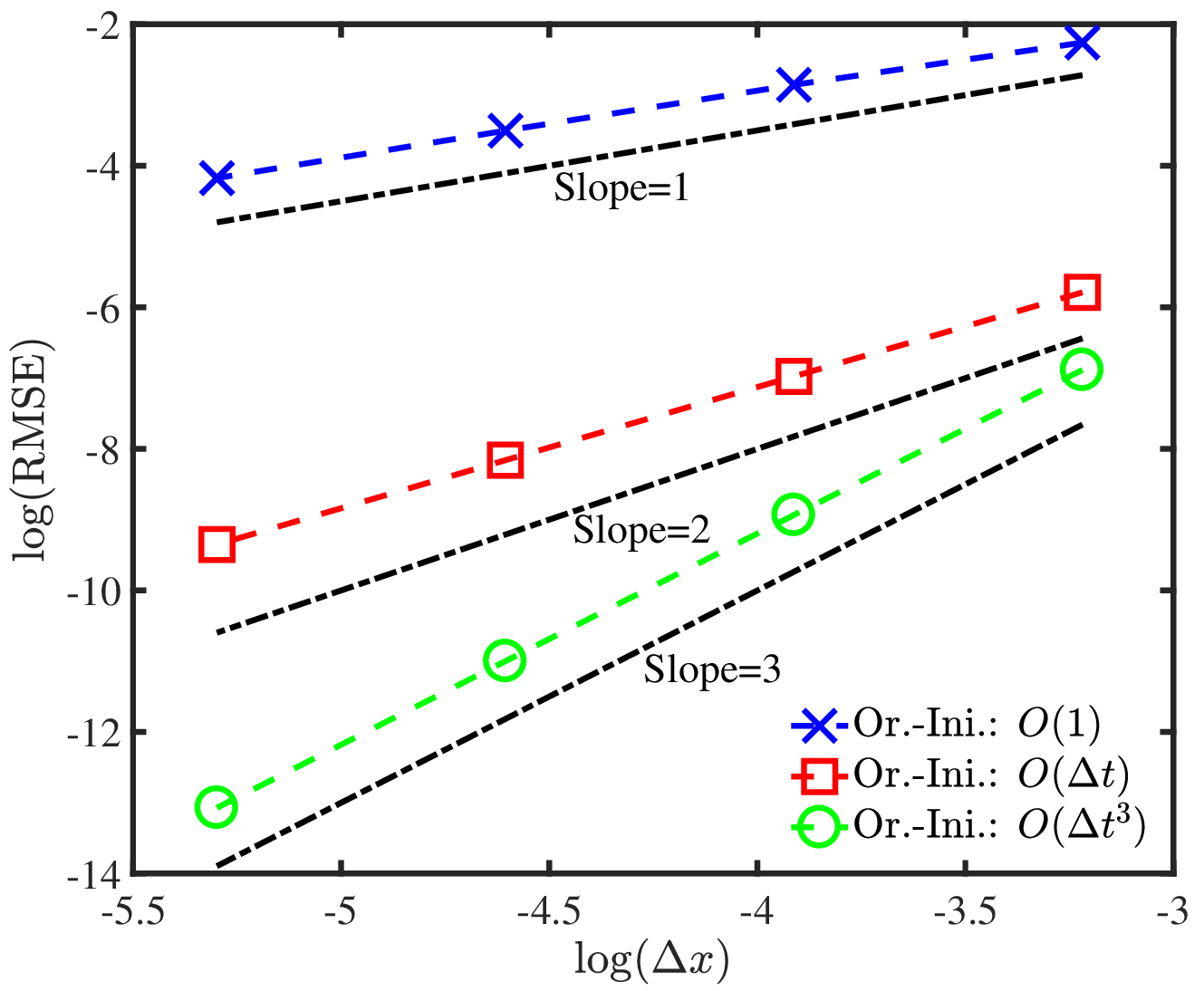} 
	} 
		\subfloat[Or.-Boun.: $O(\Delta t^4)$]    
	{
		\includegraphics[width=0.3\textwidth]{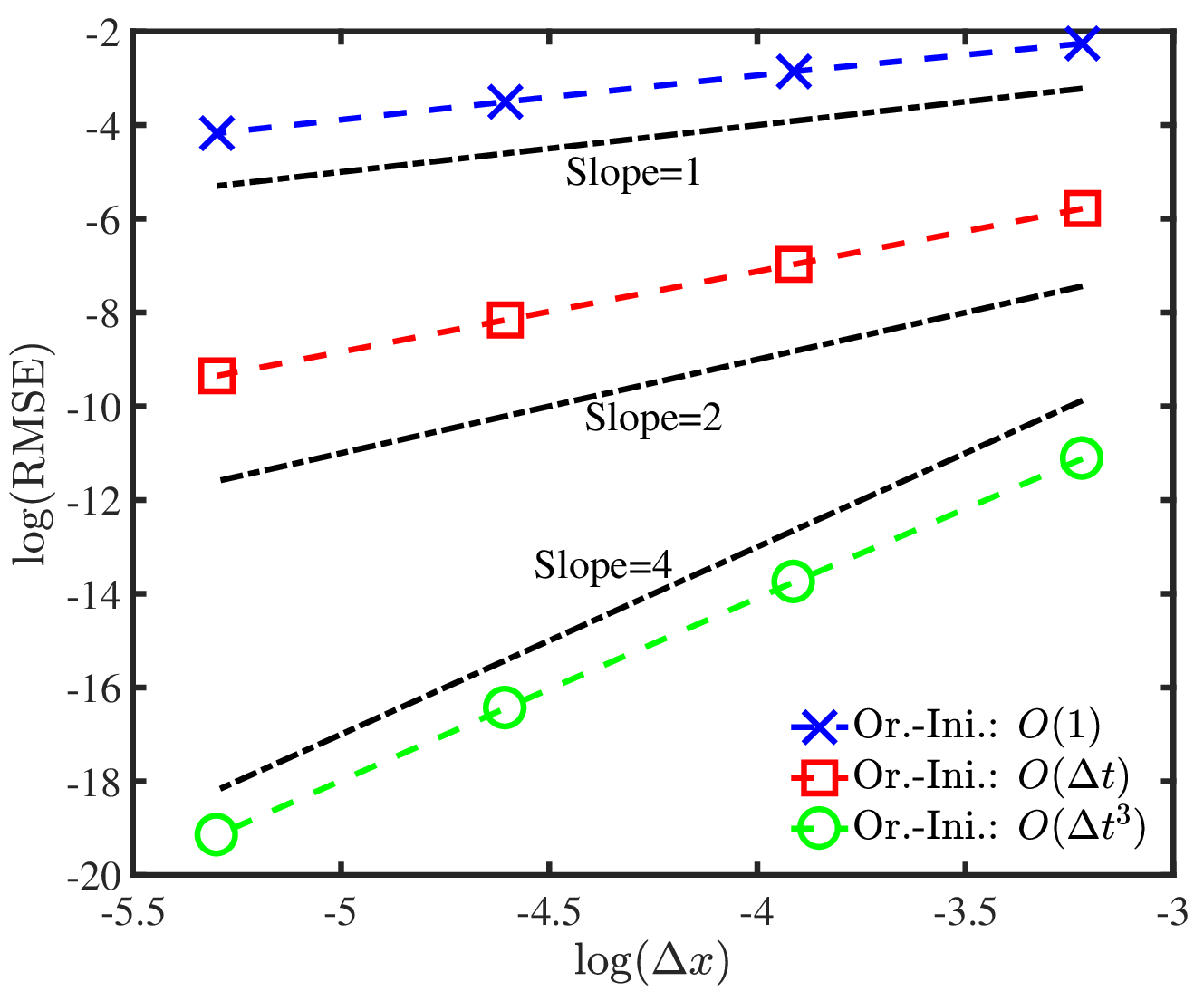} 
	} 
	\caption{The convergence rate of the  BGK-LB model  under  the transport velocity $\mathbf{u}=(0.1,0.1)$. \protect\footnotemark[4] }
	\label{fig-full-u0.10.1}
\end{figure}

\begin{figure} 
	\centering
	\subfloat[Or.-Boun.: $O(\Delta t)$]    
	{
		\includegraphics[width=0.3\textwidth]{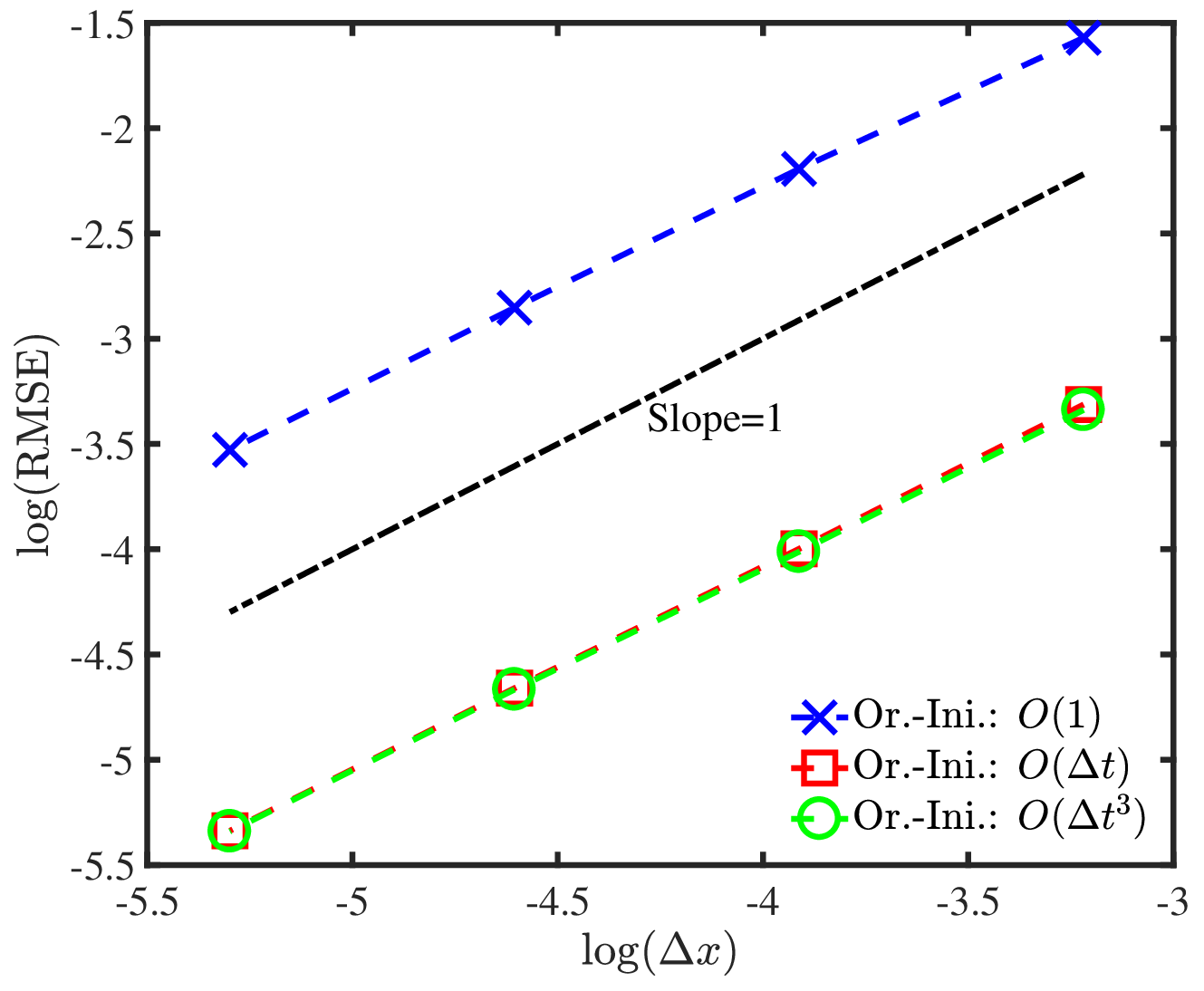} 
	} 
	\subfloat[Or.-Boun.: $O(\Delta t^3)$]    
	{
		\includegraphics[width=0.3\textwidth]{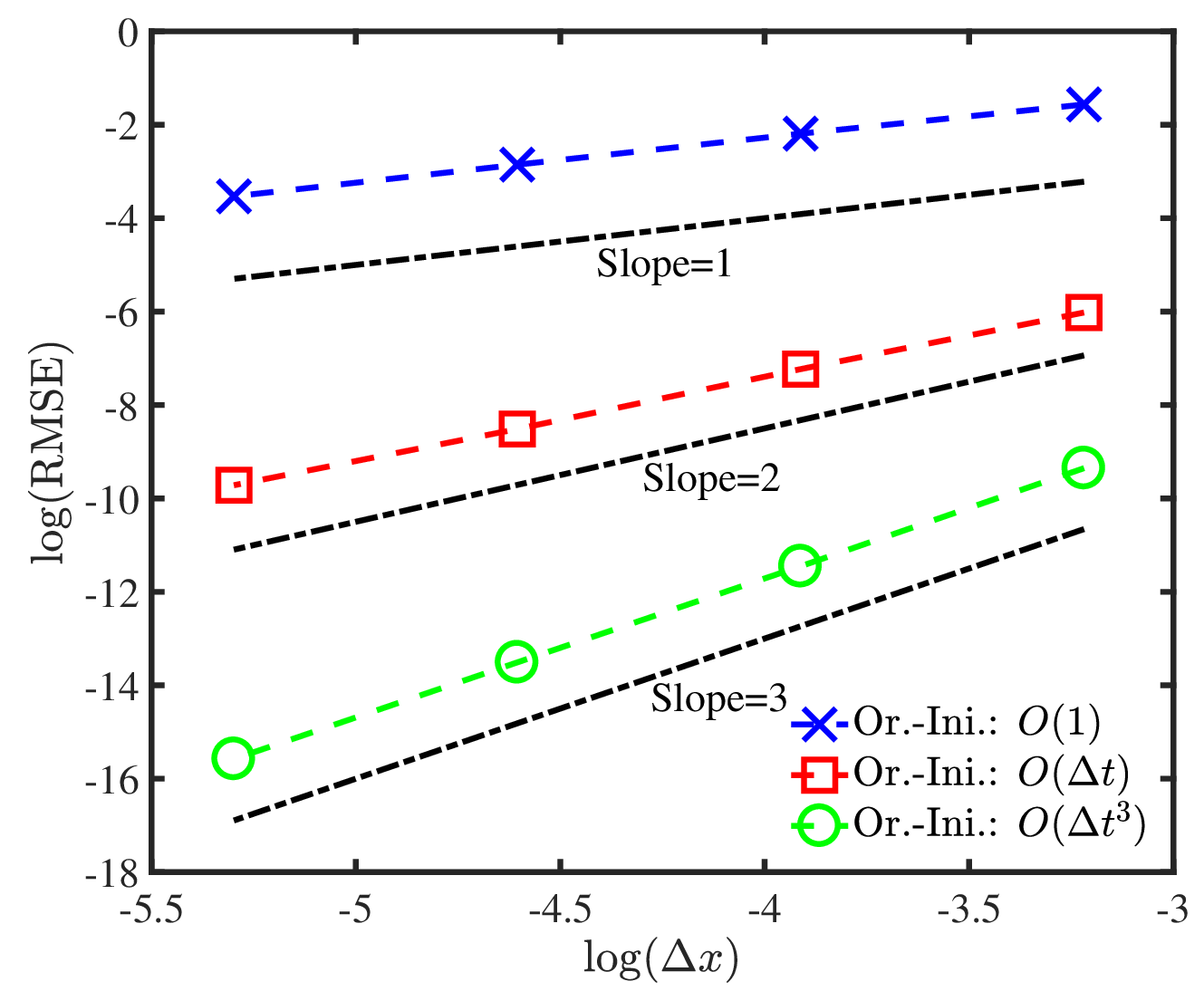} 
	} 
	\subfloat[Or.-Boun.: $O(\Delta t^4)$]    
	{
		\includegraphics[width=0.3\textwidth]{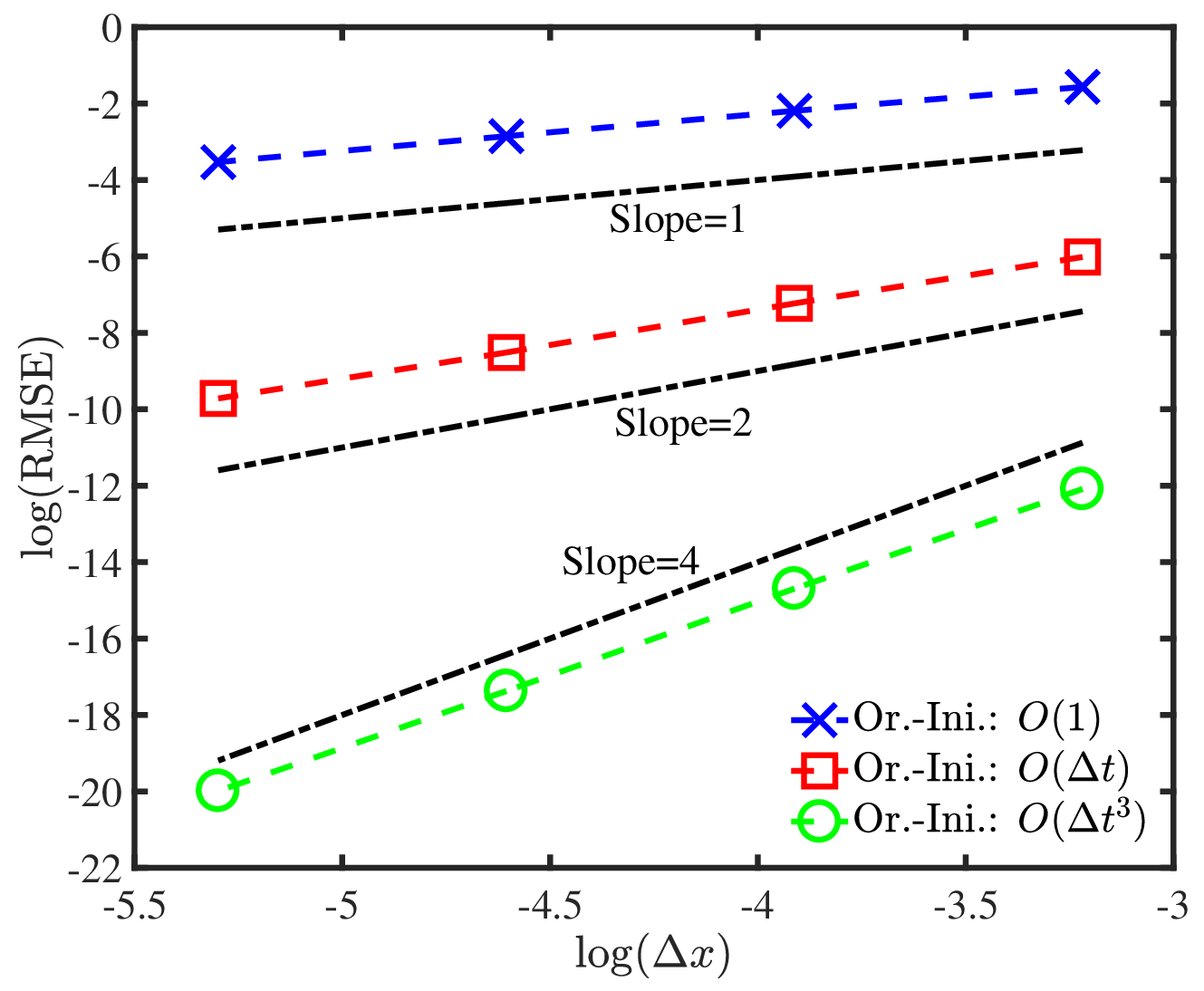} 
	} 
	\caption{The convergence rate of the BGK-LB model  under the transport velocity $\mathbf{u}=(0.01,0.01)$. } 
	\label{fig-full-u0.010.01}
\end{figure}

\begin{figure} 
	\centering
	\subfloat[Or.-Boun.: $O(\Delta t)$]
	{
		\includegraphics[width=0.35\textwidth]{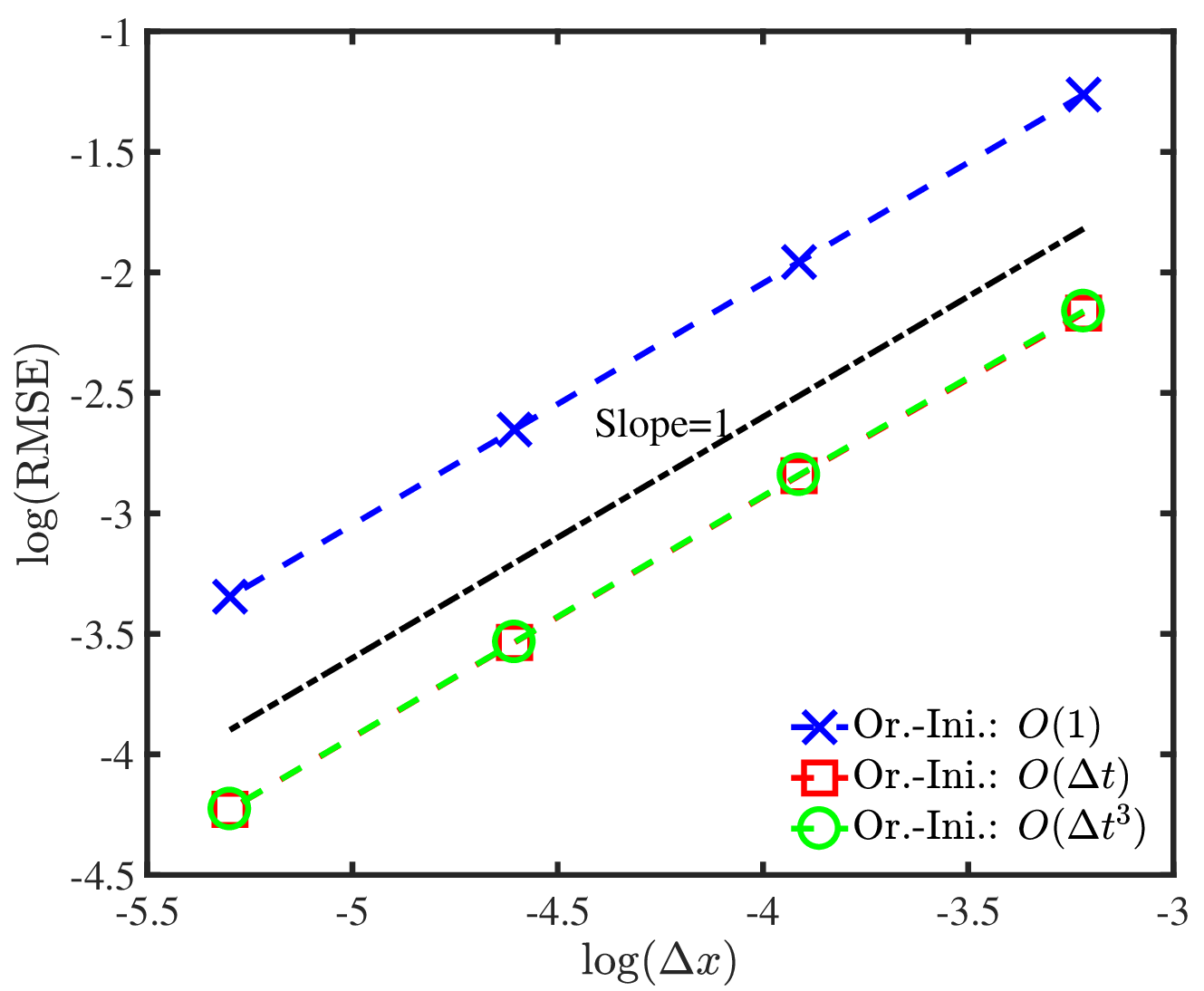} 
	} 
	\subfloat[Or.-Boun.: $O(\Delta t^2)$]    
	{
		\includegraphics[width=0.35\textwidth]{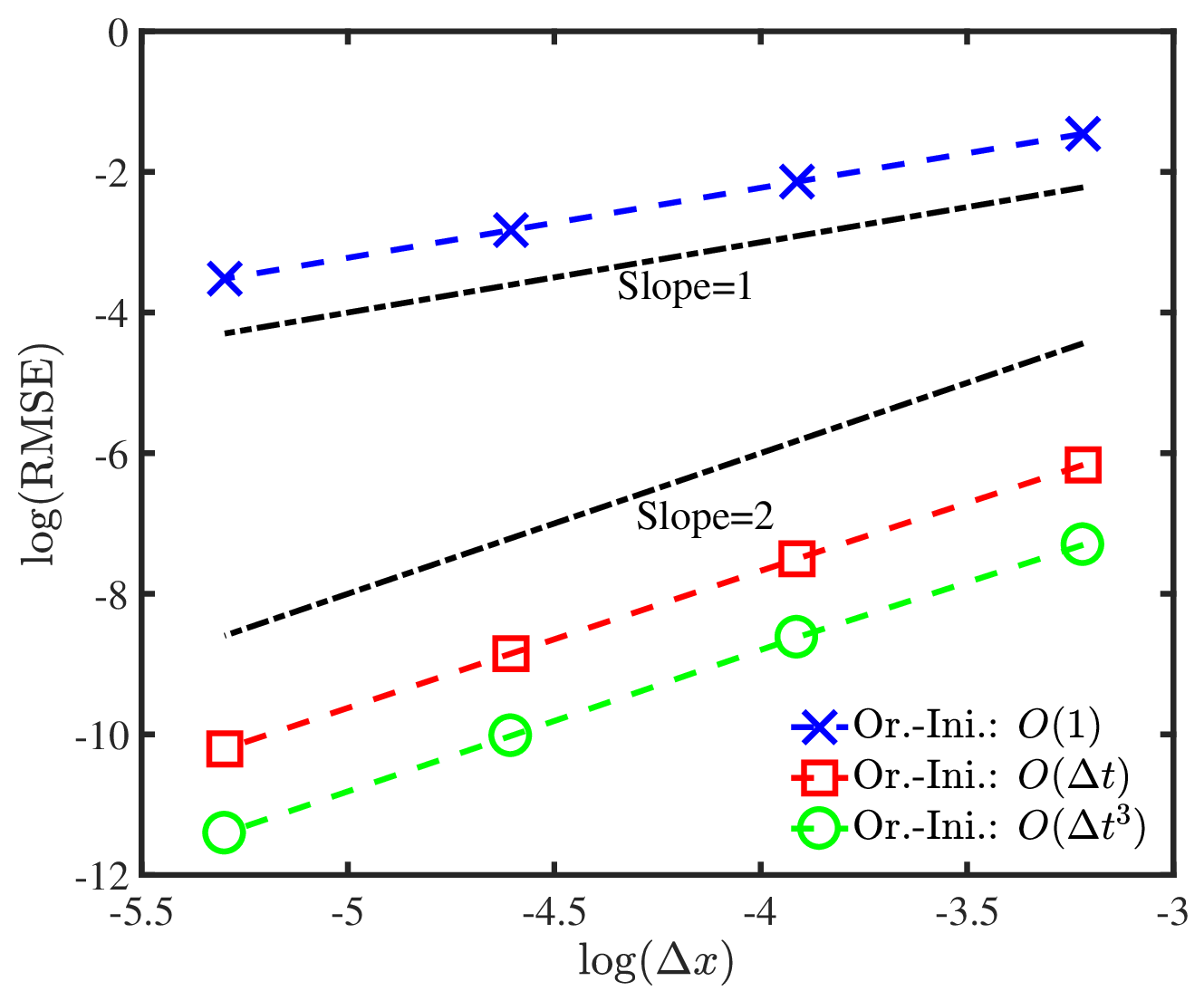} 
	} 
	
	\subfloat[Or.-Boun.: $O(\Delta t^3)$]    
	{
		\includegraphics[width=0.35\textwidth]{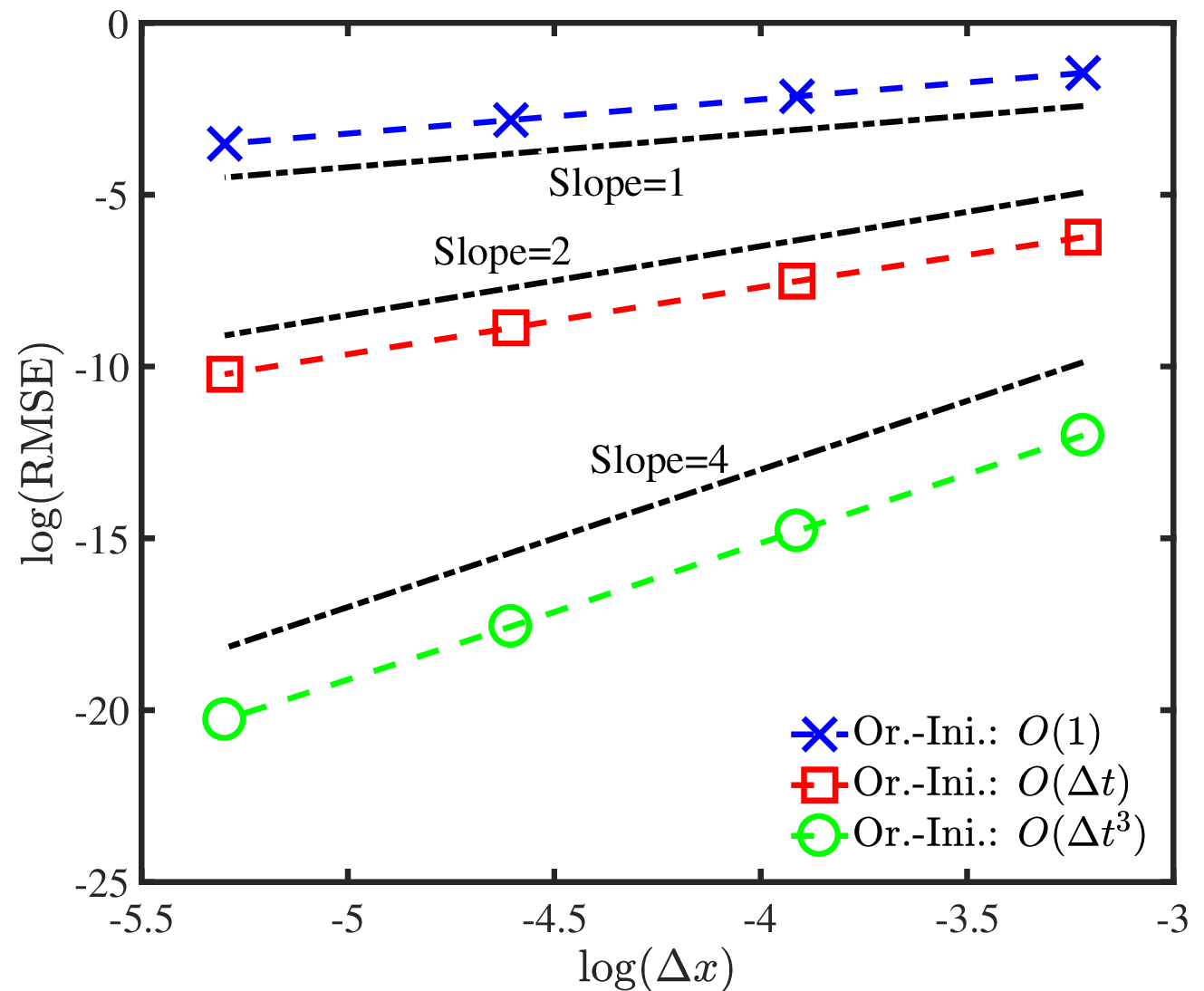} 
	} 
	\subfloat[Or.-Boun.: $O(\Delta t^4)$]    
	{
		\includegraphics[width=0.35\textwidth]{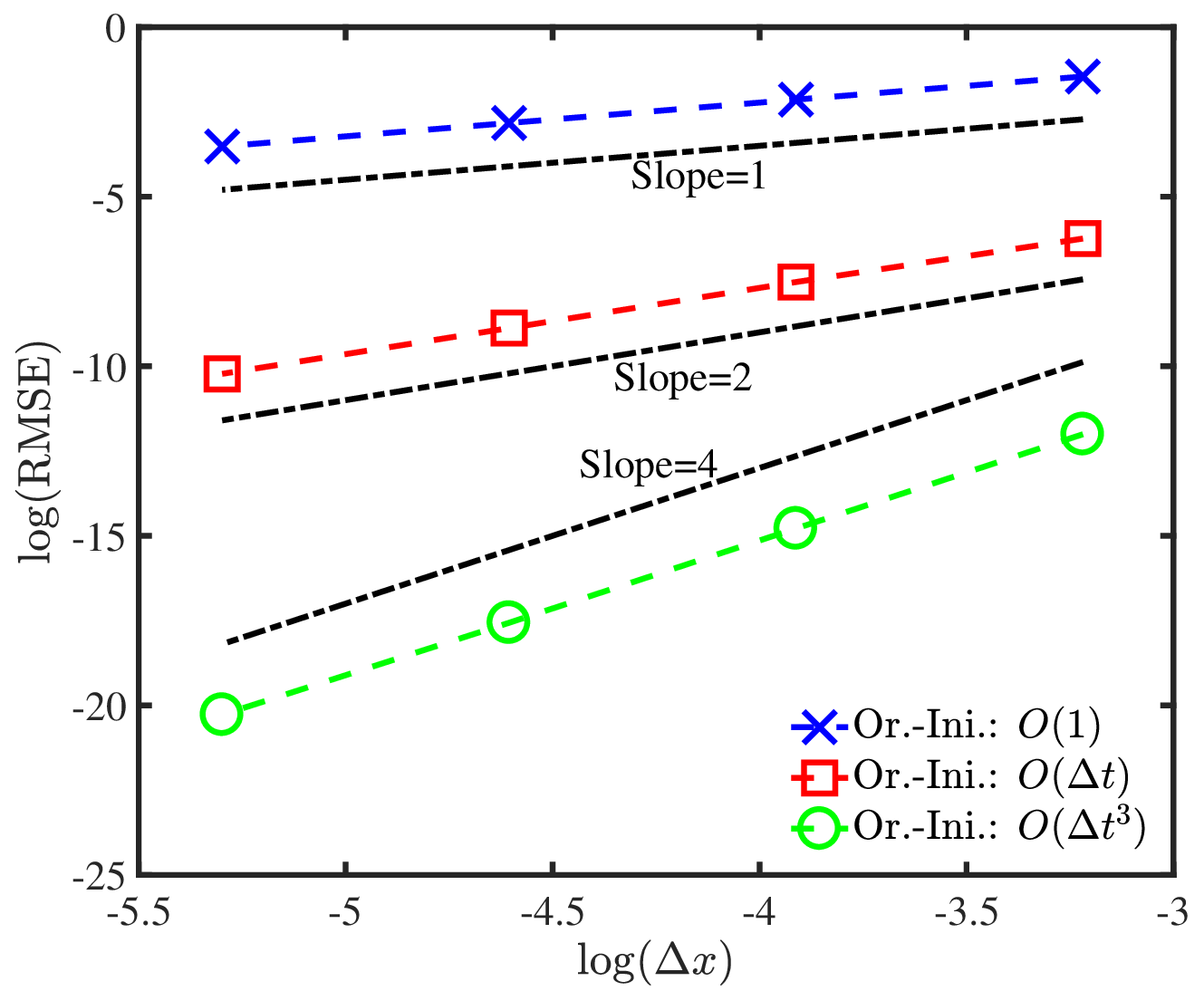} 
	} 
	\caption{The convergence rate  of the  BGK-LB model  under   the transport velocity $\mathbf{u}=(0.01,0.01)$. } 
	\label{fig-half-u0.010.01}
\end{figure}

In this example, the full-way and half-way boundary schemes for the Dirichlet condition are adopted. Firstly, we perform a number of simulations with different accuracy orders of the initialization and boundary schemes under some specified values of the transport velocity $\mathbf{u}=(u_1,u_2)$, and calculate the RMSEs and CRs under  different lattice sizes. As seen from Figs. \ref{fig-full-u0.10.1}, \ref{fig-full-u0.010.01},  and \ref{fig-half-u0.010.01}, where the simulations are suspended at the time $T=10.0$, the CR of the BGK-LB model (\ref{BGK-LB-Evolution-seq2}) is  closely related to the accuracy orders of the initialization and boundary schemes. In particular, for the full-way boundary scheme,  one can conclude  
\begin{align}\label{cr-estimation}
	Z= \min\{K+1, P\}, 
\end{align} 
where $Z$ is the order of the overall CR of the BGK-LB model (\ref{BGK-LB-Evolution-seq2}). $K$ and $P$ are the accuracy orders of the initialization and full-way boundary schemes, respectively. And for the half-way boundary scheme in Eq. (\ref{bun2}), one can observe the same relation as shown in Eq. (\ref{cr-estimation}) when $P=$1, 2, and 4 in Figs. \ref{fig-half-u0.010.01}(a), \ref{fig-half-u0.010.01}(b), and \ref{fig-half-u0.010.01}(d), however, for the case of $P = 3$ shown in  Fig. \ref{fig-half-u0.010.01}(c),  the BGK-LB model (\ref{BGK-LB-Evolution-seq2}) can reach fourth-order convergence when the third-order initialization scheme is adopted. In fact, according to the expression of $\omega_k$ given in Eq. (\ref{omega_k-2}), we would like to point out that  the terms of the fourth-order half-way boundary scheme in Eq. (\ref{bun-2}) at the order of $O(\Delta t^3)$ are actually at the order of $O(\mathbf{u}\Delta t^3)$. This means that when the transport velocity is small enough, compared to the bulk error of the BGK-LB model (\ref{BGK-LB-Evolution-seq2}), these terms can be omitted in the fourth-order half-way boundary scheme, i,e., the behavior of the third-order and fourth-order half-way boundary schemes can be considered almost the same. Thus, under the small transport velocity $\mathbf{u}=(0.01,0.01)$, the BGK-LB model (\ref{BGK-LB-Evolution-seq2}) can also be fourth-order accurate when using the third-order half-way boundary and initialization schemes.

Then, in order to compare the fourth-order full-way ($\vartheta = 1$) and half-way ($\vartheta = 0$) boundary schemes, we conduct some simulations under different values of the transport velocity that satisfy the microscopic entropy stability and $L^2$ stability conditions of the BGK-LB model (\ref{BGK-LB-Evolution-seq2}), and present the results in Table \ref{full-half}. From this table, for the case where the transport velocity $|u_{1,2}|>0.01$, although the microscopic entropy stability and $L^2$ stability hold, the numerical solution obtained from the BGK-LB model (\ref{BGK-LB-Evolution-seq2}), where the half-way boundary scheme for the Dirichlet condition is adopted, cannot converge to the analytical solution.  This means that the full-way boundary scheme, which is consistent with the bulk BGK-LB model (\ref{BGK-LB-Evolution-seq2}), is more stable than the half-way boundary scheme.
\footnotetext[4]{For brevity, we shall consistently write "the accuracy  order of the initialization scheme" as "Or.-Ini." and  "the accuracy order of the boundary scheme" as "Or.-Boun.".} 

 \begin{table} 
	\caption{The comparison in the RMSEs and CRs of the BGK-LB model with the full-way and half-way   boundary schemes   ($\Delta x=1/10$).} \label{full-half}
	\begin{tabular*}{\textwidth}{@{\extracolsep\fill}cccccccc}	\hline \hline
		$\mathbf{u}$&$\vartheta$&RMSE$_{\Delta x,\Delta t}$&RMSE$_{\Delta x/2,\Delta t/2}$&RMSE$_{\Delta x/4,\Delta t/4}$&RMSE$_{\Delta x/8,\Delta t/8}$&RMSE$_{\Delta x/16,\Delta t/16}$&CR\\ 
		\hline
		\multirow{2}{*}{$\big(0.4,0.4)$}&1&4.6014$\times 10^{-4}$&3.4988$\times 10^{-5}$&2.3994$\times 10^{-6}$&1.5688$\times 10^{-7}$&1.0063$\times 10^{-8}$&$\sim$3.8702 \\ 
		\cline{2-8}
		&0&16.7187&3.4275$\times 10^{+8}$&6.0287$\times 10^{+24}$&3.7709$\times 10^{+58}$&2.6896$\times 10^{+127}$&-\\ 
		\cline{1-8}
		\multirow{2}{*}{$\big(0.2,0.1)$}&1&2.2060$\times 10^{-5}$&1.6464$\times 10^{-6}$&1.1200$\times 10^{-7}$&7.3275$\times 10^{-9}$&4.6940$\times 10^{-10}$&$\sim$3.8801 \\ 
		\cline{2-8}
		&0&0.0011&2.7736&1.5666$\times 10^{+8}$&4.2672$\times 10^{+25}$&7.0852$\times 10^{+61}$&-\\ 
		\cline{1-8}
		\multirow{2}{*}{$\big(0.1,0.1)$}&1&1.4832$\times 10^{-5}$&1.0649$\times 10^{-6}$&7.2253$\times 10^{-8}$&4.8008$\times 10^{-9}$&3.3198$\times 10^{-10}$&$\sim$3.8618\\
		\cline{2-8}
		&0&3.4340$\times 10^{-5}$&8.1972$\times 10^{-5}$&3.3869$\times 10^{-1}$&2.1633$\times 10^{+8}$&4.0130$\times 10^{+27}$&- \\ 
		\cline{1-8}
		\multirow{2}{*}{$\big(0.05,0.05)$}&1&9.0341$\times 10^{-6}$&6.4447$\times 10^{-7}$&4.4311$\times 10^{-8}$&3.0765$\times 10^{-9}$&2.2320$\times 10^{-10}$&$\sim$3.8262\\
		\cline{2-8}
		&0&1.2374$\times 10^{-5}$&7.8645$\times 10^{-7}$&1.6739$\times 10^{-6}$&0.0102&5.5674$\times 10^{+6}$&- \\ 
		\cline{1-8}
		\multirow{2}{*}{$\big(0.01,0.01)$}&1&5.6359$\times 10^{-5}$&4.1543$\times 10^{-7}$&2.8506$\times 10^{-8}$&2.0876$\times 10^{-9}$&1.4790$\times 10^{-10}$&$\sim$3.8044\\
		\cline{2-8}
		&0&6.3881$\times 10^{-6}$&3.9229$\times 10^{-7}$&2.4429$\times 10^{-8}$&1.5211$\times 10^{-9}$&1.6180$\times 10^{-10}$&$\sim$3.8172\\
			\hline \hline
	\end{tabular*}
\end{table}
\noindent\textbf{Example 2.} We then consider the two and three-dimensional cases with the periodic boundary condition and the following initial condition:
	\begin{align}
		&\phi(x_1,\ldots, x_d,0)=\frac{1}{2}-\frac{1}{2}\tanh\frac{Rr-\sqrt{\sum_{i=1}^dx_i^2}}{2W},\:\: x\in\Big[\frac{d-4}{2},\frac{4-d}{2}\Big]^d,\:\: d\in\llbracket 2,3\rrbracket,
	\end{align}
and the analytical solution is
\begin{align}
	\phi^{\star}(x_1,\ldots,x_d,t)=\frac{1}{2}-\frac{1}{2}\tanh\frac{Rr-\sqrt{\sum_{i=1}^d\big(x_i-u_it\big)^2}}{2W},\:\: x\in\Big[\frac{d-4}{2},\frac{4-d}{2}\Big]^d,\: t>0.
\end{align}

\begin{figure} 
		\centering       
			\subfloat[$d=2$]    
			{
				\includegraphics[width=0.35\textwidth]{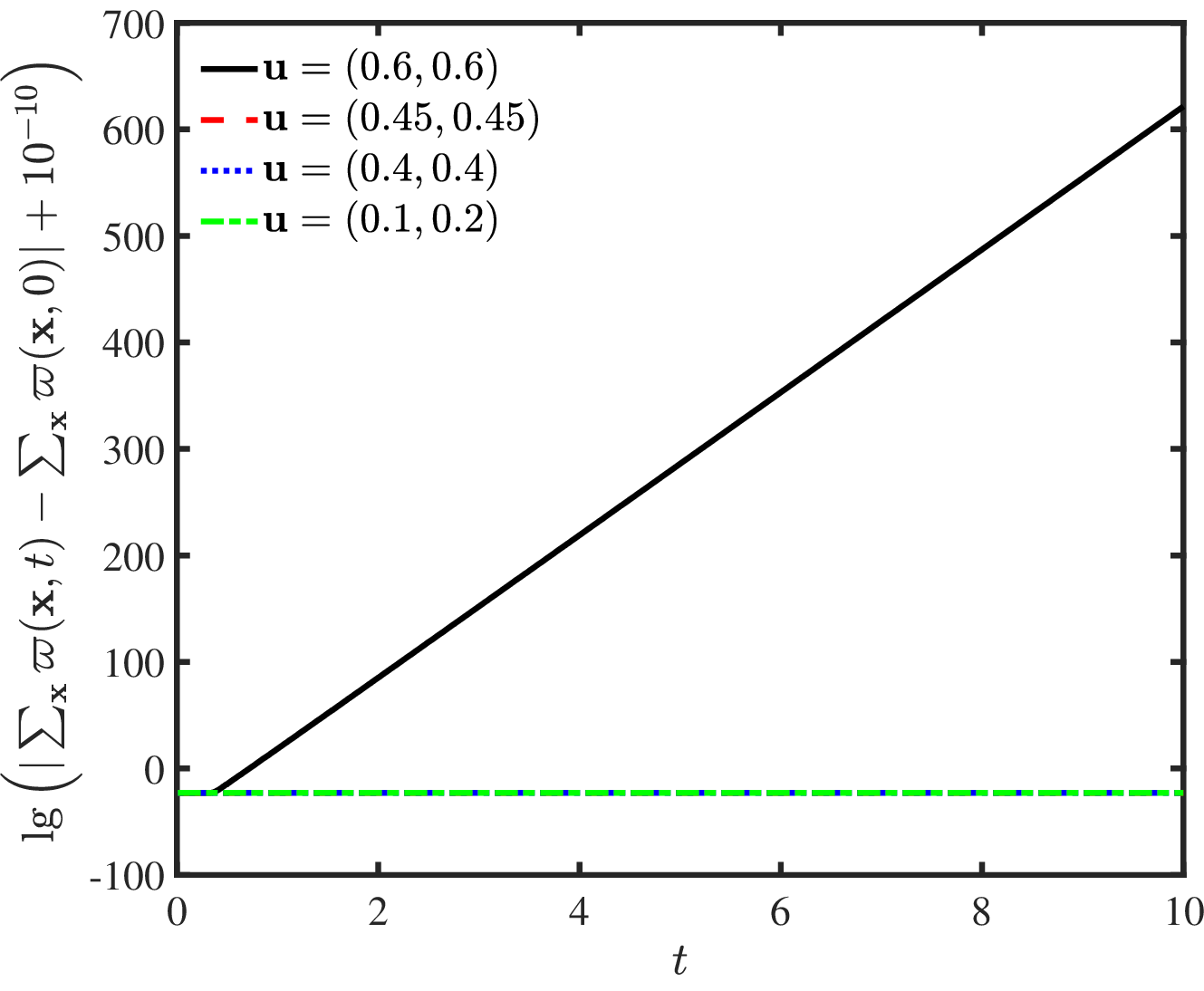} 
			}
			\subfloat[$d=3$]
			{
				\includegraphics[width=0.35\textwidth]{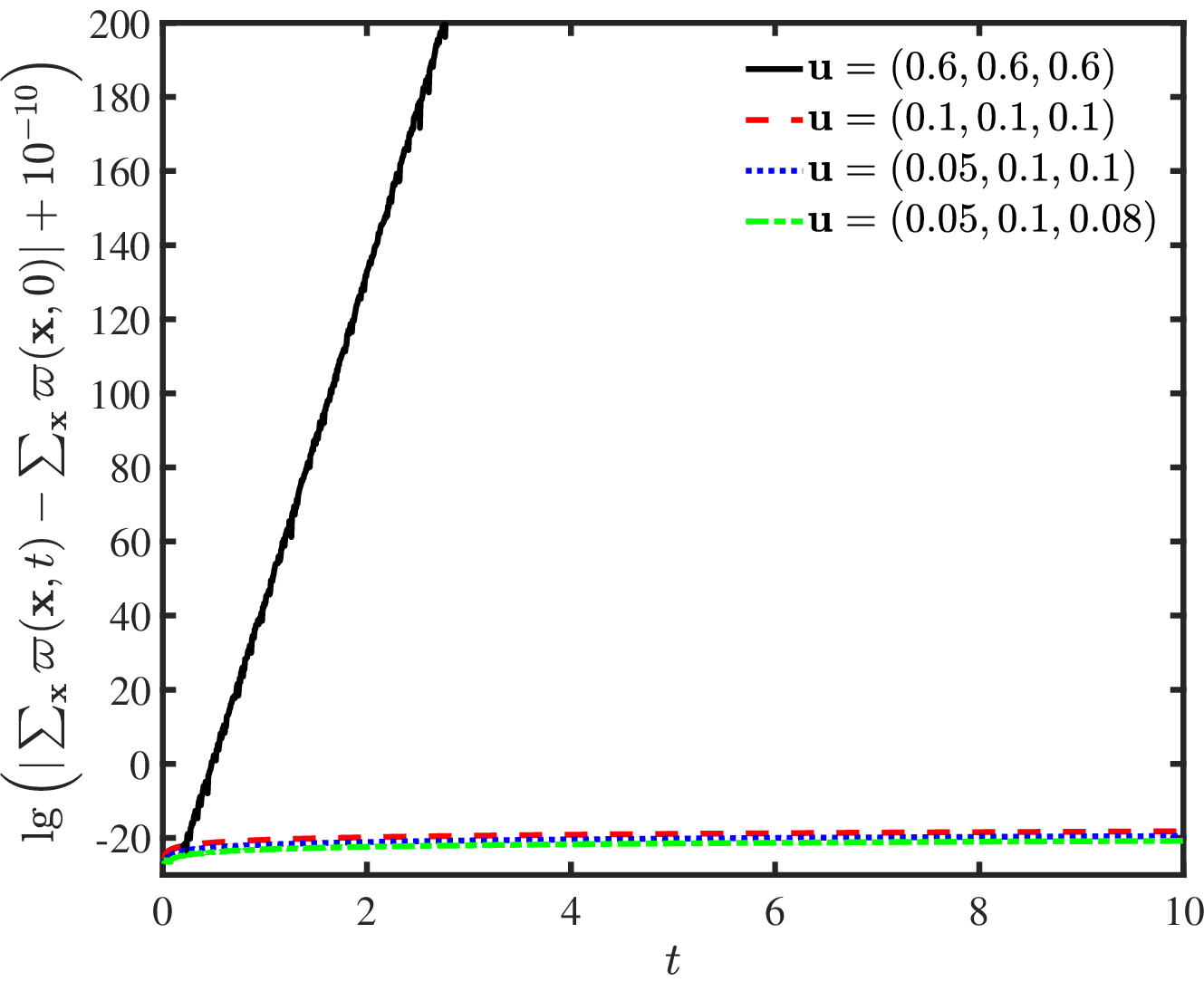} 
			}  
			\caption{The evolution of the total microscopic entropy $\varpi(\mathbf{f})$ with time. } 
			\label{Fig-tanh}  
	\end{figure}

In our simulations, the radius $Rr=1/[5(d-1)]$ and the interfacial parameter $W=(5-d)/75$. We first conduct some simulations under  different values of the transform velocity, and plot the total microscopic entropy $\varpi(\mathbf{f})$ in Fig. \ref{Fig-tanh}, where  one can observe that the BGK-LB model (\ref{BGK-LB-Evolution-seq2})  conserves the  total microscopic entropy unless the microscopic entropy stability conditions are not satisfied, which is in agreement with our theoretical analysis (see Proposition 1 and Remark 6). Then, to test the CR of the BGK-LB model (\ref{BGK-LB-Evolution-seq2}), some simulations are carried out with different accuracy orders of the initialization schemes considered. We measure the RMSEs between the numerical and analytical solutions under different values of the transport velocity at the time $T=1.0$, and calculate the corresponding CRs in Tables \ref{Table-tanh2d} and \ref{Table-tanh3d}.  From these results, one can find that at the acoustic scaling, the BGK-LB model (\ref{BGK-LB-Evolution-seq2}) has a second-order CR in space under the first-order initialization scheme, and has a fourth-order CR in space under the third-order initialization scheme, which are consistent with the conclusions in Refs. \cite{BOGHOSIAN2024,strikwerda2004finite}.\\
 \begin{table} 
	\caption{The RMSEs and CRs of the BGK-LB model  under different values of the transport velocity ($\Delta x=1/100$).} \label{Table-tanh2d}
	\begin{tabular*}{\textwidth}{@{\extracolsep\fill}clccccc}	\hline \hline
		$\mathbf{u}$&Or.-Ini.&RMSE$_{\Delta x,\Delta t}$&RMSE$_{\Delta x/2,\Delta t/4}$&RMSE$_{\Delta x/4,\Delta t/16}$&RMSE$_{\Delta x/8,\Delta t/64}$&CR\\ 
		\midrule 
		\multirow{2}{*}{(0.45,0.45)}&$O(\Delta t)$& 5.3248$\times 10^{-4}$&1.1268 $\times 10^{-4}$&2.8036 $\times 10^{-5}$& 7.0064$\times 10^{-6}$&$\sim$2.0862 \\ 
		\cline{2-7}
		&$O(\Delta t^3)$&3.8584 $\times 10^{-4}$&1.7974 $\times 10^{-5}$&9.5134 $\times 10^{-7}$& 5.7536$\times 10^{-8}$&$\sim$4.2371 \\ 
		\hline 
		\multirow{2}{*}{(0.4,0.4)}&$O(\Delta t)$&6.6315$\times 10^{-4}$&1.1029$\times 10^{-4}$&2.6940$\times 10^{-5}$&6.7289$\times 10^{-6}$&$\sim$2.2076\\ 
		\cline{2-7}
		&$O(\Delta t^3)$&6.0690$\times 10^{-4}$&3.2012$\times 10^{-5}$&1.7255$\times 10^{-6}$&1.0455$\times 10^{-7}$&$\sim$4.1677\\
		\hline
		\multirow{2}{*}{(0.1,0.2)}&$O(\Delta t^3)$&6.9373$\times 10^{-4}$&4.2835$\times 10^{-5}$&2.4055$\times 10^{-6}$&1.4660$\times 10^{-7}$&$\sim$4.0694\\
		\cline{2-7}
		&$O(\Delta t)$&7.4116$\times 10^{-4}$&1.1748$\times 10^{-4}$&2.8829$\times 10^{-5}$&7.1921$\times 10^{-6}$&$\sim$2.2291\\ 
		\hline 
		\multirow{2}{*}{(0.05,0.05)}&$O(\Delta t)$&5.1039$\times 10^{-4}$&1.1907$\times 10^{-4}$&2.9514$\times 10^{-5}$&7.3654$\times 10^{-6}$&$\sim$2.0382\\
		\cline{2-7}
		&$O(\Delta t^3)$&2.8047$\times 10^{-4}$&1.3643$\times 10^{-5}$&7.7182$\times 10^{-7}$&4.7140$\times 10^{-8}$&$\sim$4.1795\\
		\hline 
		\multirow{2}{*}{(0.005,0.005)}&$O(\Delta t)$&5.0156$\times 10^{-4}$&1.2161$\times 10^{-4}$&3.0190$\times 10^{-5}$&7.5346$\times 10^{-6}$&$\sim$2.0189\\
		\cline{2-7}
		&$O(\Delta t^3)$&8.5135$\times 10^{-5}$&2.9091$\times 10^{-6}$&1.6864$\times 10^{-7}$&1.0354$\times 10^{-8}$&$\sim$4.3351\\
			\hline \hline
	\end{tabular*}
\end{table}
	\begin{table}
	\caption{The RMSEs and CRs of the BGK-LB model    under different values of the transport velocity.} \label{Table-tanh3d}
	\begin{tabular*}{\textwidth}{@{\extracolsep\fill}cccccccccc}	\hline \hline
		\multirow{2}{*}{Or.-Ini.}&\multirow{2}{*}{$\Delta x$}&	\multicolumn{2}{c}{$\mathbf{u}=\big(0.1,0.1,0.1\big)$}&&	\multicolumn{2}{c}{$\mathbf{u}=\big(0.05,0.1,0.1\big)$}&&
		\multicolumn{2}{c}{$\mathbf{u}=\big(0.05,0.1,0.08\big)$} \\
		\cline{3-4}\cline{6-7}\cline{9-10} 
		&&RMSE&CR&&RMSE&CR&&RMSE&CR  \\
		\midrule
		\multirow{7}{*}{$O(\Delta t)$}&1/150&  
		6.5798$\times 10^{-4}$&-       
		&&6.3617$\times 10^{-4}$&-   
		&&6.1954$\times 10^{-4}$&-  \\
		&1/180&
		4.0596$\times 10^{-4}$&2.6487 
		&&3.9646$\times 10^{-4}$&2.5837 
		&&3.8952$\times 10^{-4}$&2.5453 \\
		&1/210&
		2.7631$\times 10^{-4}$&2.4959
		&&2.7234$\times 10^{-4}$&2.4362
		&&2.6957$\times 10^{-4}$&2.3878  \\
		&1/240&
		2.0225$\times 10^{-4}$&2.3367
		&&2.0057$\times 10^{-4}$&2.2906 
		&&1.9952$\times 10^{-4}$&2.2536 \\ 
		&1/270&
		1.5570$\times 10^{-4}$&2.2209
		&&1.5500$\times 10^{-4}$& 2.1886
		&&1.5465$\times 10^{-4}$&2.1630  \\ 
		&1/300&
		1.2419$\times 10^{-4}$&2.1463
		&&1.2391$\times 10^{-4}$&2.1247
		&&1.2385$\times 10^{-4}$&2.1077  \\ 
		\hline
		\multirow{7}{*}{$O(\Delta t^3)$}&1/150&  
		4.6801$\times 10^{-4}$&-       
		&&4.2799$\times 10^{-4}$&-   
		&&3.9524$\times 10^{-4}$&-  \\
		&1/180&
		2.3223$\times 10^{-4}$& 3.8428 
		&&2.1010$\times 10^{-4}$&3.9027 
		&&1.9219$\times 10^{-4}$&3.9624  \\
		&1/210&
		1.2253$\times 10^{-4}$& 4.1475
		&&1.1027$\times 10^{-4}$&4.1821
		&&1.0017$\times 10^{-4}$&4.2274  \\
		&1/240&
		6.9020$\times 10^{-5}$& 4.2986
		&&6.1832$\times 10^{-5}$&4.3320 
		&&5.5904$\times 10^{-5}$&4.3674  \\ 
		&1/270&
		4.1333$\times 10^{-5}$& 4.3532
		&&3.6930$\times 10^{-5}$&4.3758 
		&&3.3302$\times 10^{-5}$&4.3980  \\ 
		&1/300&
		2.6140$\times 10^{-5}$& 4.3490
		&&2.3319$\times 10^{-5}$&4.3637
		&&2.1002$\times 10^{-5}$&4.3754  \\ 
			\hline \hline
	\end{tabular*}
\end{table}  

\noindent\textbf{Example 3.} The last example considers the four-dimensional L-HE (\ref{HE}), under the proper initial  and boundary conditions, the analytical solution can be given by
\begin{align}
	\phi^{\star}(x_1,\ldots,x_4,t)=\sin\bigg[\pi\Big(\sum_{i=1}^dx_i-u_it\Big)\bigg],\:\: x\in [-1,1]^4,\: t>0,
\end{align}   
in the following simulations, we consider both the periodic and Dirichlet boundary conditions.

\begin{figure} 
	\centering             
		{
			\includegraphics[width=0.4\textwidth]{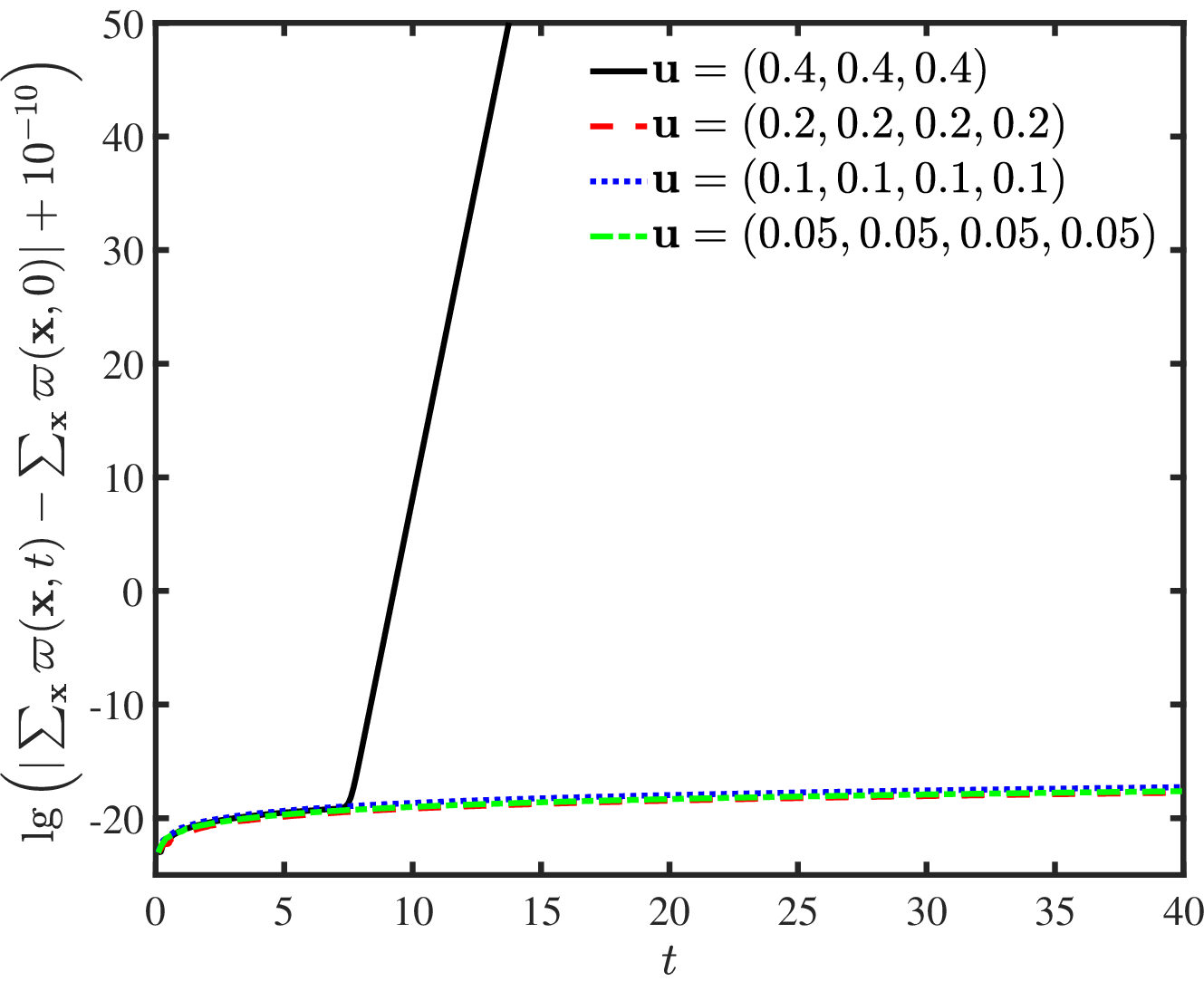} 
		} 
		\caption{The evolution of the total microscopic entropy $\varpi(\mathbf{f})$ with time. } 
		\label{Fig-sin-4d}  
\end{figure}

\begin{figure} 
 	\centering
 	\subfloat[$\mathbf{u}=(0.2,0.1,0.1,0.2)$]    
 	{
 		\includegraphics[width=0.3\textwidth]{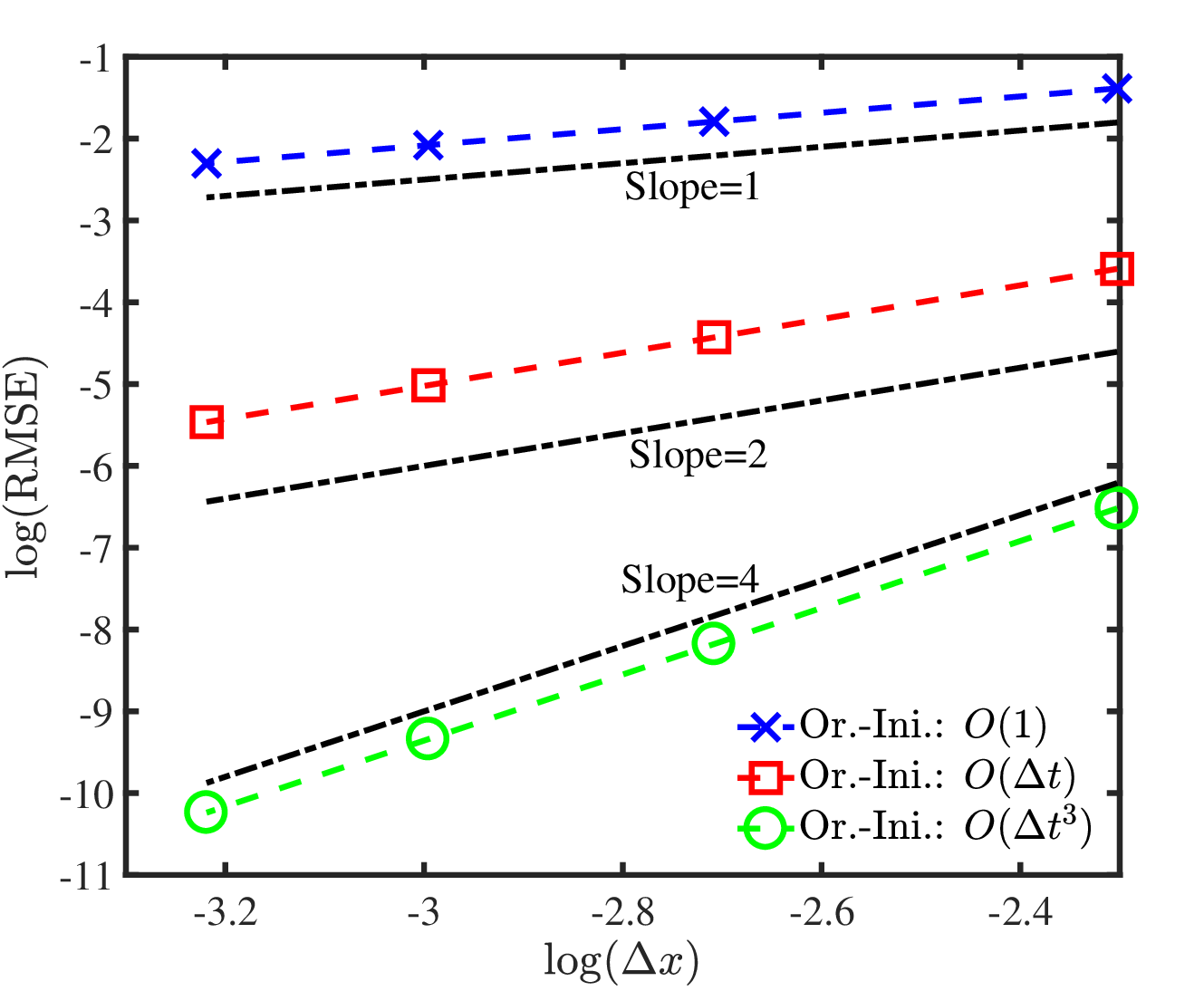} 
 	}       
 	\subfloat[$\mathbf{u}=(0.05,0.2,0.08,0.1)$]    
 	{
 		\includegraphics[width=0.3\textwidth]{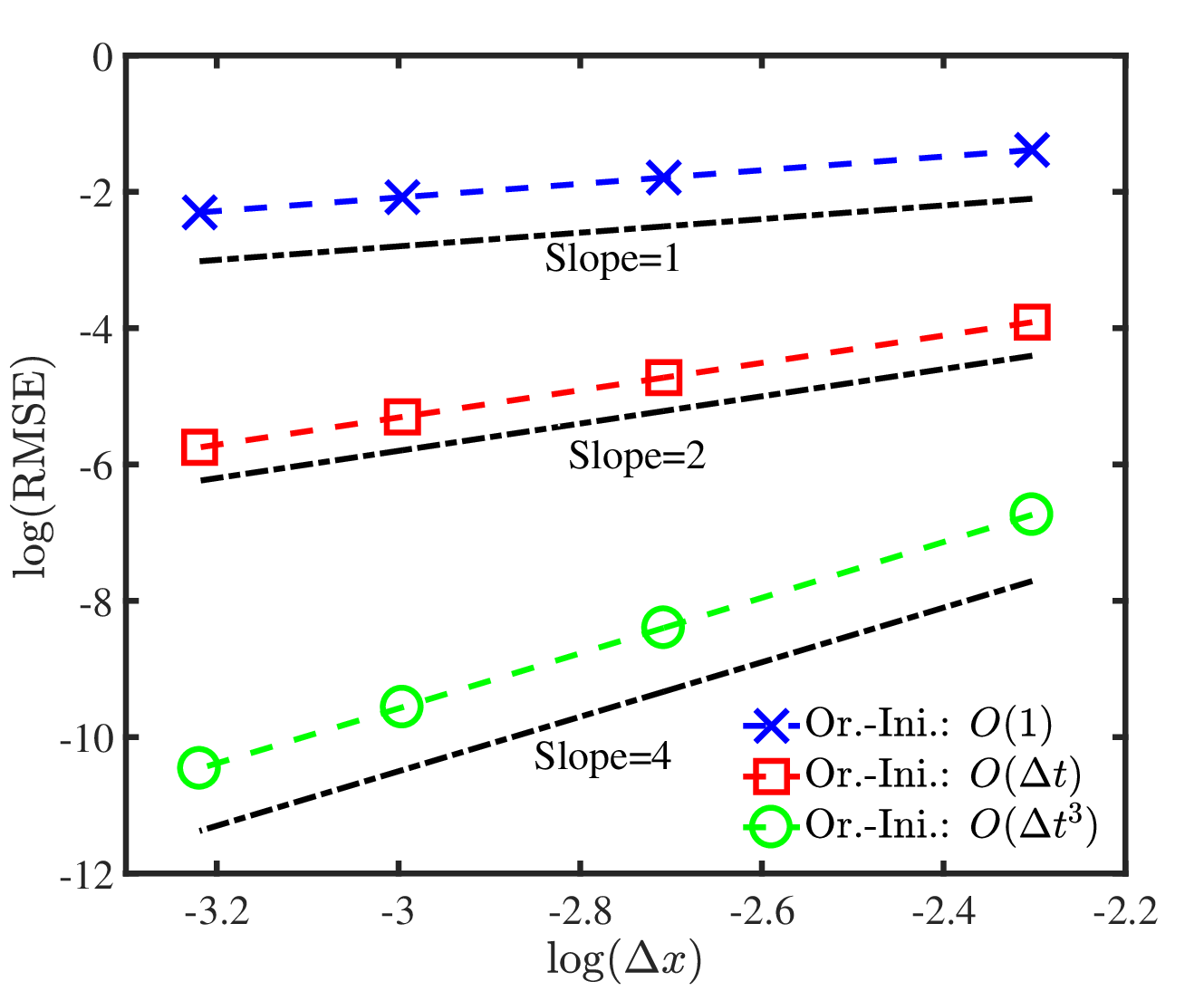} 
 	} 
 	\subfloat[$\mathbf{u}=(10^{-3},10^{-3},10^{-3},10^{-3})$]    
 	{
 		\includegraphics[width=0.3\textwidth]{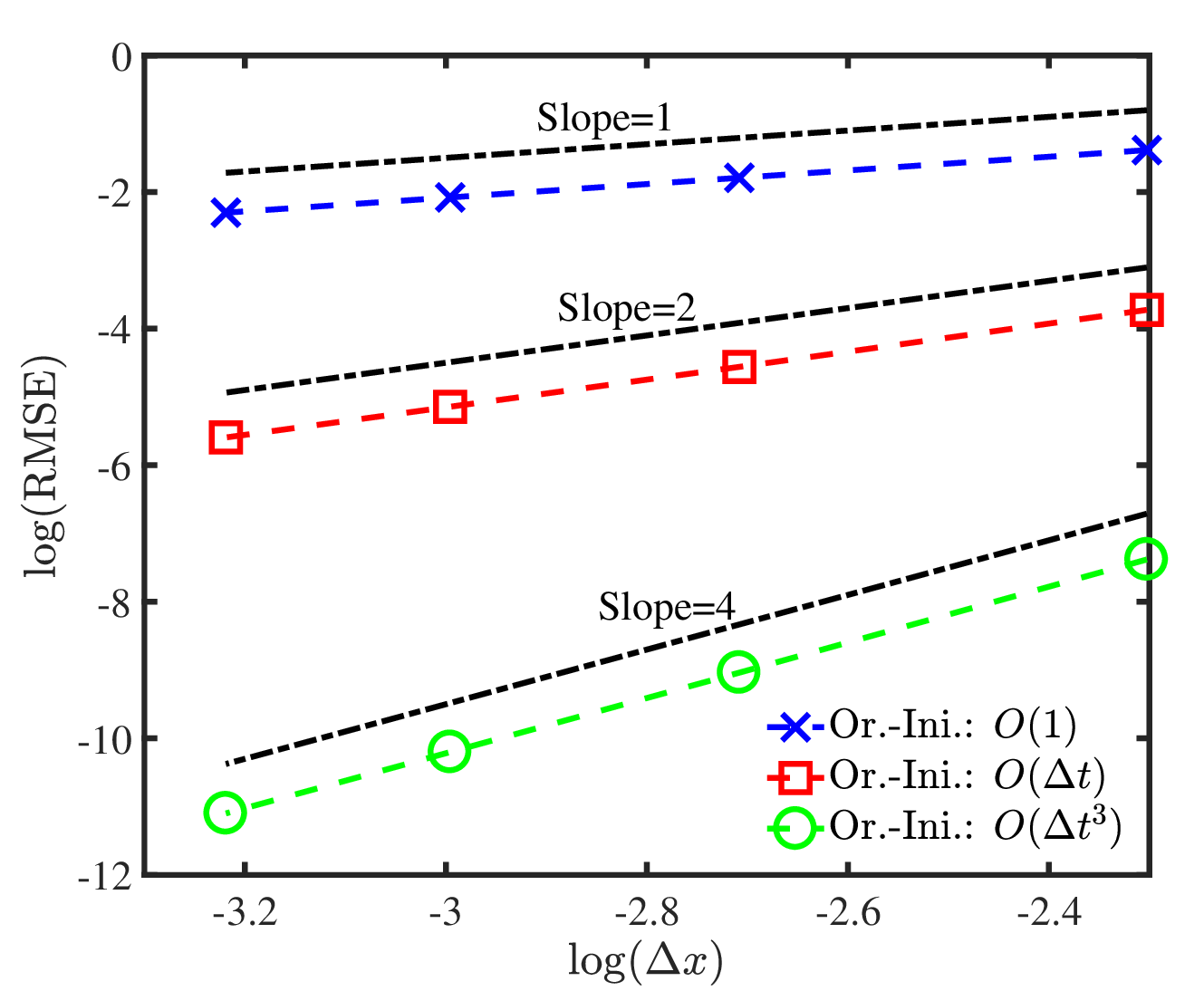} 
 	} 
 	\caption{The convergence rates of the BGK-LB model  under different values of the transport velocity. } 
 	\label{sin-4d-periodic}
 \end{figure}

\begin{figure} 
 	\centering
 	\subfloat[Or.-Boun.: $O(\Delta t)$]    
 	{
 		\includegraphics[width=0.3\textwidth]{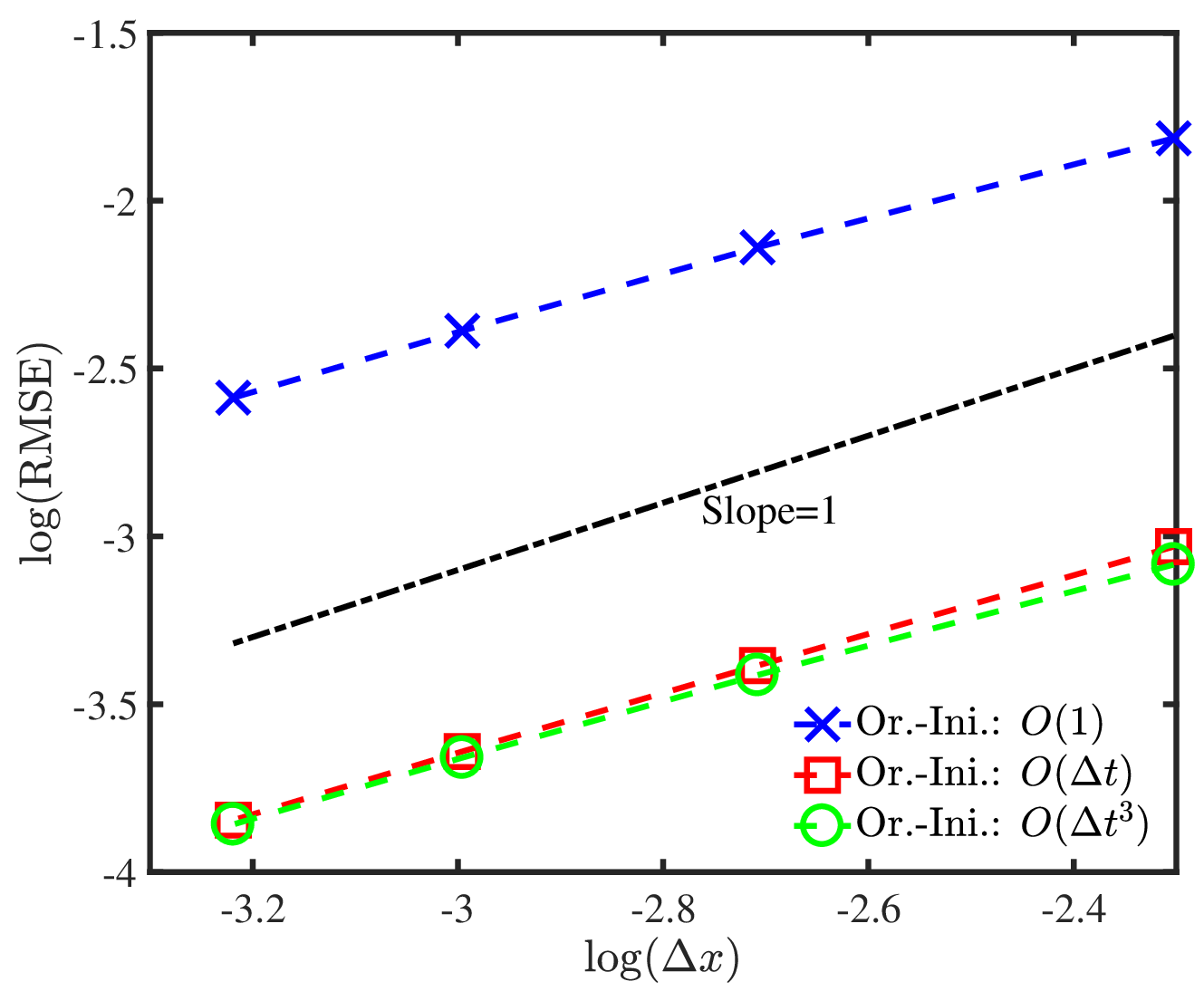} 
 	}       
 	\subfloat[Or.-Boun.: $O(\Delta t^3)$]    
 	{
 		\includegraphics[width=0.3\textwidth]{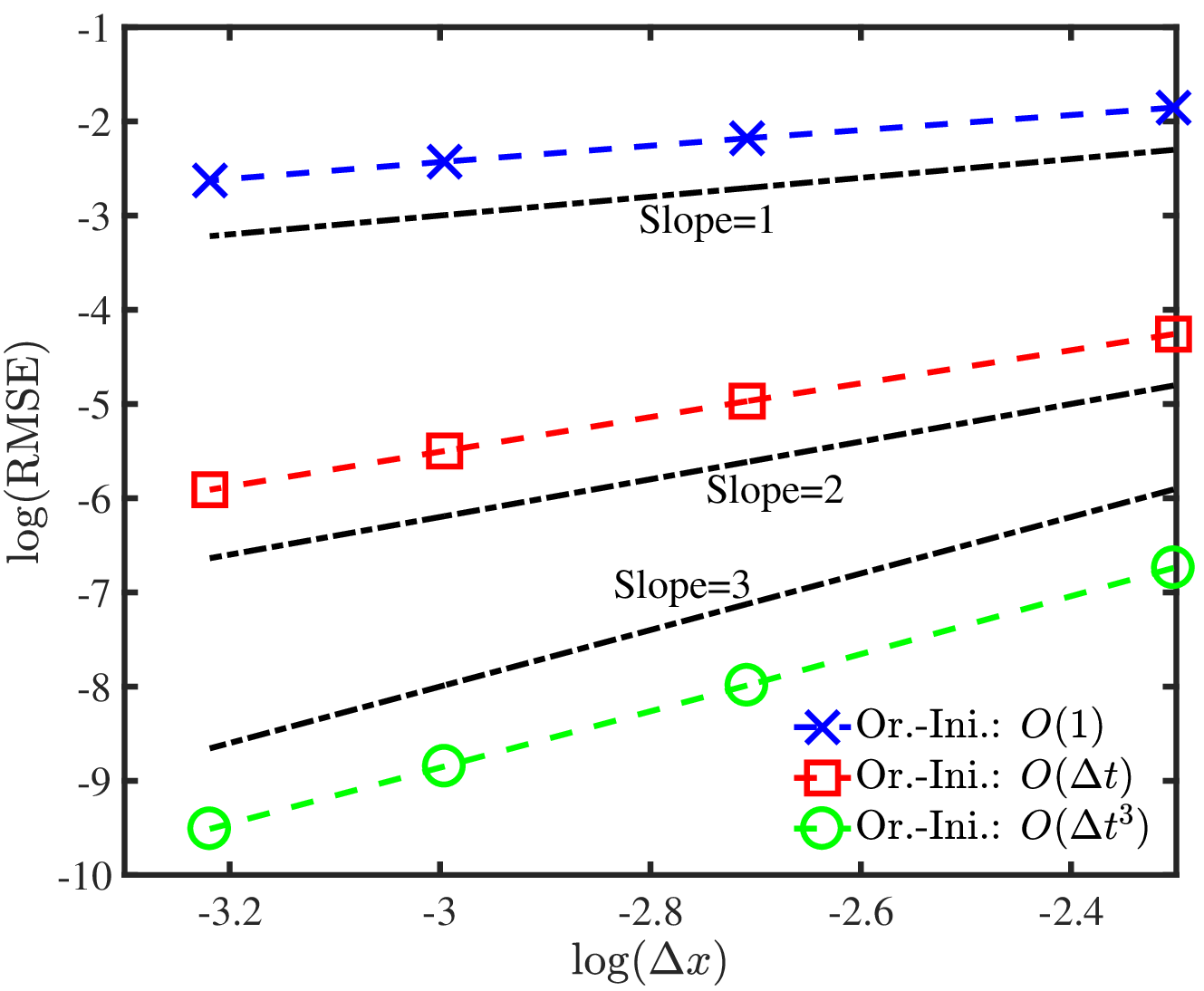} 
 	} 
 	\subfloat[Or.-Boun.: $O(\Delta t^4)$]    
 	{
 		\includegraphics[width=0.3\textwidth]{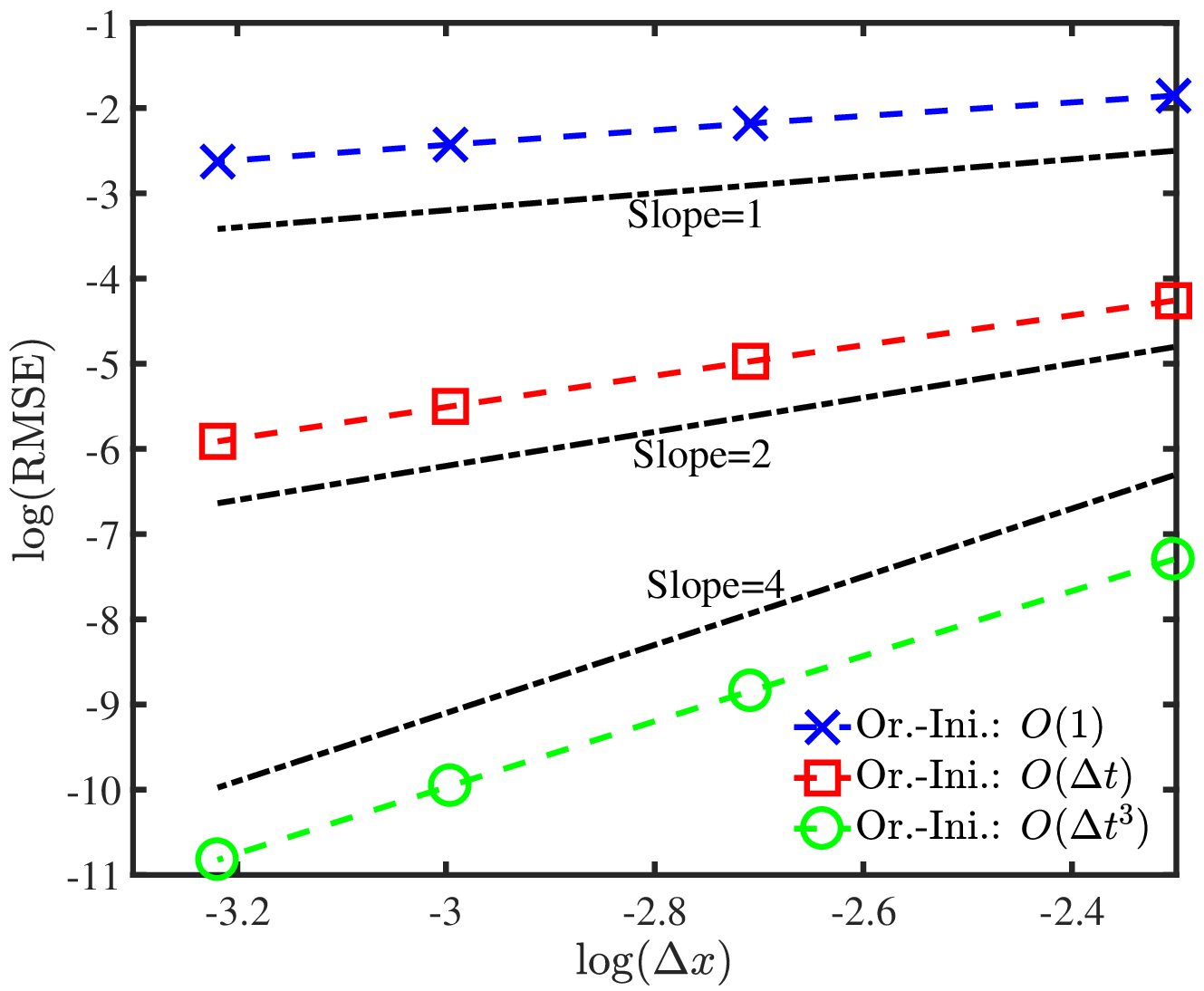} 
 	} 
 	\caption{The convergence rate  of the BGK-LB model  under the transport velocity $\mathbf{u}=(0.05,0.2,0.08,0.1)$. } 
 	\label{sin-4d-ful-diff-u}
 \end{figure}

\begin{figure} 
	\centering
	\subfloat[Or.-Boun.: $O(\Delta t)$]    
	{
		\includegraphics[width=0.3\textwidth]{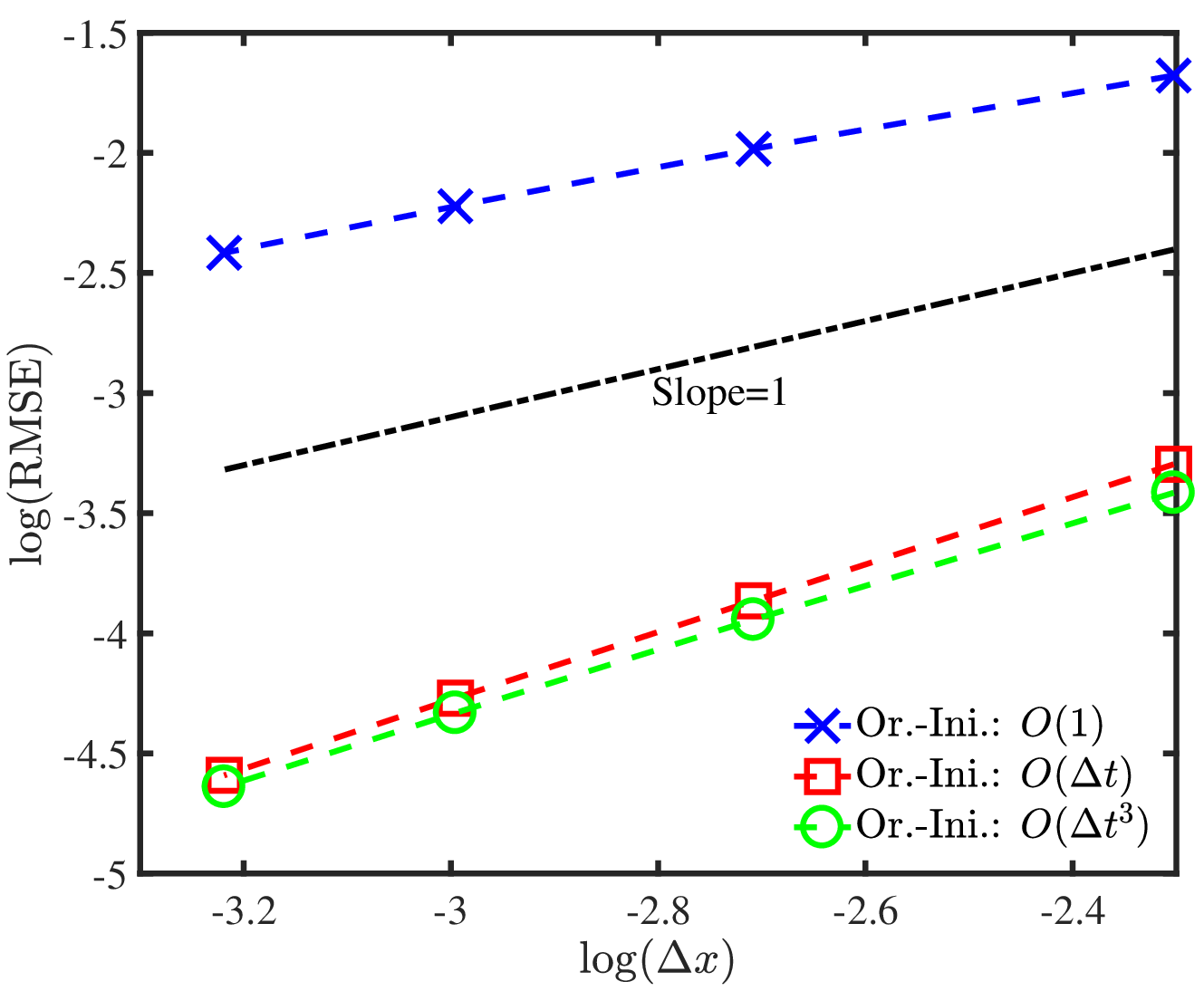} 
	}       
	\subfloat[Or.-Boun.: $O(\Delta t^3)$]    
	{
		\includegraphics[width=0.3\textwidth]{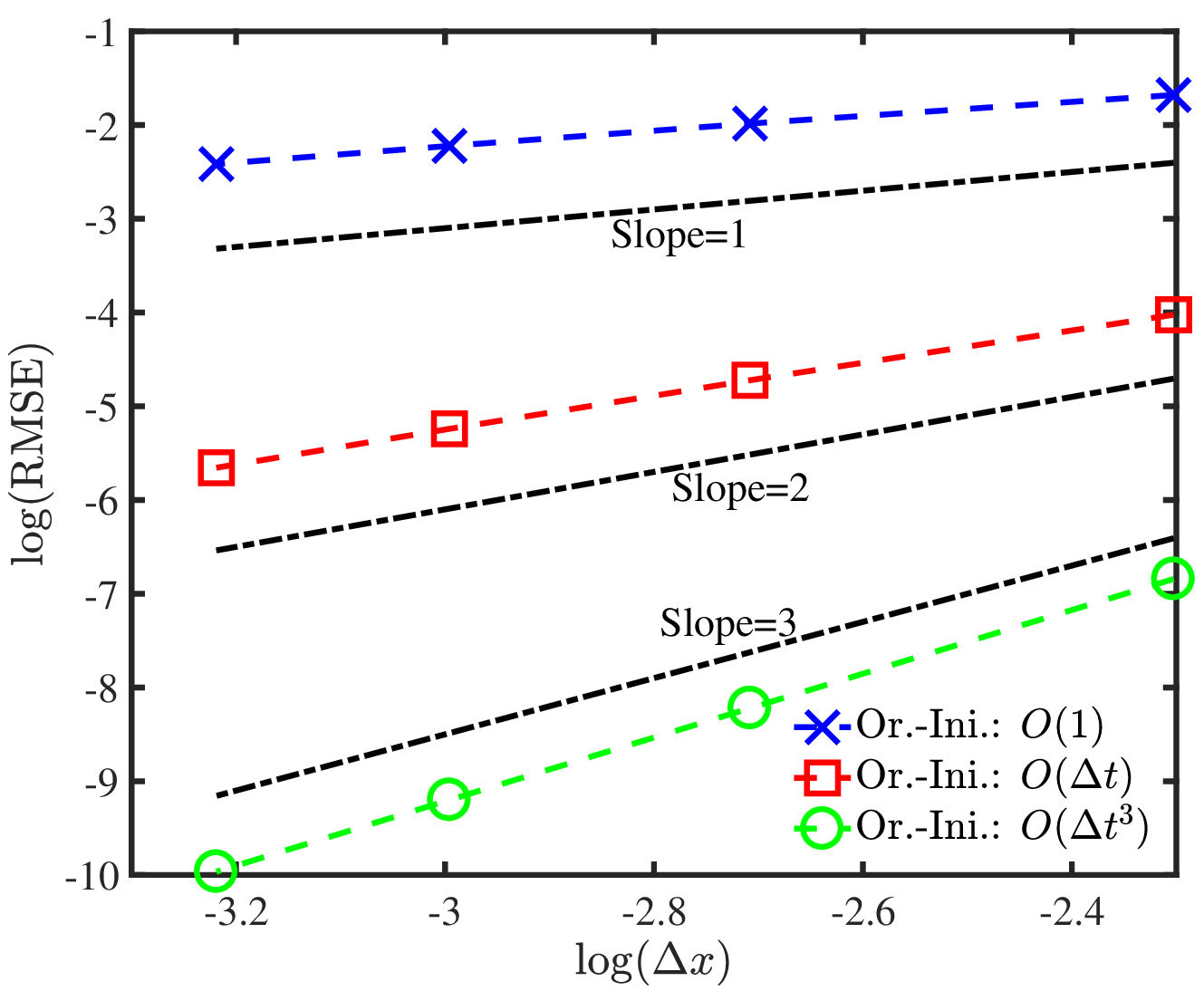} 
	} 
	\subfloat[Or.-Boun.: $O(\Delta t^4)$]    
	{
		\includegraphics[width=0.3\textwidth]{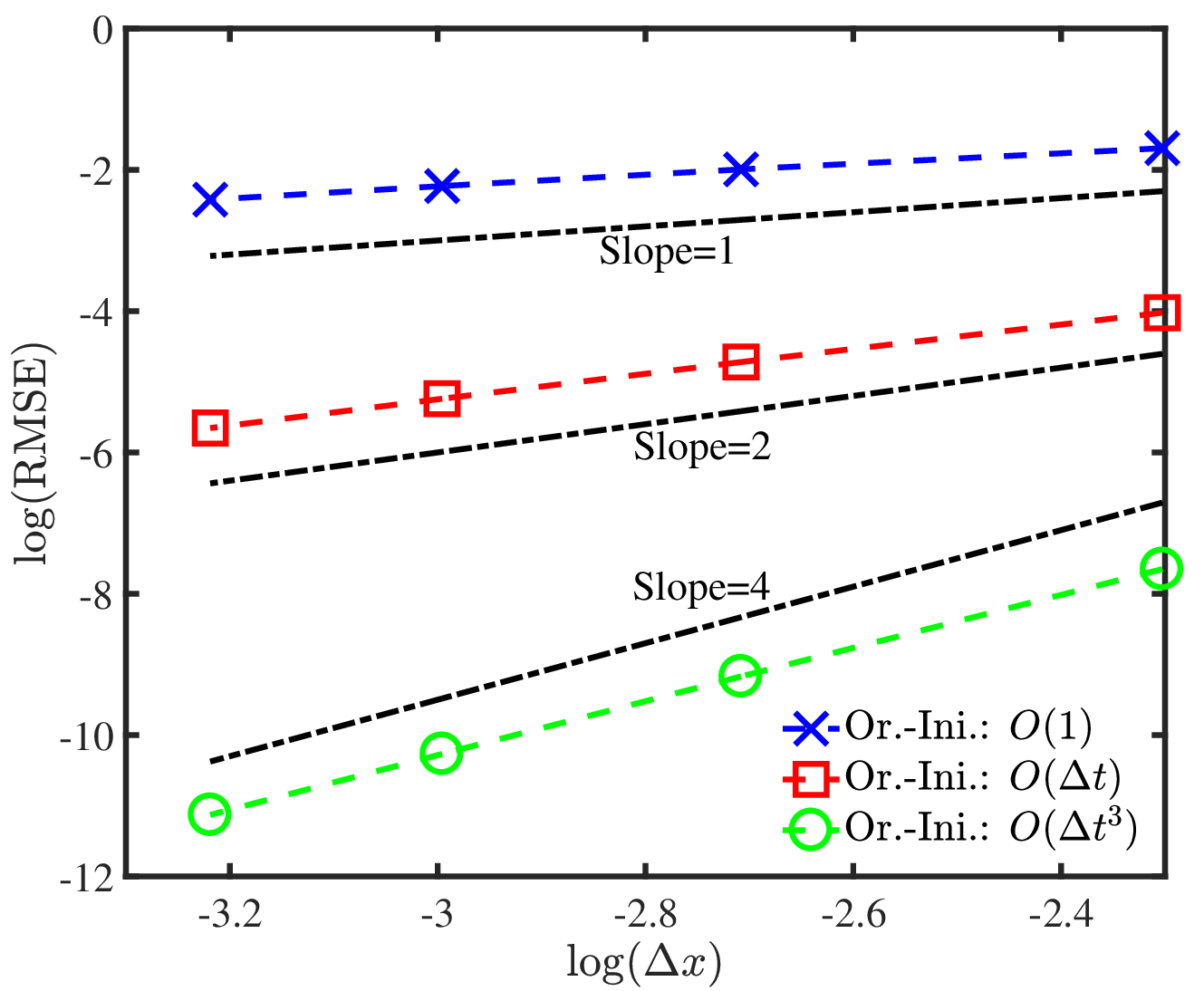} 
	} 
	\caption{The convergence rate  of the BGK-LB model  under  the transport velocity $\mathbf{u}=(10^{-3},10^{-3},10^{-3},10^{-3})$. } 
	\label{sin-4d-ful-same-u}
\end{figure}

\begin{figure} 
 	\centering
 	\subfloat[Or.-Boun.: $O(\Delta t)$]
 	{
 		\includegraphics[width=0.35\textwidth]{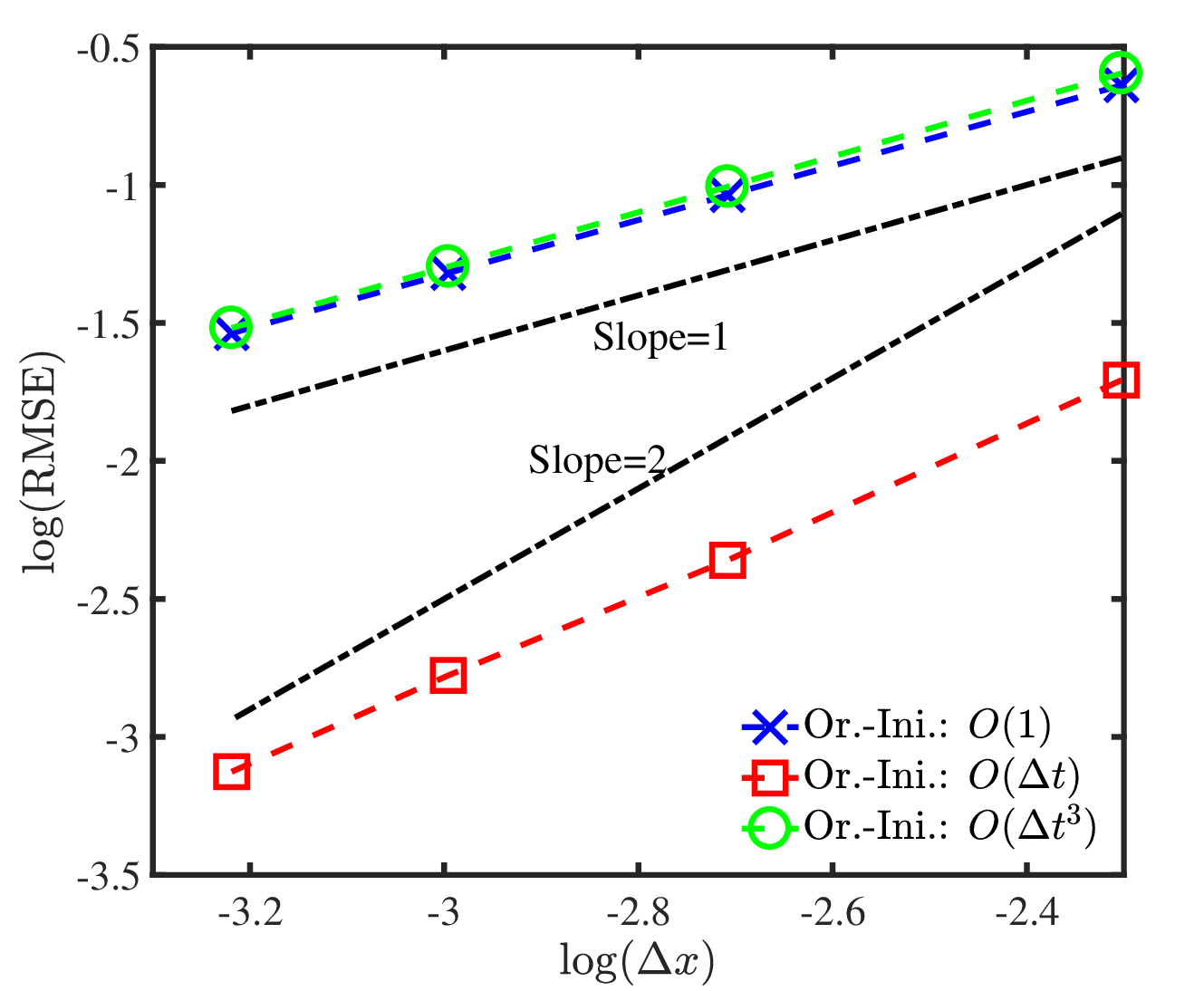} 
 	} 
 	\subfloat[Or.-Boun.: $O(\Delta t^2)$]    
 	{
 		\includegraphics[width=0.35\textwidth]{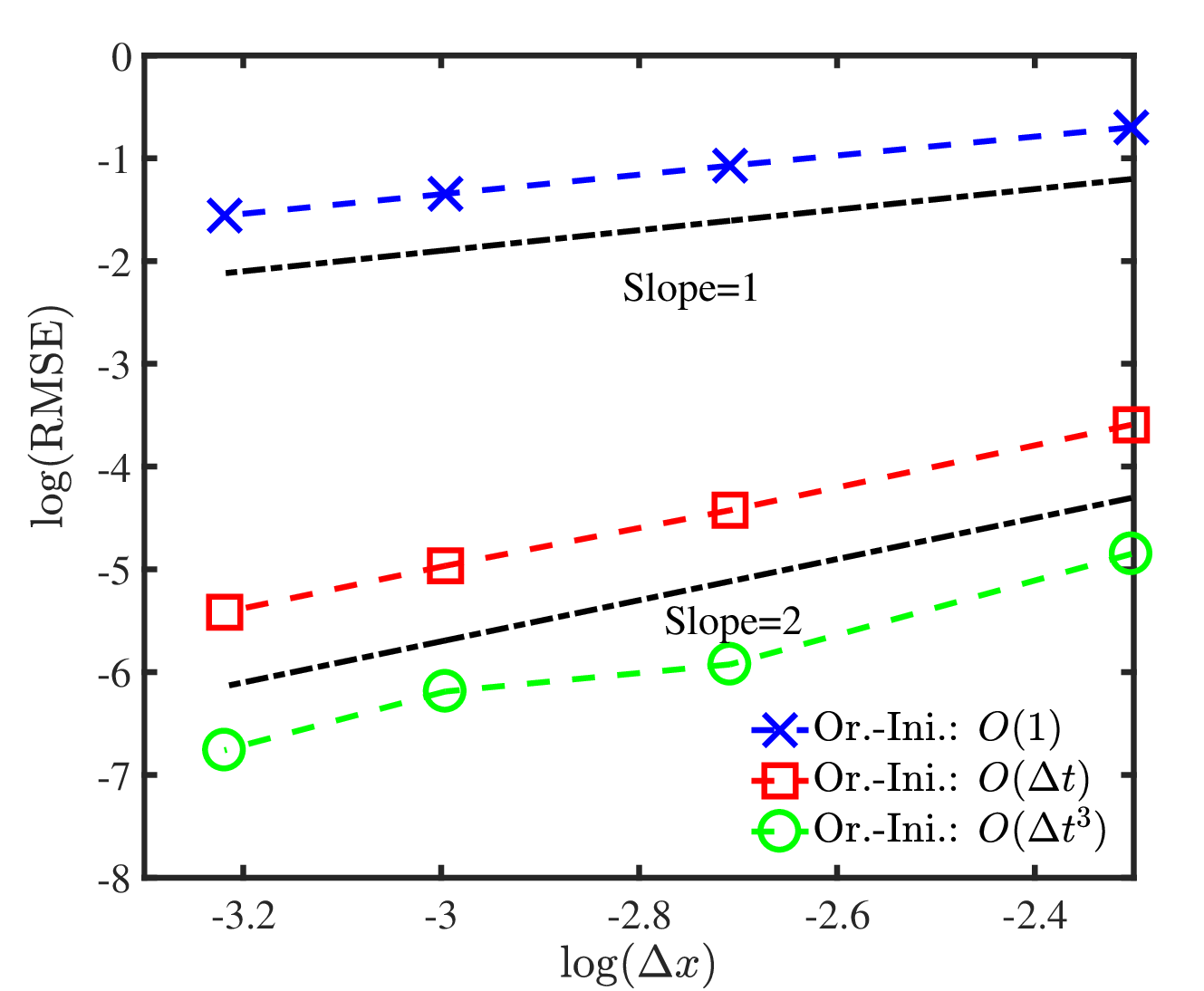} 
 	} 
 	
 	\subfloat[Or.-Boun.: $O(\Delta t^3)$]    
 	{
 		\includegraphics[width=0.35\textwidth]{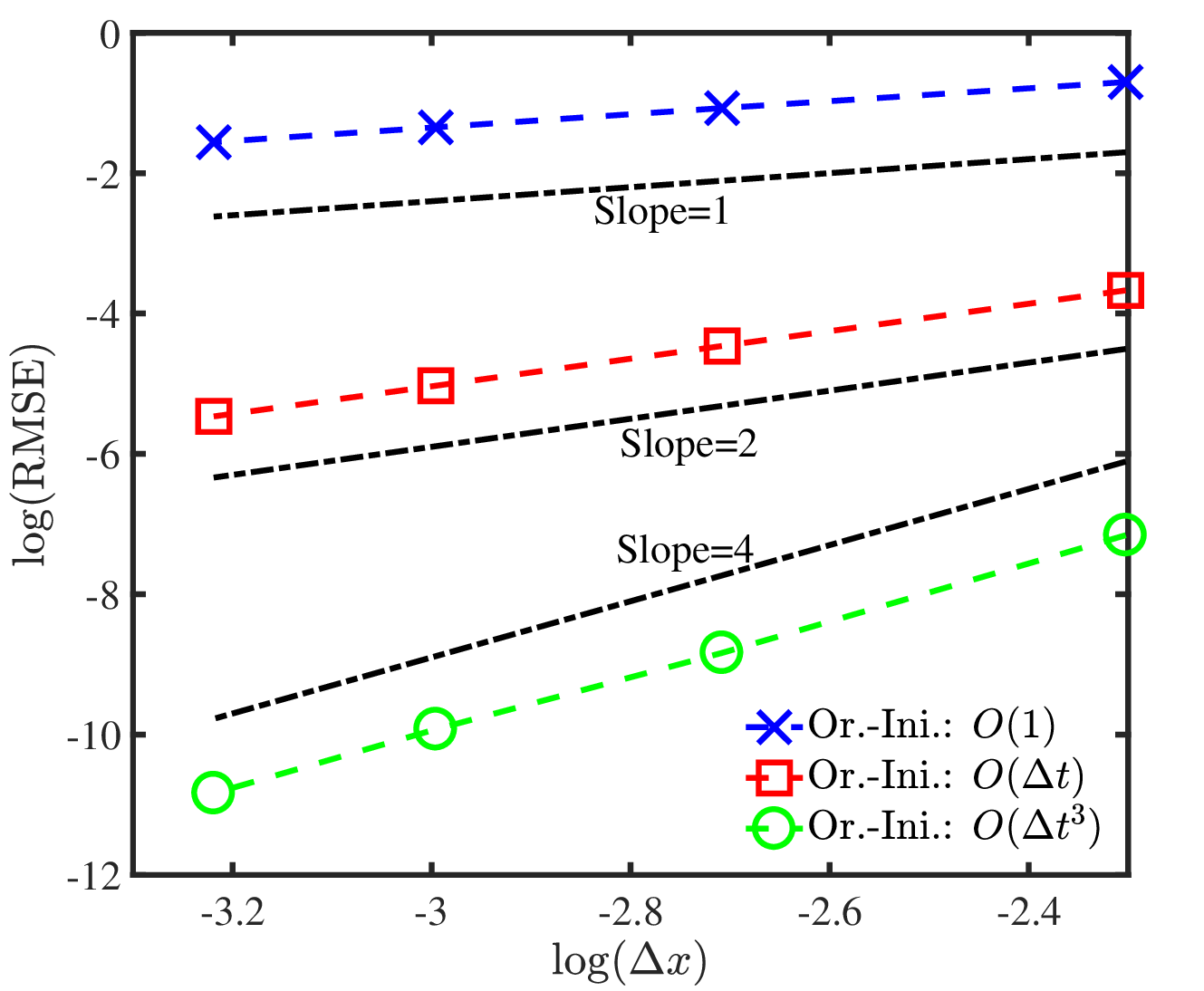} 
 	} 
 	\subfloat[Or.-Boun.: $O(\Delta t^4)$]    
 	{
 		\includegraphics[width=0.35\textwidth]{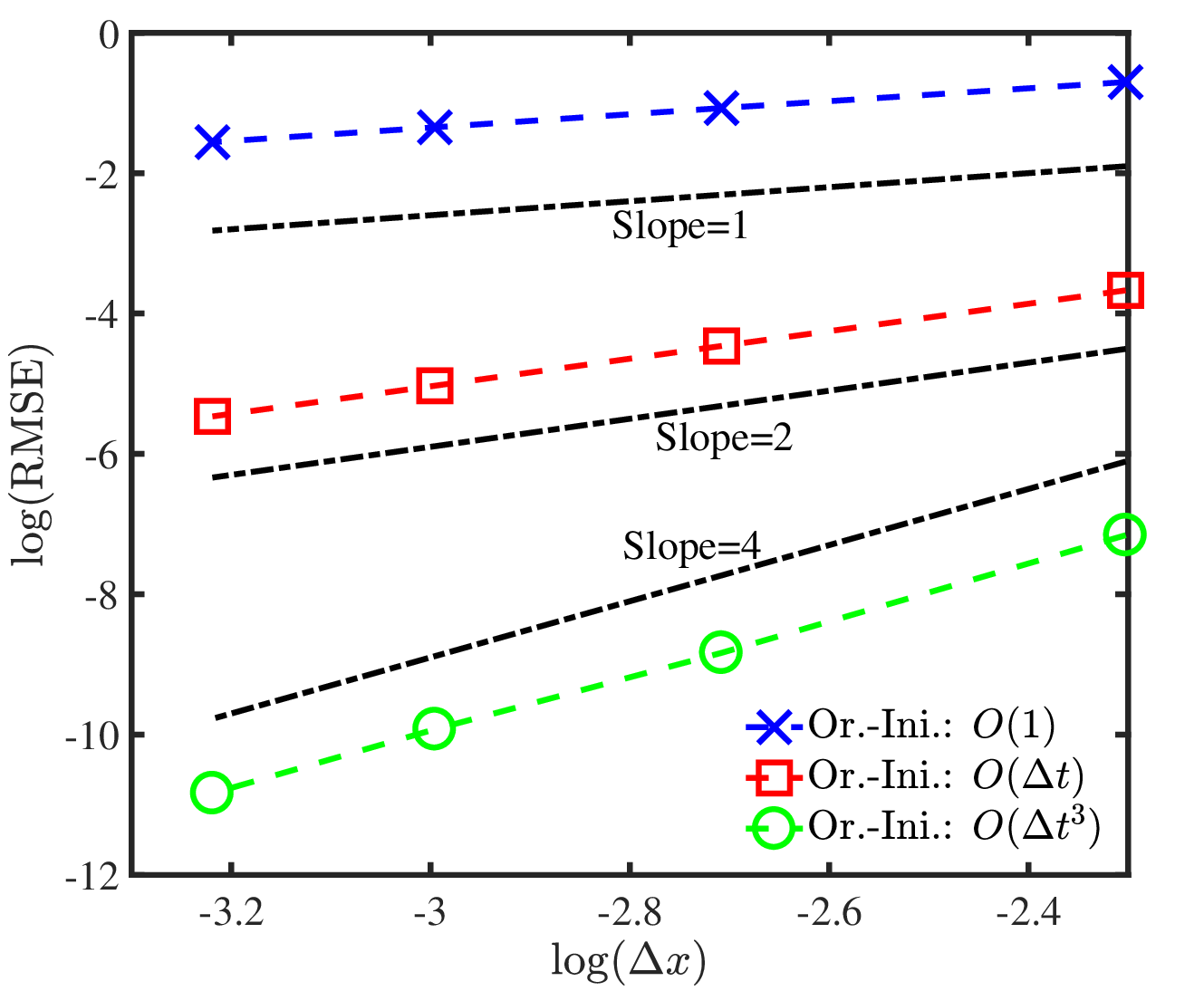} 
 	} 
 	\caption{The convergence rate  of the  BGK-LB model under   the transport velocity $\mathbf{u}=(10^{-3},10^{-3},10^{-3},10^{-3})$. } 
 	\label{sin-4d-half}
 \end{figure}

First of all, when considering the periodic boundary condition, we perform some simulations under several different values of the transport velocity $\mathbf{u}=uI_4^T$ ($u =$ 0.4, 0.2, 0.1, and 0.05), and present the evolution of the total microscopic entropy of the BGK-LB model  (\ref{BGK-LB-Evolution-seq2}) with time in Fig. \ref{Fig-sin-4d}, where the microscopic entropy stability conditions in Eq. (\ref{entropy-stable-condition}) fails to hold when  $u=0.4$. As seen from Fig. \ref{Fig-sin-4d}, it can be found that the case of $u=0.4$  does not conserve the total microscopic entropy while the others do, which is consistent with our discussion in Proposition 1 and Remark 6. Then, we calculate the RMSEs and CRs at the time $T = 2.0$ under several different values of the transport velocity and plot the results in Fig. \ref{sin-4d-periodic}. As shown in Fig. \ref{sin-4d-periodic}, with respect to the periodic boundary condition (here, we would like to point out that the commonly used periodic boundary scheme in the framework of the LB method for the periodic condition is rigorously precise), the relation between the accuracy order of the initialization scheme and the overall CR of the BGK-LB model  (\ref{BGK-LB-Evolution-seq2}) is consistent with Eq. (\ref{cr-estimation}). Finally, we consider the Dirichlet boundary condition with the full-way and half-way boundary schemes adopted, respectively. For the full-way boundary scheme, as seen from Figs.  \ref{sin-4d-ful-diff-u} and \ref{sin-4d-ful-same-u}, where the simulations are suspended at the time $T=2.0$ , it can be seen that the relation among the overall CR of the BGK-LB model (\ref{BGK-LB-Evolution-seq2}), the accuracy order of the initialization scheme, and that of the full-way boundary scheme is consistent with Eq. (\ref{cr-estimation}).  While for the half-way boundary scheme, since the BGK-LB model (\ref{BGK-LB-Evolution-seq2}) does not converge when the transport velocity  $\mathbf{u}=(0.05,0.2,0.08,0.1)$ in the simulation, we only present the RMSEs and CRs with different lattice sizes under the transport velocity $\mathbf{u}=(10^{-3},10^{-3},10^{-3},10^{-3})$ in Fig. \ref{sin-4d-half}. Thus, we can also conclude that the full-way scheme is more stable than the half-way scheme for the four-dimensional case.
 
	\section{Conclusions}\label{Conclusion} 
In this paper, we first develop a unified BGK-LB model for the $d$-dimensional L-HE by adopting the natural moments and the D$d$Q$(2d^2 + 1)$ lattice structure. Then, at the acoustic scaling, with the help of the DTE method, we derive the unified forms of the second- and third-order moments of the EDFs $f_k^{eq}$, which can ensure that the BGK-LB model is fourth-order consistent with the L-HE. Furthermore, to obtain an overall fourth-order BGK-LB model, we discuss the initialization and boundary treatments. For the initialization condition, we provide the first- to third-order schemes of the distribution functions $f_k$. For the Dirichlet boundary condition, we consider the full-way and half-way boundary schemes, respectively, and present the first- to fourth-order schemes of the distribution functions $f_k$. After that, we analyze the microscopic entropy stability of the BGK-LB model based on the kinetic entropy theory to determine the fourth-order moments of the EDFs $f_k^{eq}$. Moreover, we also discuss the $L^2$ stability of the BGK-LB model through the von Neumann stability analysis. Finally, some numerical experiments are conducted to test the accuracy and stability of the BGK-LB model, and the numerical results are consistent with our theoretical analysis. In particular, for the Dirichlet boundary condition, the developed full-way boundary scheme is more stable than the half-way boundary scheme.

\section*{Acknowledgments}
This research was supported by the National Natural Science Foundation of China (Grants No. 12072127 and No. 123B2018), the Interdisciplinary Research Program of HUST (2024JCYJ001 and 2023JCYJ002), and the Fundamental Research Funds for the Central Universities, HUST (No. YCJJ20241101 and No. 2024JYCXJJ016). The computation was completed on the HPC Platform of Huazhong University of Science and Technology.
\appendix
\section{The construction of the half-way boundary scheme for the Dirichlet boundary condition}\label{Appendix-half-way}
In order to obtain the unknown distribution functions $f_k(\mathbf{x}_f,t)$ for the Dirichlet boundary condition in the case where the boundary is located in the middle of two lattice nodes, we here construct the fourth-order half-way scheme of the distribution functions  $f_k(\mathbf{x}_f,t)$ based on the commonly used half-way anti-bounce-back scheme \cite{Zhao2019-2,ZHANG2019}. Firstly, we expand $f_k(\mathbf{x}_f,t+\Delta t)$ at position $\mathbf{x}_b$ and time $t$ as
\begin{align}\label{fi}
	f_k(\mathbf{x}_f,t+\Delta t)=&f_k(\mathbf{x}_b+\mathbf{c}_k\Delta t/2,t+\Delta t)=f_k(\mathbf{x}_b,t)+\frac{\Delta t}{2}\mathbf{c}_k\cdot\nabla f_k(\mathbf{x}_b,t)+\Delta t\partial_tf_k(\mathbf{x}_b,t)\notag\\
	&+\frac{\Delta t^2}{8}\mathbf{c}_k^{.2}\overset{2}{\cdot}\nabla^2 f_k(\mathbf{x}_b,t)+\frac{\Delta t^2}{2}\partial_t^2f_k(\mathbf{x}_b,t)+\frac{\Delta t^2}{2}\mathbf{c}_k\cdot\nabla\partial_tf_k(\mathbf{x}_b,t)\notag\\
	&+\frac{\Delta t^3}{48}\mathbf{c}_k^{.3}\overset{3}{\cdot}\nabla^3 f_k(\mathbf{x}_b,t)+\frac{\Delta t^3}{6}\partial_t^3f_k(\mathbf{x}_b,t)+\frac{\Delta t^3}{4}\mathbf{c}_k\cdot\nabla\partial_t^2f_k(\mathbf{x}_b,t)+\frac{\Delta t^3}{8}\mathbf{c}_k^{.2}\overset{2}{\cdot}\nabla^2\partial_tf_k(\mathbf{x}_b,t)+O(\Delta t^4),
\end{align}
similarly,  $f_{\overline{k}}^{\star}(\mathbf{x}_f,t)$ can be expanded as
\begin{align}\label{fi-star}
	f_{\overline{k}}^{\star}(\mathbf{x}_f,t)=&f_{\overline{k}}(\mathbf{x}_b-\mathbf{c}_k\Delta t/2,t+\Delta t)=f_{\overline{k}}(\mathbf{x}_b,t)-\frac{\Delta t}{2}\mathbf{c}_k\cdot\nabla f_{\overline{k}}(\mathbf{x}_b,t)+\Delta t\partial_tf_{\overline{k}}(\mathbf{x}_b,t)\notag\\
	&+\frac{\Delta t^2}{8}\mathbf{c}_k^{.2}\overset{2}{\cdot}\nabla^2 f_{\overline{k}}(\mathbf{x}_b,t)+\frac{\Delta t^2}{2}\partial_t^2f_{\overline{k}}(\mathbf{x}_b,t)-\frac{\Delta t^2}{2}\mathbf{c}_k\cdot\nabla\partial_tf_{\overline{k}}(\mathbf{x}_b,t)\notag\\
	&-\frac{\Delta t^3}{48}\mathbf{c}_k^{.3}\overset{3}{\cdot}\nabla^3 f_{\overline{k}}(\mathbf{x}_b,t)+\frac{\Delta t^3}{6}\partial_t^3f_{\overline{k}}(\mathbf{x}_b,t)-\frac{\Delta t^3}{4}\mathbf{c}_k\cdot\nabla\partial_t^2f_{\overline{k}}(\mathbf{x}_b,t)+\frac{\Delta t^3}{8}\mathbf{c}_k^{.2}\overset{2}{\cdot}\nabla^2\partial_tf_{\overline{k}}(\mathbf{x}_b,t)+O(\Delta t^4).
\end{align}
To simplify the following derivation, we denote  $\psi=\psi(\mathbf{x}_b,t)$ ($\psi=f_k$ and $f_k^{eq}$). With the help of the first- to fourth-order approximations of the distribution functions $f_k$ at the EDFs $f_k^{eq}$ in Eqs. (\ref{ME-1}), (\ref{ME-2}), (\ref{eq-od3}), and (\ref{eq-od4}),  Eq. (\ref{fi}) can be expressed as the following form which is only related to the EDFs $f_k^{eq}$   with a fourth-order truncation error: 
\begin{align}
	f_k(\mathbf{x}_f,t+\Delta t)=&f_k^{eq}-\frac{\Delta t}{2}\mathbf{c}_k\cdot\nabla f_k^{eq}-\frac{\Delta t}{2}\partial_tf_k^{eq}\notag\\
	&+\frac{\Delta t^3}{8}\mathbf{c}_k\cdot\nabla \partial_t^2f_k^{eq}+\frac{\Delta t^3}{8}\mathbf{c}_k^{.2}\overset{2}{\cdot}\nabla^2\partial_tf_k^{eq}+\frac{\Delta t^3}{24}\partial_t^3f_k^{eq}+\frac{\Delta t^3}{24}\mathbf{c}_k^{.3}\overset{3}{\cdot} \nabla^3 f_k^{eq}+\frac{\Delta t}{2}R_k\notag\\
	&+\frac{\Delta t}{2}\mathbf{c}_k\cdot\nabla \Big(f_k^{eq}-\frac{\Delta t}{2}\mathbf{c}_k\cdot\nabla f_k^{eq}-\frac{\Delta t}{2}\partial_tf_k^{eq}\Big)+\Delta t\partial_t\Big(f_k^{eq}-\frac{\Delta t}{2}\mathbf{c}_k\cdot\nabla f_k^{eq}-\frac{\Delta t}{2}\partial_tf_k^{eq}\Big)\notag\\
	&+\frac{\Delta t^2}{8}\mathbf{c}_k^{.2}\overset{2}{\cdot}\nabla^2 \Big(f_k^{eq}-\frac{\Delta t}{2}\mathbf{c}_k\cdot\nabla f_k^{eq}-\frac{\Delta t}{2}\partial_tf_k^{eq}\Big)+\frac{\Delta t^2}{2}\partial_t^2\Big(f_k^{eq}-\frac{\Delta t}{2}\mathbf{c}_k\cdot\nabla f_k^{eq}-\frac{\Delta t}{2}\partial_tf_k^{eq}\Big)\notag\\
	&+\frac{\Delta t^2}{2}\mathbf{c}_k\cdot\nabla\partial_t\Big(f_k^{eq}-\frac{\Delta t}{2}\mathbf{c}_k\cdot\nabla f_k^{eq}-\frac{\Delta t}{2}\partial_tf_k^{eq}\Big)\notag\\
	&+\frac{\Delta t^3}{48}\mathbf{c}_k^{.3}\overset{3}{\cdot}\nabla^3 f_k^{eq}+\frac{\Delta t^3}{6}\partial_t^3f_k^{eq}+\frac{\Delta t^3}{4}\mathbf{c}_k\cdot\nabla\partial_t^2f_k^{eq}+\frac{\Delta t^3}{8}\mathbf{c}_k^{.2}\overset{2}{\cdot}\nabla^2\partial_tf_k^{eq}+O(\Delta t^4),
\end{align}
i.e., 
\begin{align}\label{fi-expand}
	f_k(\mathbf{x}_f,t+\Delta t) =&f_k^{eq}-\frac{\Delta t}{2}\mathbf{c}_k\cdot\nabla f_k^{eq}-\frac{\Delta t}{2}\partial_tf_k^{eq} +\frac{\Delta t}{2}\mathbf{c}_k\cdot\nabla f_k^{eq}+\Delta t\partial_tf_k^{eq}\notag\\
	&-\frac{\Delta t^2}{2}\mathbf{c}_k\cdot\nabla \partial_tf_k^{eq}-\frac{\Delta t^2}{4}\mathbf{c}_k\cdot\nabla\partial_tf_k^{eq}+\frac{\Delta t^2}{2}\mathbf{c}_k\cdot\nabla\partial_tf_k^{eq}\notag\\
	&-\frac{\Delta t^2}{4}\mathbf{c}_k^{.2}\overset{2}{\cdot}\nabla^2 f_k^{eq}+\frac{\Delta t^2}{8}\mathbf{c}_k^{.2}\overset{2}{\cdot}\nabla^2 f_k^{eq}+\frac{\Delta t^2}{2}\partial_t^2f_k^{eq}-\frac{\Delta t^2}{2}\partial_t^2f_k^{eq}\notag\\
	&-\frac{\Delta t^3}{16} \mathbf{c}_k^{.3}\overset{2}{\cdot}\nabla^3  f_k^{eq}+\frac{\Delta t^3}{48}\mathbf{c}_k^{.3}\overset{3}{\cdot}\nabla^3 f_k^{eq}+\frac{\Delta t^3}{24}\mathbf{c}_k^{.3}\overset{3}{\cdot} \nabla^3 f_k^{eq}\notag\\
	&-\frac{\Delta t^3}{16} \mathbf{c}_k^{.2}\overset{2}{\cdot}\nabla^2 \partial_tf_k^{eq} -\frac{\Delta t^3}{4} \mathbf{c}_k^{.2}\overset{2}{\cdot}\nabla^2\partial_t f_k^{eq}+\frac{\Delta t^3}{8}\mathbf{c}_k^{.2}\overset{2}{\cdot}\nabla^2\partial_tf_k^{eq}+\frac{\Delta t^3}{8}\mathbf{c}_k^{.2}\overset{2}{\cdot}\nabla^2\partial_tf_k^{eq}\notag\\
	&-\frac{\Delta t^3}{4}\mathbf{c}_k\cdot\nabla \partial_t^2f_k^{eq}+\frac{\Delta t^3}{4}\mathbf{c}_k\cdot\nabla\partial_t^2f_k^{eq}+\frac{\Delta t^3}{8}\mathbf{c}_k\cdot\nabla \partial_t^2f_k^{eq} -\frac{\Delta t^3}{4} \mathbf{c}_k\cdot\nabla\partial_t^2f_k^{eq}\notag\\
	&+\frac{\Delta t^3}{6}\partial_t^3f_k^{eq} +\frac{\Delta t^3}{24}\partial_t^3f_k^{eq}-\frac{\Delta t^3}{4}\partial_t^3f_k^{eq}+\frac{\Delta t}{2}R_k+O(\Delta t^4),
\end{align}
simplifying  Eq. (\ref{fi-expand}) yields
\begin{align}\label{fi-2}
	f_k(\mathbf{x}_f,t+\Delta t)= &f_k^{eq} +\frac{\Delta t}{2}\partial_tf_k^{eq}-\frac{\Delta t^2}{8}\mathbf{c}_k^{.2}\overset{2}{\cdot}\nabla^2 f_k^{eq}-\frac{\Delta t^2}{4}\mathbf{c}_k\cdot\nabla\partial_tf_k^{eq} \notag\\
	&-\frac{\Delta t^3}{8}\mathbf{c}_k\cdot\nabla \partial_t^2f_k^{eq}-\frac{\Delta t^3}{16}\mathbf{c}_k^{.2}\overset{2}{\cdot}\nabla^2\partial_tf_k^{eq}-\frac{\Delta t^3}{24}\partial_t^3f_k^{eq} +\frac{\Delta t}{2}R_k +O(\Delta t^4).
\end{align}
Similarly, substituting Eqs. (\ref{ME-1}), (\ref{ME-2}), (\ref{eq-od3}), and (\ref{eq-od4}) into Eq. (\ref{fi-star}) without degrading the fourth-order accuracy, $f_{\overline{k}}^{\star}$ can be also given by
\begin{align}\label{App-eq1}
	f_{\overline{k}}^{\star}(\mathbf{x}_f,t)=&f_{\overline{k}}^{eq}-\frac{\Delta t}{2}\mathbf{c}_{\overline{k}}\cdot\nabla f_{\overline{k}}^{eq}-\frac{\Delta t}{2}\partial_tf_{\overline{k}}^{eq}\notag\\
	&+\frac{\Delta t^3}{8}\mathbf{c}_{\overline{k}}\cdot\nabla \partial_t^2f_{\overline{k}}^{eq}+\frac{\Delta t^3}{8}\mathbf{c}_{\overline{k}}^{.2}\overset{2}{\cdot}\nabla^2\partial_tf_{\overline{k}}^{eq}+\frac{\Delta t^3}{24}\partial_t^3f_{\overline{k}}^{eq}+\frac{\Delta t^3}{24}\mathbf{c}_{\overline{k}}^{.3}\overset{3}{\cdot} \nabla^3 f_{\overline{k}}^{eq}+\frac{\Delta t}{2}R_k\notag\\
	&-\frac{\Delta t}{2}\mathbf{c}_k\cdot\nabla f_{\overline{k}}^{eq}+\frac{\Delta t^2}{4}\mathbf{c}_k\mathbf{c}_{\overline{k}}\overset{2}{\cdot}\nabla^2 f_{\overline{k}}^{eq}+\frac{\Delta t^2}{4}\mathbf{c}_k\cdot\nabla\partial_tf_{\overline{k}}^{eq}\notag\\
	&+\Delta t\partial_tf_{\overline{k}}^{eq}-\frac{\Delta t^2}{2}\mathbf{c}_{\overline{k}}\cdot\nabla \partial_tf_{\overline{k}}^{eq}-\frac{\Delta t^2}{2}\partial_t^2f_{\overline{k}}^{eq}\notag\\
	&+\frac{\Delta t^2}{8}\mathbf{c}_k^{.2}\overset{2}{\cdot}\nabla^2 f_{\overline{k}}^{eq}-\frac{\Delta t^3}{16} \mathbf{c}_k^{.2}\mathbf{c}_{\overline{k}}\overset{2}{\cdot}\nabla^3  f_{\overline{k}}^{eq}-\frac{\Delta t^3}{16} \mathbf{c}_k^{.2}\overset{2}{\cdot}\nabla^2 \partial_tf_{\overline{k}}^{eq}\notag\\
	&+\frac{\Delta t^2}{2}\partial_t^2f_{\overline{k}}^{eq}-\frac{\Delta t^3}{4}\mathbf{c}_{\overline{k}}\cdot\nabla \partial_t^2f_{\overline{k}}^{eq}-\frac{\Delta t^3}{4}\partial_t^3f_{\overline{k}}^{eq}\notag\\
	&-\frac{\Delta t^2}{2}\mathbf{c}_k\cdot\nabla\partial_tf_{\overline{k}}^{eq}+\frac{\Delta t^3}{4} \mathbf{c}_k\mathbf{c}_{\overline{k}}\overset{2}{\cdot}\nabla^2\partial_t f_{\overline{k}}^{eq}+\frac{\Delta t^3}{4} \mathbf{c}_k\cdot\nabla\partial_t^2f_{\overline{k}}^{eq}\notag\\
	&-\frac{\Delta t^3}{48}\mathbf{c}_k^{.3}\overset{3}{\cdot}\nabla^3 f_{\overline{k}}^{eq}+\frac{\Delta t^3}{6}\partial_t^3f_{\overline{k}}^{eq}-\frac{\Delta t^3}{4}\mathbf{c}_k\cdot\nabla\partial_t^2f_{\overline{k}}^{eq}+\frac{\Delta t^3}{8}\mathbf{c}_k^{.2}\overset{2}{\cdot}\nabla^2\partial_tf_{\overline{k}}^{eq}+O(\Delta t^4),
\end{align}
with the relation of $\mathbf{c}_{\overline{k}}=-\mathbf{c}_k$, the above equation becomes
\begin{align}\label{App-eq2}
	f_{\overline{k}}^{\star}(\mathbf{x}_f,t)=&f_{\overline{k}}^{eq}+\frac{\Delta t}{2}\mathbf{c}_k\cdot\nabla f_{\overline{k}}^{eq}-\frac{\Delta t}{2}\partial_tf_{\overline{k}}^{eq}\notag\\
	&-\frac{\Delta t^3}{8}\mathbf{c}_k\cdot\nabla \partial_t^2f_{\overline{k}}^{eq}+\frac{\Delta t^3}{8}\mathbf{c}_k^{.2}\overset{2}{\cdot}\nabla^2\partial_tf_{\overline{k}}^{eq}+\frac{\Delta t^3}{24}\partial_t^3f_{\overline{k}}^{eq}-\frac{\Delta t^3}{24}\mathbf{c}_k^{.3}\overset{3}{\cdot} \nabla^3 f_{\overline{k}}^{eq}+\frac{\Delta t}{2}R_k\notag\\
	&-\frac{\Delta t}{2}\mathbf{c}_k\cdot\nabla f_{\overline{k}}^{eq}-\frac{\Delta t^2}{4}\mathbf{c}_k\mathbf{c}_k\overset{2}{\cdot}\nabla^2 f_{\overline{k}}^{eq}+\frac{\Delta t^2}{4}\mathbf{c}_k\cdot\nabla\partial_tf_{\overline{k}}^{eq}\notag\\
	&+\Delta t\partial_tf_{\overline{k}}^{eq}+\frac{\Delta t^2}{2}\mathbf{c}_k\cdot\nabla \partial_tf_{\overline{k}}^{eq}-\frac{\Delta t^2}{2}\partial_t^2f_{\overline{k}}^{eq}\notag\\
	&+\frac{\Delta t^2}{8}\mathbf{c}_k^{.2}\overset{2}{\cdot}\nabla^2 f_{\overline{k}}^{eq}+\frac{\Delta t^3}{16} \mathbf{c}_k^{.3}\overset{2}{\cdot}\nabla^3  f_{\overline{k}}^{eq}-\frac{\Delta t^3}{16} \mathbf{c}_k^{.2}\overset{2}{\cdot}\nabla^2 \partial_tf_{\overline{k}}^{eq}\notag\\
	&+\frac{\Delta t^2}{2}\partial_t^2f_{\overline{k}}^{eq}+\frac{\Delta t^3}{4}\mathbf{c}_k\cdot\nabla \partial_t^2f_{\overline{k}}^{eq}-\frac{\Delta t^3}{4}\partial_t^3f_{\overline{k}}^{eq}\notag\\
	&-\frac{\Delta t^2}{2}\mathbf{c}_k\cdot\nabla\partial_tf_{\overline{k}}^{eq}-\frac{\Delta t^3}{4} \mathbf{c}_k^{.2}\overset{2}{\cdot}\nabla^2\partial_t f_{\overline{k}}^{eq}+\frac{\Delta t^3}{4} \mathbf{c}_k\cdot\nabla\partial_t^2f_{\overline{k}}^{eq}\notag\\
	&-\frac{\Delta t^3}{48}\mathbf{c}_k^{.3}\overset{3}{\cdot}\nabla^3 f_{\overline{k}}^{eq}+\frac{\Delta t^3}{6}\partial_t^3f_{\overline{k}}^{eq}-\frac{\Delta t^3}{4}\mathbf{c}_k\cdot\nabla\partial_t^2f_{\overline{k}}^{eq}+\frac{\Delta t^3}{8}\mathbf{c}_k^{.2}\overset{2}{\cdot}\nabla^2\partial_tf_{\overline{k}}^{eq}+O(\Delta t^4),
\end{align}
and can be simplified as 
\begin{align}\label{fi-star-2}
	f_{\overline{k}}^{\star}(\mathbf{x}_f,t)= &f_{\overline{k}}^{eq}+\frac{\Delta t}{2}\partial_tf_{\overline{k}}^{eq}-\frac{\Delta t^2}{8}\mathbf{c}_k^{.2}\overset{2}{\cdot}\nabla^2 f_{\overline{k}}^{eq} \notag\\
	&+\frac{\Delta t^3}{8}\mathbf{c}_k\cdot\nabla \partial_t^2f_{\overline{k}}^{eq}-\frac{\Delta t^3}{16}\mathbf{c}_k^{.2}\overset{2}{\cdot}\nabla^2\partial_tf_{\overline{k}}^{eq}-\frac{\Delta t^3}{24}\partial_t^3f_{\overline{k}}^{eq} +\frac{\Delta t}{2}R_k +O(\Delta t^4).
\end{align} 
Finally,  the first- to fourth-order half-way boundary schemes (\ref{bun2}) of the distribution functions $f_k(\mathbf{x}_f,t)$ for the Dirichlet boundary condition can be obtained after summing Eqs. (\ref{fi-2}) and   (\ref{fi-star-2}). 
\section{The derivation of the inverse of the transform matrix $\mathbf{M}$}\label{inverse-M}
For the unified form of the transform matrix $\mathbf{M}$ in Eq. (\ref{DdQq-M}), we now derive the unified form of the inverse of  $\mathbf{M}$ and denote it as $\mathbf{\tilde{M}}$.

According to the rotational symmetry of the Cartesian coordinates, without loss of generality, here we only consider the $1_{st}$ and $2_{nd}$ spatial directions for the $d$-dimensional case. To this end, let us take the derivation of the $\big[(3d^2+d+4)/2\big]_{th}$ column of the matrix $\mathbf{\tilde{M}}$, i.e., $\mathbf{\tilde{M}}^{\dot{m}^{eq|4}_{x_1^2x_2^2}}_{(3d^2+d+4)/2}$,  as an example. In particular, according to $\mathbf{M}\mathbf{\tilde{M}}=\mathbf{I}$, we have
\begin{subequations}
	\begin{align}
		&I_{2d^2+1}^T\mathbf{\tilde{M}}^{\dot{m}^{eq|4}_{x_1^2x_2^2}}_{(3d^2+d+4)/2}=0,\label{x0y0}\\
		&(\mathbf{c}^{x_i})^T\mathbf{\tilde{M}}^{\dot{m}^{eq|4}_{x_1^2x_2^2}}_{(3d^2+d+4)/2}=0,\:\: i\in  \llbracket{1,d}\rrbracket,\label{x1}\\
		&\big[(\mathbf{c}^{x_i})^{.2}\big]^T\mathbf{\tilde{M}}^{\dot{m}^{eq|4}_{x_1^2x_2^2}}_{(3d^2+d+4)/2}=0,\:\:  i\in\llbracket{1,d}\rrbracket,\label{x2}\\
		&(\mathbf{c}^{x_i}.\mathbf{c}^{x_j})^T\mathbf{\tilde{M}}^{\dot{m}^{eq|4}_{x_1^2x_2^2}}_{(3d^2+d+4)/2}=0,\:\:  i,j\in\llbracket{1,d}\rrbracket,\: (i< j),\label{x1y1}\\
		&\big[(\mathbf{c}^{x_i})^{.2}.\mathbf{c}^{x_j}\big]^T\mathbf{\tilde{M}}^{\dot{m}^{eq|4}_{x_1^2x_2^2}}_{(3d^2+d+4)/2}=0,\:\:  i,j\in\llbracket{1,d}\rrbracket,\: (i\neq j),\label{x2y1}\\
		&\big[(\mathbf{c}^{x_1})^{.2}.(\mathbf{c}^{x_2})^{.2}\big]^T\mathbf{\tilde{M}}^{\dot{m}^{eq|4}_{x_1^2x_2^2}}_{(3d^2+d+4)/2}=1,\label{x2y2-12}\\
		&\big[(\mathbf{c}^{x_i})^{.2}.(\mathbf{c}^{x_j})^{.2}\big]^T\mathbf{\tilde{M}}^{\dot{m}^{eq|4}_{x_1^2x_2^2}}_{(3d^2+d+4)/2}=0,\:\:  i,j\in\llbracket{1,d}\rrbracket,\: (i< j,\:i\neq 1,\: j\neq 2),\label{x2y2}
	\end{align}
\end{subequations}
where $\mathbf{\tilde{M}}^{\dot{m}^{eq|4}_{x_1^2x_2^2}}_{(3d^2+d+4)/2}$ can be written as
\begin{align}
	\mathbf{\tilde{M}}^{\dot{m}^{eq|4}_{x_1^2x_2^2}}_{(3d^2+d+4)/2}=\Bigg[&\tilde{M}^{\dot{m}^{eq|4}_{x_1^2x_2^2}}_{0}
	,\tilde{M}^{\dot{m}^{eq|4}_{x_1^2x_2^2}}_{1},\tilde{M}^{\dot{m}^{eq|4}_{x_1^2x_2^2}}_{2},\ldots,\tilde{M}^{\dot{m}^{eq|4}_{x_1^2x_2^2}}_{d},\tilde{M}^{\dot{m}^{eq|4}_{x_1^2x_2^2}}_{\overline{1}},\tilde{M}^{\dot{m}^{eq|4}_{x_1^2x_2^2}}_{\overline{2}},\ldots,\tilde{M}^{\dot{m}^{eq|4}_{x_1^2x_2^2}}_{\overline{d}},\notag\\
	&\tilde{M}^{\dot{m}^{eq|4}_{x_1^2x_2^2}}_{1,2},\tilde{M}^{\dot{m}^{eq|4}_{x_1^2x_2^2}}_{\overline{1},2}\tilde{M}^{\dot{m}^{eq|4}_{x_1^2x_2^2}}_{\overline{1},\overline{2}},\tilde{M}^{\dot{m}^{eq|4}_{x_1^2x_2^2}}_{1,\overline{2}}, \tilde{M}^{\dot{m}^{eq|4}_{x_1^2x_2^2}}_{1,3},\tilde{M}^{\dot{m}^{eq|4}_{x_1^2x_2^2}}_{\overline{1},3}\tilde{M}^{\dot{m}^{eq|4}_{x_1^2x_2^2}}_{\overline{1},\overline{3}},\tilde{M}^{\dot{m}^{eq|4}_{x_1^2x_2^2}}_{1,\overline{3}}, \ldots, \tilde{M}^{\dot{m}^{eq|4}_{x_1^2x_2^2}}_{1,d},\tilde{M}^{\dot{m}^{eq|4}_{x_1^2x_d^2}}_{\overline{1},d}\tilde{M}^{\dot{m}^{eq|4}_{x_1^2x_d^2}}_{\overline{1},\overline{d}},\tilde{M}^{\dot{m}^{eq|4}_{x_1^2x_2^2}}_{1,\overline{d}},\notag\\
	&\tilde{M}^{\dot{m}^{eq|4}_{x_1^2x_2^2}}_{2,3},\tilde{M}^{\dot{m}^{eq|4}_{x_1^2x_2^2}}_{\overline{2},3}\tilde{M}^{\dot{m}^{eq|4}_{x_1^2x_2^2}}_{\overline{2},\overline{3}},\tilde{M}^{\dot{m}^{eq|4}_{x_1^2x_2^2}}_{2,\overline{3}}, \tilde{M}^{\dot{m}^{eq|4}_{x_1^2x_2^2}}_{2,4},\tilde{M}^{\dot{m}^{eq|4}_{x_1^2x_2^2}}_{\overline{2},4}\tilde{M}^{\dot{m}^{eq|4}_{x_1^2x_2^2}}_{\overline{2},\overline{4}},\tilde{M}^{\dot{m}^{eqv}_{x_1^2x_2^2}}_{2,\overline{4}},\ldots,  \tilde{M}^{\dot{m}^{eq|4}_{x_1^2x_2^2}}_{2,d},\tilde{M}^{\dot{m}^{eq|4}_{x_1^2x_2^2}}_{\overline{2},d}\tilde{M}^{\dot{m}^{eq|4}_{x_1^2x_2^2}}_{\overline{2},\overline{d}},\tilde{M}^{\dot{m}^{eq|4}_{x_1^2x_2^2}}_{2,\overline{d}},\notag\\
	&\ldots,\\
	&\tilde{M}^{\dot{m}^{eq|4}_{x_1^2x_2^2}}_{d-1,d},\tilde{M}^{\dot{m}^{eq|4}_{x_1^2x_2^2}}_{\overline{d-1},d}\tilde{M}^{\dot{m}^{eq|4}_{x_1^2x_2^2}}_{\overline{d-1},\overline{d}},\tilde{M}^{\dot{m}^{eq|4}_{x_1^2x_2^2}}_{d-1,\overline{d}},
	\Bigg]^T.
\end{align}

From Eq. (\ref{x2y1}), we can first obtain
\begin{align}\label{con-1}
	\begin{cases}
		\tilde{M}^{\dot{m}^{eq|4}_{x_1^2x_2^2}}_{1,2}=\tilde{M}^{\dot{m}^{eq|4}_{x_1^2x_2^2}}_{\overline{1},\overline{2}},\:\: \tilde{M}^{\dot{m}^{eq|4}_{x_1^2x_2^2}}_{\overline{1},2}=\tilde{M}^{\dot{m}^{eq|4}_{x_1^2x_2^2}}_{1,\overline{2}},&\\
		\tilde{M}^{\dot{m}^{eq|4}_{x_1^2x_2^2}}_{1,3}=\tilde{M}^{\dot{m}^{eq|4}_{x_1^2x_2^2}}_{\overline{1},\overline{3}},\:\: \tilde{M}^{\dot{m}^{eq|4}_{x_1^2x_2^2}}_{\overline{1},3}=\tilde{M}^{\dot{m}^{eq|4}_{x_1^2x_2^2}}_{1,\overline{3}},&\\
	\qquad\qquad\qquad	\ldots\qquad\qquad\qquad\:\:,&\\
		\tilde{M}^{\dot{m}^{eq|4}_{x_1^2x_2^2}}_{d-1,d}=\tilde{M}^{\dot{m}^{eq|4}_{x_1^2x_2^2}}_{\overline{d-1},\overline{d}},\:\: \tilde{M}^{\dot{m}^{eq|4}_{x_1^2x_2^2}}_{\overline{d-1},d}=\tilde{M}^{\dot{m}^{eq|4}_{x_1^2x_2^2}}_{d-1,\overline{d}},
	\end{cases}
\end{align}
  from Eqs. (\ref{x2y2-12}) and (\ref{x2y2}), and in combination with the relations shown in Eq. (\ref{con-1}), we have
\begin{align}
	\begin{cases}
		\tilde{M}^{\dot{m}^{eq|4}_{x_1^2x_2^2}}_{1,2}+\tilde{M}^{\dot{m}^{eq|4}_{x_1^2x_2^2}}_{\overline{1},3}=1/(2c^4),\:\: \tilde{M}^{\dot{m}^{eq|4}_{x_1^2x_2^2}}_{1,2}=\tilde{M}^{\dot{m}^{eq|4}_{x_1^2x_2^2}}_{\overline{1},\overline{2}},\:\: \tilde{M}^{\dot{m}^{eq|4}_{x_1^2x_2^2}}_{\overline{1},2}=\tilde{M}^{\dot{m}^{eq|4}_{x_1^2x_2^2}}_{1,\overline{2}},&\\
		\tilde{M}^{\dot{m}^{eq|4}_{x_1^2x_2^2}}_{1,3}+\tilde{M}^{\dot{m}^{eq|4}_{x_1^2x_2^2}}_{\overline{1},3}=0,\:\: \tilde{M}^{\dot{m}^{eq|4}_{x_1^2x_2^2}}_{1,3}=\tilde{M}^{\dot{m}^{eq|4}_{x_1^2x_2^2}}_{\overline{1},\overline{3}},\:\: \tilde{M}^{\dot{m}^{eq|4}_{x_1^2x_2^2}}_{\overline{1},3}=\tilde{M}^{\dot{m}^{eq|4}_{x_1^2x_2^2}}_{1,\overline{3}},&\\
		\qquad\qquad\qquad\qquad\qquad	\ldots\qquad\qquad\qquad\qquad\qquad\:\:,&\\
		\tilde{M}^{\dot{m}^{eq|4}_{x_1^2x_2^2}}_{d-1,d}+\tilde{M}^{\dot{m}^{eq|4}_{x_1^2x_2^2}}_{\overline{d-1},d}=0,\:\: \tilde{M}^{\dot{m}^{eq|4}_{x_1^2x_2^2}}_{d-1,d}=\tilde{M}^{\dot{m}^{eq|4}_{x_1^2x_2^2}}_{\overline{d-1},\overline{d}},\:\: \tilde{M}^{\dot{m}^{eq|4}_{x_1^2x_2^2}}_{\overline{d-1},d}=\tilde{M}^{\dot{m}^{eq|4}_{x_1^2x_2^2}}_{d-1,\overline{d}},
	\end{cases}
\end{align}
then, from Eq. (\ref{x1y1}), we can further obtain
\begin{subequations}\label{4rd}
	\begin{align}
	&	\tilde{M}^{\dot{m}^{eq|4}_{x_1^2x_2^2}}_{1,2}=\tilde{M}^{\dot{m}^{eq|4}_{x_1^2x_2^2}}_{\overline{1},2}=\tilde{M}^{\dot{m}^{eq|4}_{x_1^2x_2^2}}_{\overline{1},\overline{2}}=\tilde{M}^{\dot{m}^{eq|4}_{x_1^2x_2^2}}_{1,\overline{2}}=1/(4c^4),\\
	& \tilde{M}^{\dot{m}^{eq|4}_{x_1^2x_2^2}}_{1,3}=\tilde{M}^{\dot{m}^{eq|4}_{x_1^2x_2^2}}_{\overline{1},\overline{3}}= \tilde{M}^{\dot{m}^{eq|4}_{x_1^2x_2^2}}_{\overline{1},3}=\tilde{M}^{\dot{m}^{eq|4}_{x_1^2x_2^2}}_{1,\overline{3}},\\
	&	\qquad\qquad\qquad	\ldots\qquad\qquad\qquad\:\:,\\
	&	\tilde{M}^{\dot{m}^{eq|4}_{x_1^2x_2^2}}_{d-1,d}=\tilde{M}^{\dot{m}^{eq|4}_{x_1^2x_2^2}}_{\overline{d-1},\overline{d}}=\tilde{M}^{\dot{m}^{eq|4}_{x_1^2x_2^2}}_{\overline{d-1},d}=\tilde{M}^{\dot{m}^{eq|4}_{x_1^2x_2^2}}_{d-1,\overline{d}}=0.
	\end{align}
\end{subequations}

Similar to the derivation of the terms $\tilde{M}^{\dot{m}^{eq|4}_{x_1^2x_2^2}}_{i,j}$, $\tilde{M}^{\dot{m}^{eq|4}_{x_1^2x_2^2}}_{\overline{i},j}$, $\tilde{M}^{\dot{m}^{eq|4}_{x_1^2x_2^2}}_{\overline{i},\overline{j}}$ and $\tilde{M}^{\dot{m}^{eq|4}_{x_1^2x_2^2}}_{i,\overline{j}}$ ($i,j\in\llbracket{1,d}\rrbracket,\: j>i$) above,  according to Eqs. (\ref{x1}), (\ref{x2}), and  (\ref{x1y1}), the terms $\tilde{M}^{\dot{m}^{eq|4}_{x_1^2x_2^2}}_{i}$ and $\tilde{M}^{\dot{m}^{eq|4}_{x_1^2x_2^2}}_{\overline{i}}$ ($i\in \llbracket{1,d}\rrbracket$) can be determined as
\begin{subequations}\label{1rd}
	\begin{align} 
			&\tilde{M}^{\dot{m}^{eq|4}_{x_1^2x_2^2}}_{1}=\tilde{M}^{m^{x_1^2x_2^2}_4}_{\overline{1}}=-1/(2c^4),\\
			&\tilde{M}^{\dot{m}^{eq|4}_{x_1^2x_2^2}}_{2}=\tilde{M}^{m^{x_1^2x_2^2}_4}_{\overline{2}}=-1/(2c^4),\\
			&\tilde{M}^{\dot{m}^{eq|4}_{x_1^2x_2^2}}_{3}=\tilde{M}^{m^{x_1^2x_2^2}_4}_{\overline{3}}=0,\\
			&\quad\qquad	\ldots\quad\qquad\:\:,\\
			&\tilde{M}^{\dot{m}^{eq|4}_{x_1^2x_2^2}}_{d}=\tilde{M}^{m^{x_1^2x_2^2}_4}_{\overline{d}}=0. 
	\end{align}
\end{subequations}
Finally, according to Eq. (\ref{x0y0}), and with the aid of Eqs. (\ref{4rd}) and (\ref{1rd}), we can deduce that $\tilde{M}^{\dot{m}^{eq|4}_{x_1^2x_2^2}}_{0}=1/c^4$. The other columns of the matrix $\mathbf{\tilde{M}}$ can be similarly obtained, and we do not present the detailed derivation here.

\section{The analysis on Eqs. (\ref{eq-2}) and (\ref{eq-4}) when $\sigma_{1,2}\in C_2\cup C_3\cup C_4$}\label{Appen-entropy-2D}
Here we will prove that Eqs. (\ref{eq-2}) and (\ref{eq-4}) do not hold when 
$\sigma_{1,2}\in C_2\cup C_3\cup C_4$, and the proof can be divided into the following three parts.
\subsection{The proof of $\sigma_{1,2}\notin C_2$} \label{app-case-1}
Let $\sigma_1=\pm1/2\pm\delta_1$ and $\sigma_2=\pm1/2\pm\delta_2$ with $1/2>\delta_1\geq \delta_2>0$, we have
\begin{align}
	2\sigma_1^2-2|\sigma_1|+1-|\sigma_1-\sigma_2|-2\sigma_1\sigma_2=&2\delta_1^2-\delta_1+\delta_2-2\delta_1\delta_2-\delta_1-\delta_2\notag\\
	=&2\delta_1(\delta_1-1-\delta_2)<0,
\end{align}
which contradicts with Eq. (\ref{eq-2}).

Similarly, let $\sigma_1=\pm1/2\pm\delta_1$ and $\sigma_2=\pm1/2\pm\delta_2$ with $1/2>\delta_2\geq \delta_1>0$, we have
\begin{align}
	2\sigma_2^2-2|\sigma_2|+1-|\sigma_1-\sigma_2|-2\sigma_1\sigma_2=&2\delta_2^2-\delta_2+\delta_1-2\delta_1\delta_2-\delta_1-\delta_2\notag\\
	=&2\delta_2(\delta_2-1-\delta_1)<0,
\end{align}
which contradicts with Eq. (\ref{eq-4}). Thus $\sigma_{1,2}\notin C_2$.
\subsection{The proof $\sigma_{1,2}\notin C_3$}\label{app-case-2}
For the case of $\sigma_2=1-\sigma_1$ with $0<\sigma_1<1/2$, we have
\begin{align}
	2\sigma_1^2-2|\sigma_1|+1-|\sigma_1-\sigma_2|-2\sigma_1\sigma_2=2\sigma_1(2\sigma_1-1)<0,
	\end{align}
	which contradicts with Eq. (\ref{eq-2}). 
	
For the case of $\sigma_2=1-\sigma_1$ with $1/2<\sigma_1<1$, we have
	\begin{align}
		2\sigma_1^2-2|\sigma_1|+1-|\sigma_1-\sigma_2|-2\sigma_1\sigma_2=4\sigma_1^2-6\sigma_1+2<0,
	\end{align}
	which contradicts with Eq. (\ref{eq-2}).  Thus $\sigma_{1,2}\notin C_3$.
\subsection{The proof $\sigma_{1,2}\notin C_4$}\label{app-case-3}
For the case of $\sigma_1=-1-\sigma_2$ with $-1/2<\sigma_2<0$, we have
 \begin{align}
 	2\sigma_2^2-2|\sigma_2|+1-|\sigma_1-\sigma_2|-2\sigma_1\sigma_2=2\sigma_2(2\sigma_2-1)<0,
 \end{align}
 which contradicts with Eq. (\ref{eq-4}). 
 
For the case of $\sigma_1=-1-\sigma_2$ with $-1<\sigma_2<-1/2$, we have
 \begin{align}
 	2\sigma_2^2-2|\sigma_2|+1-|\sigma_1-\sigma_2|-2\sigma_1\sigma_2=4\sigma_2^2-6\sigma_2+2<0,
 \end{align}
 which contradicts with Eq. (\ref{eq-4}). Thus $\sigma_{1,2}\notin C_4$.
 
 In combination with the conclusions in  \ref{app-case-1}, \ref{app-case-2}, and \ref{app-case-3}, the proof of $\sigma_{1,2}\notin C1\cup C_2\cup C_3 \cup C_4$ is completed.
 \section{The proof that the BGK-LB model (\ref{BGK-LB-Evolution-seq2}) is $L^2$ stable for the case of $|\sigma_1|=1/2$}\label{app-sigma=1/2}
 For simplicity, here we only consider the case of $\sigma_1=-1/2$, and the analysis for the other case where $\sigma_1=1/2$ is similar. In this case, the amplification matrix $\mathbf{G}$  in Eq. (\ref{G}) is similar to the following matrix $\mathbf{\tilde{G}}$
 \begin{align}
 	\mathbf{\tilde{G}}=\left(\begin{matrix}
 		0&1&1\\
 		e^{-\xi}&0&e^{-\xi}\\
 		0&0&-e^{\xi}
 	\end{matrix}\right).
 \end{align}
 When $\xi=0$, the discussion is identical to that in Sec. \ref{l2-stability}. When $\xi\neq 0$, it is easy to check that the three distinct roots of the characteristic polynomial of the  matrix $\mathbf{\tilde{G}}$ are $-e^{\xi}$, $e^{-\xi/2}$, and $-e^{-\xi/2}$, respectively. Thus, for the case of $|\sigma_1|=1/2$, the minimal polynomial of the amplification matrix $\mathbf{G}$  in Eq. (\ref{G}) is still a simple von Neumann polynomial, i.e., the BGK-LB model (\ref{BGK-LB-Evolution-seq2}) is $L^2$ stable.
  \bibliographystyle{elsarticle-num-names} 
 \bibliography{reference}
\end{document}